\documentclass[12pt]{amsart}
\usepackage{amssymb}
\usepackage[all]{xy}

\usepackage{verbatim}
\usepackage{eucal}
\usepackage{latexsym}
\usepackage{amsfonts,amssymb,amsxtra}
\usepackage{relsize}
\usepackage[dvips]{color}
\usepackage{latexsym}
\usepackage{mathptmx}
\usepackage{amssymb}
\usepackage{graphicx}
\usepackage{amsthm}
\usepackage{amsfonts}
\usepackage{amscd}
\usepackage{bm}
\usepackage[utf8]{inputenc}

\usepackage{verbatim} 
\usepackage{eucal} 
\usepackage{color}

\usepackage[draft]{todonotes}
\usepackage{pgf}
\usepackage{tikz}
\usetikzlibrary{arrows,automata}
\usepackage[
ordering=Kac,
edge/.style=black,
indefinite-edge={draw=green,fill=white,densely dashed},
indefinite-edge-ratio=5,
mark=o,
root-radius=.06cm]
{dynkin-diagrams}

\usepackage{pgf, tikz}

\usetikzlibrary{arrows, automata}
\usepackage{xcolor}
\usetikzlibrary{arrows}
\usetikzlibrary {positioning}

\definecolor {processblue}{cmyk}{0.96,0,0,0}
\input xy
\xyoption{all}
\usepackage{mathtools}
\usepackage{amsfonts}
\usepackage{amsthm}
\usepackage{mathrsfs}

\usepackage{amscd}
 \numberwithin{equation}{section}
  \newtheorem{The}[equation]{Theorem}
  \newtheorem{Pro}[equation]{Proposition}
  \newtheorem{Lem}[equation]{Lemma}
  \newtheorem{Cor}[equation]{Corollary}

  \newtheorem{Examp}[equation]{Example}
  
 \numberwithin{equation}{section}

\newcommand{\bsm}{\begin{smallmatrix}}
\newcommand{\esm}{\end{smallmatrix}}
\newcommand{\bbm}{\begin{matrix}}
\newcommand{\ebm}{\end{matrix}}

\begin{document}
\title{ Coxeter diagrams and the K\"othe's problem}

\author{Ziba Fazelpour}
\address{ School of Mathematics, Institute for Research in Fundamental Sciences (IPM), P.O. Box: 19395-5746, Tehran, Iran}
\email{z.fazelpour@ipm.ir}
\author{Alireza Nasr-Isfahani}
\address{Department of Mathematics\\
University of Isfahan\\
P.O. Box: 81746-73441, Isfahan, Iran and School of Mathematics, Institute for Research in Fundamental Sciences (IPM), P.O. Box: 19395-5746, Tehran, Iran}
\email{nasr$_{-}$a@sci.ui.ac.ir / nasr@ipm.ir}



\begin{abstract}
A ring $\Lambda$ is called right K\"othe if every right $\Lambda$-module is a direct sum of cyclic modules. In this paper, we give a characterization of basic hereditary right K\"othe rings in terms of their Coxeter valued quivers. Also we give a characterization of basic right K\"othe rings with radical square zero. Therefore we give a solution of the K\"othe's problem in these two cases.
\end{abstract}

\maketitle
\section{Introduction}
It is known that every finitely generated $\mathbb{Z}$-module is a direct sum of cyclic modules. The idea of this important property of abelian groups go back to Pr\"{u}fer \cite{Pu}. K\"othe showed that artinian principal ideal rings have this property. He also proved that if a commutative artinian ring $\Lambda$ has the property that each of its $\Lambda$-modules is a direct sum of cyclic modules, then it is a principal ideal ring and posed the question to classify the noncommutative rings with this property \cite{K}. The K\"othe's problem is one of the old problems in rings and modules theory which is not solved yet. A ring for which any right module is a direct sum of cyclic modules, is now called a {\it right K\"othe} ring. Nakayama gave an example of a right K\"{o}the ring $\Lambda$ which is not a principal right ideal ring (see \cite[page 289]{Na}). Later, Cohen and Kaplansky proved that if a commutative ring $\Lambda$ is K\"othe, then $\Lambda$ is an artinian principal ideal ring \cite{CK}. Combining results of Cohen and Kaplansky \cite{CK} and K\"othe \cite{K} one obtains that, a commutative ring $\Lambda$ is K\"othe if and only if $\Lambda$ is an artinian principal ideal ring. A right artinian ring $\Lambda$ is called {\it representation-finite} provided $\Lambda$ has, up to isomorphism, only finitely many finitely generated indecomposable right $\Lambda$-modules. Following \cite{Sim6}, we call the ring $\Lambda$ {\it right pure semisimple if every right $\Lambda$-module is a direct sum of finitely generated right $\Lambda$-modules.} It is known that a commutative ring $\Lambda$ is pure semisimple, if and only if, $\Lambda$ is a representation-finite ring, if and only if, $\Lambda$ is a K\"othe ring \cite{Grif}.  A ring $\Lambda$ is a representation-finite ring if and only if $\Lambda$ is right pure semisimple and left pure semisimple \cite{Auslander}. The problem of whether right pure semisimple rings are representation-finite, known as the pure semisimplicity conjecture, still remains open (see \cite{Auslander}, \cite{Sim6} and \cite{Sim7}) . It seems that there is a strong connection between pure semisimplicity conjecture and K\"othe's problem.

Kawada completely solved the K\"othe's problem for the basic finite-dimensional $K$-algebras \cite{Kawada1, Kawada2, Kawada3} (see also \cite{Ringel2}). Kawada's papers contain a set of 19 conditions which characterize Kawada algebras, as well as, the list of all possible finitely generated indecomposable modules. Ringel by using the multiplicity-free of top and soc of finitely generated indecomposable  modules, gave a characterization of Kawada algebras \cite{Ringel2}. Behboodi et al. proved that if $\Lambda$ is a right K\"othe ring in which all idempotents are central, then $\Lambda$ is an artinian principal left ideal ring \cite{BGMS}. Recently the authors have studied the K\"othe's problem \cite{FN}. In fact, all known results related to the characterization of right K\"othe rings follow from Corollary 3.2 in \cite{FN}.

In representation theory, representation-finite algebras are of particular importance since in this case one has a
complete combinatorial description of the module category in terms of the Auslander-Reiten quiver. By \cite[Theorem 4.4]{chase}, any right K\"othe ring $\Lambda$ is right artinian, then there is a finite upper bound for the lengths of the finitely
generated indecomposable right $\Lambda$-modules. Thus by \cite[Proposition 54.3]{wi}, any right K\"othe ring is an artinian representation-finite ring. It seems that a solution of the K\"othe's problem need a classification of all representation-finite rings and some further information about the structure of the modules over representation-finite rings, which is a rather difficult problem. In this paper, by using the representation theory techniques and classifications of representation-finite hereditary rings \cite{DR1, DRS, DS1}, we solve the K\"othe's problem in this case. As a consequence we solve the K\"othe's problem for the class of rings with radical square zero.

We recall that a unitary ring $\Lambda$ is defined to be {\it right hereditary} if every right ideal of $\Lambda$ is projective. In \cite{4}, Gabriel proved that a hereditary finite-dimensional $K$-algebra $\Lambda$, over an algebraically closed field $K$, is representation-finite if and only if the underlying graph of its quiver $Q_{\Lambda}$ (see \cite{assem}) is a disjoint union of the Dynkin diagrams $\mathbb{A}_n$, $\mathbb{D}_n$, $\mathbb{E}_6$, $\mathbb{E}_7$ and $\mathbb{E}_8$ presented in Table A, that appear also in Lie theory. It is also shown in \cite{4} that there is a bijection between the isomorphism classes of finite-dimensional indecomposable representations and the positive roots of the corresponding Dynkin diagrams, see also \cite[Ch. VII]{assem}.
\begin{center}
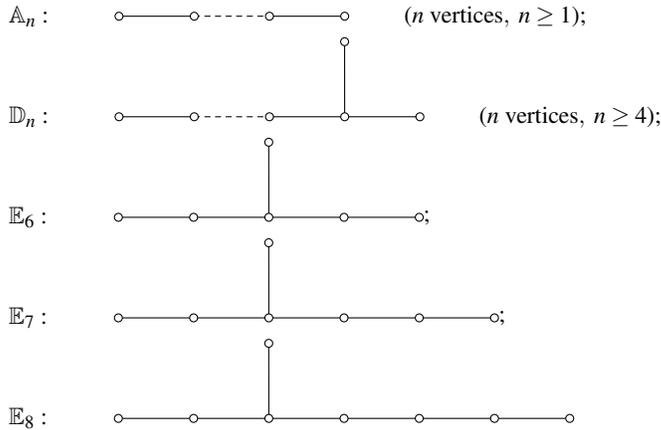
\textbf{Table A.} Dynkin diagrams
\end{center}

	\definecolor{qqqqff}{rgb}{0.,0.,1.}
	\begin{tikzpicture}[line cap=round,line join=round,>=triangle 45,x=1.0cm,y=1.0cm]
	\draw (2.,9.)-- (3.,9.);
	\draw [dash pattern=on 2pt off 2pt](3.,9.)--(4.,9.);
	\draw (4.,9.)-- (5.,9.);
		\begin{scriptsize}
	\draw[color=black] (.8,9.) node {$\mathbb{A}_n:$};
	\draw [fill=white] (2.,9.) circle (1.5pt);
	\draw [fill=white] (3.,9.) circle (1.5pt);
	\draw [fill=white] (4.,9.) circle (1.5pt);
	\draw [fill=white] (5.,9.) circle (1.5pt);
	\draw[color=black] (7.,9.) node {($n ~~{\rm vertices},~ n\geq 1$);};
	\end{scriptsize}
	\end{tikzpicture}

	\definecolor{qqqqff}{rgb}{0.,0.,1.}
	\begin{tikzpicture}[line cap=round,line join=round,>=triangle 45,x=1.0cm,y=1.0cm]
	\draw (2.,9.)-- (3.,9.);
	\draw [dash pattern=on 2pt off 2pt](3.,9.)--(4.,9.);
	\draw (4.,9.)-- (5.,9.);
	\draw (5.,9.)-- (6.,9.);
	\draw (5.,9.)-- (6.,10.9999999999987,2.6);
		\begin{scriptsize}
	\draw[color=black] (.8,9.) node {$\mathbb{D}_n:$};
	\draw [fill=white] (2.,9.) circle (1.5pt);
	\draw [fill=white] (3.,9.) circle (1.5pt);
	\draw [fill=white] (4.,9.) circle (1.5pt);
	\draw [fill=white] (5.,9.) circle (1.5pt);
	\draw [fill=white] (6.,9.) circle (1.5pt);
	\draw [fill=white] (6.,10.9999999999987,2.6) circle (1.5pt);
	\draw[color=black] (8.,9.) node {($n ~~{\rm vertices}, ~n\geq 4$);};
	
	\end{scriptsize}
	\end{tikzpicture}

	\definecolor{qqqqff}{rgb}{0.,0.,1.}
	\begin{tikzpicture}[line cap=round,line join=round,>=triangle 45,x=1.0cm,y=1.0cm]
	\draw (2.,9.)-- (3.,9.);
	\draw (3.,9.)-- (4.,9.);
	\draw (4.,9.)-- (5.,9.);
	\draw (4.,9.)-- (5.,10.9999999999987,2.6);
	\draw (5.,9.)-- (6.,9.);
		\begin{scriptsize}
	\draw[color=black] (.8,9.) node {$\mathbb{E}_6:$};
	\draw [fill=white] (2.,9.) circle (1.5pt);
	\draw [fill=white] (3.,9.) circle (1.5pt);
	\draw [fill=white] (4.,9.) circle (1.5pt);
	\draw [fill=white] (5.,9.) circle (1.5pt);
	\draw [fill=white] (5.,10.9999999999987,2.6) circle (1.5pt);
	\draw [fill=white] (6.,9.) circle (1.5pt);
	\draw[color=black] (6.1,9.) node {;};
	\end{scriptsize}
	\end{tikzpicture}

	\definecolor{qqqqff}{rgb}{0.,0.,1.}
	\begin{tikzpicture}[line cap=round,line join=round,>=triangle 45,x=1.0cm,y=1.0cm]
	\draw (2.,9.)-- (3.,9.);
	\draw (3.,9.)-- (4.,9.);
	\draw (4.,9.)-- (5.,9.);
	\draw (4.,9.)-- (5.,10.9999999999987,2.6);
	\draw (5.,9.)-- (6.,9.);
	\draw (6.,9.)-- (7.,9.);
		\begin{scriptsize}
	\draw[color=black] (.8,9.) node {$\mathbb{E}_7:$};
	\draw [fill=white] (2.,9.) circle (1.5pt);
	\draw [fill=white] (3.,9.) circle (1.5pt);
	\draw [fill=white] (4.,9.) circle (1.5pt);
	\draw [fill=white] (5.,9.) circle (1.5pt);
	\draw [fill=white] (5.,10.9999999999987,2.6) circle (1.5pt);
	\draw [fill=white] (6.,9.) circle (1.5pt);
	\draw [fill=white] (7.,9.) circle (1.5pt);
	\draw[color=black] (7.1,9.) node {;};
	\end{scriptsize}
	\end{tikzpicture}

	\definecolor{qqqqff}{rgb}{0.,0.,1.}
	\begin{tikzpicture}[line cap=round,line join=round,>=triangle 45,x=1.0cm,y=1.0cm]
	\draw (2.,9.)-- (3.,9.);
	\draw (3.,9.)-- (4.,9.);
	\draw (4.,9.)-- (5.,9.);
	\draw (4.,9.)-- (5.,10.9999999999987,2.6);
	\draw (5.,9.)-- (6.,9.);
	\draw (6.,9.)-- (7.,9.);
	\draw (7.,9.)-- (8.,9.);
		\begin{scriptsize}
	\draw[color=black] (.8,9.) node {$\mathbb{E}_8:$};
	\draw [fill=white] (2.,9.) circle (1.5pt);
	\draw [fill=white] (3.,9.) circle (1.5pt);
	\draw [fill=white] (4.,9.) circle (1.5pt);
	\draw [fill=white] (5.,9.) circle (1.5pt);
	\draw [fill=white] (5.,10.9999999999987,2.6) circle (1.5pt);
	\draw [fill=white] (6.,9.) circle (1.5pt);
	\draw [fill=white] (7.,9.) circle (1.5pt);
	\draw [fill=white] (8.,9.) circle (1.5pt);
	\end{scriptsize}
	\end{tikzpicture}
	
We recall from \cite{3} that a {\it species} $\mathscr{M}= \left( F_i,{_iM_j} \right)_{i, j \in I} $ is a finite set of division rings $F_i$ and $F_i$-$F_j$-bimodules $_iM_j$, $i \neq j$. To an arbitrary basic ring $\Lambda$ one attaches its species ${\mathscr{M}}_{\Lambda}$ as follows. Let $\Lambda/J \cong \bigoplus_{i=1}^nF_i$, where $n \in {\Bbb{N}}$, each $F_i$ is a division ring and $J$ is the Jacobson radical of $\Lambda$. We can write $J/ J^2 = \bigoplus_{1 \leq i, j \leq n}{_iM_j}$, where each $_iM_j= F_i(J/J^2)F_j$ is an $F_i$-$F_j$-bimodule. ${\mathscr{M}}_{\Lambda}= \left( F_i,{_iM_j} \right)_{i, j \in I} $ is called {\it the species of} $\Lambda$. Let $\Lambda$ be a basic hereditary ring and let ${\mathscr{M}}_{\Lambda}= \left( F_i,{_iM_j} \right)_{i, j \in I} $ be the species of $\Lambda$.  We recall from \cite{simson} that a {\it Coxeter valued quiver} $(\mathscr{C}_{\Lambda}, \mathbf{m})$ of $\Lambda$ is a quiver with vertices $1, 2, \ldots, n$ corresponding to the division rings $F_1, F_2, \ldots, F_n$. There exists a valued arrow
\begin{center}
	\definecolor{qqqqff}{rgb}{0.,0.,1.}
	\begin{tikzpicture}[line cap=round,line join=round,>=triangle 45,x=1.0cm,y=1.0cm]
	\draw (2.,9.)-- (3.,9.);
		\begin{scriptsize}
	\draw[color=black] (.8,9.) node {};
	\draw [fill=white] (2.,9.) circle (1.5pt);
	\draw[color=black] (2.,8.7) node {$i$};
	\draw[color=black] (2.5,9.2) node {$m_{ij}$};
	\draw[color=black] (2.85,9.) node {$>$};
	\draw [fill=white] (3.,9.) circle (1.5pt);
	\draw[color=black] (3.,8.7) node {$j$};
	
	\end{scriptsize}
	\end{tikzpicture}
\end{center}
 in $(\mathscr{C}_{\Lambda}, \mathbf{m})$ if and only if the $F_i$-$F_j$-bimodule ${_iM_j}$ is non-zero and there exists exactly $m_{ij} \geq 3$ pairwise non-isomorphic indecomposable right $\tiny{\begin{pmatrix}
F_i & {_iM_j}\\
0& F_j
\end{pmatrix}}$-modules. If $m_{ij} = 3$, we write simply \begin{center}
	\definecolor{qqqqff}{rgb}{0.,0.,1.}
	\begin{tikzpicture}[line cap=round,line join=round,>=triangle 45,x=1.0cm,y=1.0cm]
	\draw (2.,9.)-- (3.,9.);
		\begin{scriptsize}
	\draw[color=black] (.8,9.) node {};
	\draw [fill=white] (2.,9.) circle (1.5pt);
	\draw[color=black] (2.,8.7) node {$i$};
	\draw[color=black] (2.85,9.) node {$>$};
	\draw [fill=white] (3.,9.) circle (1.5pt);
	\draw[color=black] (3.,8.7) node {$j$};
	
	\end{scriptsize}
	\end{tikzpicture}
\end{center}

In \cite{DRS}, Dowbor, Ringel and Simson proved that a hereditary artinian ring $\Lambda$ is representation-finite if and only if the underlying Coxeter valued graph of the Coxeter valued quiver $(\mathscr{C}_{\Lambda}, \mathbf{m})$ of $\Lambda$ is a disjoint union of the Coxeter diagrams presented in Table B, where the valued edge \dynkin[Coxeter, edge length=1.0cm, gonality=3]{G}{2} is identified with \dynkin[Coxeter, edge length=1.0cm]{A}{2}, see \cite{hum}:
\begin{center}
\textbf{Table B.} Coxeter diagrams
\end{center}

	\definecolor{qqqqff}{rgb}{0.,0.,1.}
	\begin{tikzpicture}[line cap=round,line join=round,>=triangle 45,x=1.0cm,y=1.0cm]
	\draw (2.,9.)-- (3.,9.);
	\draw [dash pattern=on 2pt off 2pt](3.,9.)--(4.,9.);
	\draw (4.,9.)-- (5.,9.);
		\begin{scriptsize}
	\draw[color=black] (.8,9.) node {$\mathbb{A}_n:$};
	\draw [fill=white] (2.,9.) circle (1.5pt);
	\draw [fill=white] (3.,9.) circle (1.5pt);
	\draw [fill=white] (4.,9.) circle (1.5pt);
	\draw [fill=white] (5.,9.) circle (1.5pt);
	\draw[color=black] (7.,9.) node {($n ~~{\rm vertices},~ n\geq 1$);};
	\end{scriptsize}
	\end{tikzpicture}

\definecolor{qqqqff}{rgb}{0.,0.,1.}
	\begin{tikzpicture}[line cap=round,line join=round,>=triangle 45,x=1.0cm,y=1.0cm]
	\draw (2.,9.)-- (3.,9.);
	\draw [dash pattern=on 2pt off 2pt](3.,9.)--(4.,9.);
	\draw (4.,9.)-- (5.,9.);
		\begin{scriptsize}
	\draw[color=black] (.8,9.) node {$\mathbb{B}_n:$};
	\draw [fill=white] (2.,9.) circle (1.5pt);
	\draw[color=black] (2.5,9.2) node {$4$};
	\draw [fill=white] (3.,9.) circle (1.5pt);
	\draw [fill=white] (4.,9.) circle (1.5pt);
	\draw [fill=white] (5.,9.) circle (1.5pt);
	\draw[color=black] (7.,9.) node {($n ~~{\rm vertices},~ n\geq 2$);};
	
	\end{scriptsize}
	\end{tikzpicture}

	\definecolor{qqqqff}{rgb}{0.,0.,1.}
	\begin{tikzpicture}[line cap=round,line join=round,>=triangle 45,x=1.0cm,y=1.0cm]
	\draw (2.,9.)-- (3.,9.);
	\draw [dash pattern=on 2pt off 2pt](3.,9.)--(4.,9.);
	\draw (4.,9.)-- (5.,9.);
	\draw (5.,9.)-- (6.,9.);
	\draw (5.,9.)-- (6.,10.9999999999987,2.6);
		\begin{scriptsize}
	\draw[color=black] (.8,9.) node {$\mathbb{D}_n:$};
	\draw [fill=white] (2.,9.) circle (1.5pt);
	\draw [fill=white] (3.,9.) circle (1.5pt);
	\draw [fill=white] (4.,9.) circle (1.5pt);
	\draw [fill=white] (5.,9.) circle (1.5pt);
	\draw [fill=white] (6.,9.) circle (1.5pt);
	\draw [fill=white] (6.,10.9999999999987,2.6) circle (1.5pt);
	\draw[color=black] (8.,9.) node {($n ~~{\rm vertices}, ~n\geq 4$);};
	
	\end{scriptsize}
	\end{tikzpicture}

	\definecolor{qqqqff}{rgb}{0.,0.,1.}
	\begin{tikzpicture}[line cap=round,line join=round,>=triangle 45,x=1.0cm,y=1.0cm]
	\draw (2.,9.)-- (3.,9.);
	\draw (3.,9.)-- (4.,9.);
	\draw (4.,9.)-- (5.,9.);
	\draw (4.,9.)-- (5.,10.9999999999987,2.6);
	\draw (5.,9.)-- (6.,9.);
		\begin{scriptsize}
	\draw[color=black] (.8,9.) node {$\mathbb{E}_6:$};
	\draw [fill=white] (2.,9.) circle (1.5pt);
	\draw [fill=white] (3.,9.) circle (1.5pt);
	\draw [fill=white] (4.,9.) circle (1.5pt);
	\draw [fill=white] (5.,9.) circle (1.5pt);
	\draw [fill=white] (5.,10.9999999999987,2.6) circle (1.5pt);
	\draw [fill=white] (6.,9.) circle (1.5pt);
	\draw[color=black] (6.1,9.) node {;};
	\end{scriptsize}
	\end{tikzpicture}

	\definecolor{qqqqff}{rgb}{0.,0.,1.}
	\begin{tikzpicture}[line cap=round,line join=round,>=triangle 45,x=1.0cm,y=1.0cm]
	\draw (2.,9.)-- (3.,9.);
	\draw (3.,9.)-- (4.,9.);
	\draw (4.,9.)-- (5.,9.);
	\draw (4.,9.)-- (5.,10.9999999999987,2.6);
	\draw (5.,9.)-- (6.,9.);
	\draw (6.,9.)-- (7.,9.);
		\begin{scriptsize}
	\draw[color=black] (.8,9.) node {$\mathbb{E}_7:$};
	\draw [fill=white] (2.,9.) circle (1.5pt);
	\draw [fill=white] (3.,9.) circle (1.5pt);
	\draw [fill=white] (4.,9.) circle (1.5pt);
	\draw [fill=white] (5.,9.) circle (1.5pt);
	\draw [fill=white] (5.,10.9999999999987,2.6) circle (1.5pt);
	\draw [fill=white] (6.,9.) circle (1.5pt);
	\draw [fill=white] (7.,9.) circle (1.5pt);
	\draw[color=black] (7.1,9.) node {;};
	\end{scriptsize}
	\end{tikzpicture}

	\definecolor{qqqqff}{rgb}{0.,0.,1.}
	\begin{tikzpicture}[line cap=round,line join=round,>=triangle 45,x=1.0cm,y=1.0cm]
	\draw (2.,9.)-- (3.,9.);
	\draw (3.,9.)-- (4.,9.);
	\draw (4.,9.)-- (5.,9.);
	\draw (4.,9.)-- (5.,10.9999999999987,2.6);
	\draw (5.,9.)-- (6.,9.);
	\draw (6.,9.)-- (7.,9.);
	\draw (7.,9.)-- (8.,9.);
		\begin{scriptsize}
	\draw[color=black] (.8,9.) node {$\mathbb{E}_8:$};
	\draw [fill=white] (2.,9.) circle (1.5pt);
	\draw [fill=white] (3.,9.) circle (1.5pt);
	\draw [fill=white] (4.,9.) circle (1.5pt);
	\draw [fill=white] (5.,9.) circle (1.5pt);
	\draw [fill=white] (5.,10.9999999999987,2.6) circle (1.5pt);
	\draw [fill=white] (6.,9.) circle (1.5pt);
	\draw [fill=white] (7.,9.) circle (1.5pt);
	\draw [fill=white] (8.,9.) circle (1.5pt);
	\draw[color=black] (8.1,9.) node {;};
	\end{scriptsize}
	\end{tikzpicture}
	
	\definecolor{qqqqff}{rgb}{0.,0.,1.}
	\begin{tikzpicture}[line cap=round,line join=round,>=triangle 45,x=1.0cm,y=1.0cm]
	\draw (2.,9.)-- (3.,9.);
	\draw (3.,9.)-- (4.,9.);
	\draw (4.,9.)-- (5.,9.);
		\begin{scriptsize}
	\draw[color=black] (.8,9.) node {$\mathbb{F}_4:$};
	\draw [fill=white] (2.,9.) circle (1.5pt);
	\draw [fill=white] (3.,9.) circle (1.5pt);
	\draw[color=black] (3.5,9.2) node {$4$};
	\draw [fill=white] (4.,9.) circle (1.5pt);
	\draw [fill=white] (5.,9.) circle (1.5pt);
	\draw[color=black] (5.1,9.) node {;};
	
	\end{scriptsize}
	\end{tikzpicture}
	
	\definecolor{qqqqff}{rgb}{0.,0.,1.}
	\begin{tikzpicture}[line cap=round,line join=round,>=triangle 45,x=1.0cm,y=1.0cm]
	\draw (2.,9.)-- (3.,9.);
		\begin{scriptsize}
	\draw[color=black] (.8,9.) node {$\mathbb{G}_2:$};
	\draw [fill=white] (2.,9.) circle (1.5pt);
	\draw[color=black] (2.5,9.2) node {$6$};
	\draw [fill=white] (3.,9.) circle (1.5pt);
	\draw[color=black] (3.1,9.) node {;};
	
	\end{scriptsize}
	\end{tikzpicture}	
	
	\definecolor{qqqqff}{rgb}{0.,0.,1.}
	\begin{tikzpicture}[line cap=round,line join=round,>=triangle 45,x=1.0cm,y=1.0cm]
	\draw (2.,9.)-- (3.,9.);
	\draw (3.,9.)-- (4.,9.);
		\begin{scriptsize}
	\draw[color=black] (.8,9.) node {$\mathbb{H}_3:$};
	\draw [fill=white] (2.,9.) circle (1.5pt);
	\draw [fill=white] (3.,9.) circle (1.5pt);
	\draw [fill=white] (4.,9.) circle (1.5pt);
	\draw[color=black] (3.5,9.2) node {$5$};
	\draw[color=black] (4.1,9.) node {;};
	
	\end{scriptsize}
	\end{tikzpicture}	

\definecolor{qqqqff}{rgb}{0.,0.,1.}
	\begin{tikzpicture}[line cap=round,line join=round,>=triangle 45,x=1.0cm,y=1.0cm]
	\draw (2.,9.)-- (3.,9.);
	\draw (3.,9.)-- (4.,9.);
	\draw (4.,9.)-- (5.,9.);
		\begin{scriptsize}
	\draw[color=black] (.8,9.) node {$\mathbb{H}_4:$};
	\draw [fill=white] (2.,9.) circle (1.5pt);
	\draw [fill=white] (3.,9.) circle (1.5pt);
	\draw [fill=white] (4.,9.) circle (1.5pt);
	\draw [fill=white] (5.,9.) circle (1.5pt);
	\draw[color=black] (4.5,9.2) node {$5$};
	\draw[color=black] (5.1,9.) node {;};
	
	\end{scriptsize}
	\end{tikzpicture}	
	
	\definecolor{qqqqff}{rgb}{0.,0.,1.}
	\begin{tikzpicture}[line cap=round,line join=round,>=triangle 45,x=1.0cm,y=1.0cm]
	\draw (2.,9.)-- (3.,9.);
		\begin{scriptsize}
	\draw[color=black] (.9,9.) node {$\mathbb{I}_2(p):$};
	\draw [fill=white] (2.,9.) circle (1.5pt);
	\draw[color=black] (2.,8.7) node {$1$};
	\draw[color=black] (2.5,9.2) node {$p$};
	\draw [fill=white] (3.,9.) circle (1.5pt);
	\draw[color=black] (3.,8.7) node {$2$};
	\draw[color=black] (5.,9.) node {($p = 5 ~{\rm or}~ 7 \leq p < \infty$).};
	
	\end{scriptsize}
	\end{tikzpicture}
	
Let $\Lambda$ be a basic hereditary ring and $\mathcal{D}$ be the underlying Coxeter valued graph of the Coxeter valued quiver $(\mathscr{C}_{\Lambda}, \mathbf{m})$ of $\Lambda$.  The ring $\Lambda$ is called {\it of the Dynkin type} $\mathcal{D}$ if $\mathcal{D}$ is one of the Coxeter diagrams $\mathbb{A}_n$, $\mathbb{B}_n$, $\mathbb{D}_n$, $\mathbb{E}_6$, $\mathbb{E}_7$, $\mathbb{E}_8$, $\mathbb{F}_4$ and $\mathbb{G}_2$  presented in Table B. Also $\Lambda$ is called {\it of the Coxeter type} $\mathcal{D}$ if $\mathcal{D}$ is one of the Coxeter diagrams $\mathbb{H}_3$, $\mathbb{H}_4$ and $\mathbb{I}_2(p)$ presented in Table B.\\

It is proved by Schofield \cite{Schofield2, Schofield1} that there exist hereditary bimodule rings of the form ${\Lambda=\tiny{\begin{pmatrix}
F & M\\
0& G
\end{pmatrix}}}$ of the Coxeter type $\mathbb{I}_2(5)$.  However, the existence of such hereditary rings of the Coxeter type $\mathbb{I}_2(p)$ with $p \geq 7$ remains open. It depends on rather difficult questions concerning division rings extensions and leads to a generalized Artin problems, see \cite{Sim5} and \cite{simson}.

One of the main aims of this paper is to get a diagrammatic characterization of right K\"othe rings $\Lambda$ that are basic hereditary, or the square of the Jacobson radical of $\Lambda$ is zero.

The paper is organized as follows. In Section 2, we prove some preliminary results that will be needed later in the paper. In Section 3, we collect some results of hereditary representation-finite rings that we need in the rest of the paper.  In Section 4, we give a characterization of basic hereditary right K\"othe rings of Dynkin type. In Section 5, we give a characterization of basic hereditary right K\"othe rings in terms of their Coxeter valued quivers. Finally in Section 6, we give a characterization of basic right K\"othe rings with radical square zero in terms of their separated quivers.

\subsection{Notation }
 Throughout this paper $\Lambda$ is an associative ring with unit and all modules are unital.  We denote by Mod-$\Lambda$ (resp., mod-$\Lambda$) the category of all right $\Lambda$-modules (resp., finitely generated right $\Lambda$-modules) and by $J$ the Jacobson radical of $\Lambda$.  A ring $\Lambda$ is said to be {\it basic} if $\Lambda/J$ is a direct product of division rings. For a right $\Lambda$-module $M$, we denote by ${\rm top}(M)$ and ${\rm rad}(M)$ its top and radical, respectively. Let $X$ be a representation of a species ${\mathscr{M}}$ (see Section 2) and $V$ be a right module over a division ring $G$. We denote by $\textbf{dim~} X$ and ${\rm dim}(V)_G$ the dimension vector of $X$ and dimension of $V$, respectively. Let $F$ and $G$ be division rings and $M$ be an $F$-$G$-bimodule. We denote ${\rm dim}{_F}{(M)}$ by ${\rm l.dim}~M$ and ${\rm dim}(M)_G$ by r.dim~$M$. Also we denote $M^L:={\rm Hom}_F(M, F)$ and $M^R:= {\rm Hom}_G(M, G)$. We denote the left and right dualisations of the $F$-$G$-bimodule $M$ by setting $M^{(0)} = M$,  $M^{(j)}=(M^{(j-1)})^L$ for $j\geq 1$ and $M^{(j)}= (M^{(j+1)})^R$ for  $j\leq -1$, respectively. Moreover for each $m \geq 1$, we denote by $d_m(M)$ the sequence $\left(d^M_0, d^M_1, \ldots, d^M_{m-1}\right)$, where $d_j^M =$ r.dim $M^{(j)}$ for each $j$. Let $X$ and $Y$ be two right $\Lambda$-modules and $f: X \longrightarrow Y$ be a homomorphism. We write $\xymatrix{ X \ar @{->>}[r] & Y}$ (resp., $\xymatrix{ X \ar@{^{(}->}[r] & Y}$) when $f$ is an epimorphism (resp., $f$ is a monomorphism). Let $A$ be a $\Lambda$-submodule of $X$. We denote by $f|_A$ the restriction of $f$ to $A$. Let $M$ be a right $\Lambda$-module and $n \in {\Bbb{N}}$. Then $M^n$ denotes the direct sum of $n$ copies of $M$. Let $Q$ be a quiver and $i$ be a vertex of $Q$. We denote by $i^+$ and $i^-$ the set of direct successors of $i$ and the set of direct predecessors of $i$, respectively. Also, we denote by $|i^+|$ and  $|i^-|$ the cardinal number of $i^+$ and $i^-$, respectively. Let $K$ be a field, $n \in {\Bbb{N}}$ and $1 \leq i \leq n$. We denote by $e_i$ the vector in $K^n$ with a  1 in the $i$-th coordinate and zero in the $j$-th coordinate, for each $j\neq i$. Throughout the paper we use standard quiver representation and path algebra terminology as applied in the monographs \cite{assem} and \cite{APR}. We also refer to \cite[Ch. VII]{assem} for a detailed explanation of the reflection functors technique introduced by Bernstein, Gelfand, Ponomariev in \cite{bgp}, later developed by Dlab and Ringel in \cite{DR1} and Auslander-Platzeck-Reiten in \cite{APR}.

\section{Preliminaries}
 We recall from \cite{DR1, DRS, DS1} that a {\it valued quiver} $(\Gamma, \mathbf{d})$ of a species ${\mathscr{M}}=(F_i, {_iM_j})_{i,j\in I}$ is a finite quiver $\Gamma= (\Gamma_0, \Gamma_1)$, with the finite set $\Gamma_0$ of vertices corresponding to the division rings $F_i$ together with non-negative integers $d_{ij} = {\rm r.dim}~{_iM_j}$ and $d_{ji} = {\rm l.dim}~{_iM_j}$ for $i \neq j$ and the set $\Gamma_1$ of valued arrows defined as follows.  There exists a valued arrow \begin{center}
	\definecolor{qqqqff}{rgb}{0.,0.,1.}
	\begin{tikzpicture}[line cap=round,line join=round,>=triangle 45,x=1.0cm,y=1.0cm]
	\draw (2.,9.)-- (3.,9.);
		\begin{scriptsize}
	\draw[color=black] (.8,9.) node {};
	\draw [fill=white] (2.,9.) circle (1.5pt);
	\draw[color=black] (2.,8.7) node {$i$};
	\draw[color=black] (2.5,9.3) node {$(d_{ij}, d_{ji})$};
	\draw[color=black] (2.85,9.) node {$>$};
	\draw [fill=white] (3.,9.) circle (1.5pt);
	\draw[color=black] (3.,8.7) node {$j$};
	
	\end{scriptsize}
	\end{tikzpicture}
\end{center}
if and only if ${_iM_j} \neq 0$. If $d_{ij} = d_{ji} = 1$, we write simply
\begin{center}
	\definecolor{qqqqff}{rgb}{0.,0.,1.}
	\begin{tikzpicture}[line cap=round,line join=round,>=triangle 45,x=1.0cm,y=1.0cm]
	\draw (2.,9.)-- (3.,9.);
		\begin{scriptsize}
	\draw[color=black] (.8,9.) node {};
	\draw [fill=white] (2.,9.) circle (1.5pt);
	\draw[color=black] (2.,8.7) node {$i$};
	\draw[color=black] (2.85,9.) node {$>$};
	\draw [fill=white] (3.,9.) circle (1.5pt);
	\draw[color=black] (3.,8.7) node {$j$};
	
	\end{scriptsize}
	\end{tikzpicture}
\end{center}
The valued quiver of the species of a basic ring $\Lambda$ is denoted by $(\Gamma_{\Lambda}, \mathbf{d})$. A valued quiver is called {\it connected} if the underlying valued graph of the valued quiver is connected. A valued quiver is said to be {\it acyclic} if it has no cycles. In this paper, unless otherwise stated, we assume that ${\mathscr{M}}=(F_i, {_iM_j})_{i,j\in I}$ is a species with the property ``${{_iM_j}\neq 0}$ implies that ${_jM_i} = 0$". Note that the species of any  basic hereditary right artinian ring has this property, see \cite{Drozd}. Following \cite{DR1}, \cite{ringel} and \cite{DRS, DS1, Sim3}, a {\it representation} of ${\mathscr{M}}$ is a family $\left( X_i,{_j\varphi_i} \right)_{i, j \in I}$ of right $F_i$-modules $X_i$ and $F_j$-linear maps ${_{j}{\varphi}_{i}}: X_i \otimes_{F_i}{_iM_j} \rightarrow X_j$ for each arrow $i \rightarrow  j$ of $\Gamma$.  A representation $\left( X_i,{_j\varphi_i} \right)$ is called {\it finite-dimensional} provided that all $X_i$ are finite-dimensional right $F_i$-modules. A morphism $\alpha: \left( X_i,{_j\varphi_i} \right) \rightarrow \left( Y_i,{_j\psi_i} \right)$ is given by right $F_i$-linear maps $\alpha_i: X_i \rightarrow Y_i$ such that $_j\psi_i(\alpha_i \otimes id_{{_iM_j}})= \alpha_j{_j\varphi_i}$ for each arrow $i \rightarrow  j$ of $\Gamma$. We denote by ${\rm Rep}({\mathscr{M}})$ (resp.,  by ${\rm rep}({\mathscr{M}})$) the category of all representations of ${\mathscr{M}}$ (resp., finite-dimensional representations of ${\mathscr{M}}$) \cite{3}. Obviously representations of species are a generalization of representations of quivers.\\

Following \cite{Anderson}, a submodule $K$ of a right $\Lambda$-module $N$ is called {\it small submodule} if for every submodule $L$ of $N$, $K+L=N$ implies $L=N$. Let ${\mathscr{M}}= \left( F_i,{_iM_j} \right)_{i, j \in I} $ be a species with the valued quiver $(\Gamma, \mathbf{d})$ and $\left( X_i,{_j\varphi_i} \right)$ be a representation of ${\mathscr{M}}$. Suppose that $\left( Y_i,{_j{\psi}_i} \right)$ is a subrepresentation of the representation $\left( X_i,{_j\varphi_i} \right)$ of ${\mathscr{M}}$.  Then $\left( X_i,{_j\varphi_i} \right)/ \left( {Y}_i, {_j{\psi}_i} \right) = \left( X_i/{Y}_i, {_j{\gamma}_i} \right)$, where for each arrow $ i \rightarrow j$ of $\Gamma$ the maps ${_j{\gamma}_i}$ are defined by the following commutative diagram
\begin{displaymath}
\xymatrix{
0 \ar[r] & {Y}_i \otimes_{F_i}{_iM_j} \ar[r]^{\imath_i \otimes id_{_iM_j}} \ar[d]_{{_j{\psi}_i}} &
{X}_i \otimes_{F_i}{_iM_j} \ar[r]^{\pi_i \otimes id_{_iM_j}}  \ar[d]_{{_j{\varphi}_i}} & X_i/{Y}_i \otimes_{F_i}{_iM_j} \ar[d]_{_j{\gamma}_i} \ar [r] & 0\\
0 \ar [r] & {Y}_j \ar[r]^{\imath_j} & X_j \ar[r]^{\pi_j} & X_j/{Y}_j \ar[r]& 0}
\end{displaymath}
in which $\imath_i: {Y}_i \rightarrow X_i$ denote the inclusion map and $\pi_i: X_i \rightarrow X_i/{Y}_i$ projection. Set $\pi=(\pi_i)_i: \left( X_i,{_j\varphi_i} \right) \rightarrow \left( X_i/{Y}_i, {_j{\gamma}_i} \right)$. The representation $\left( Y_i,{_j{\psi}_i} \right)$ is called {\it a small subrepresentation} of the representation $\left( X_i,{_j\varphi_i} \right)$ if every morphism $\alpha: \left( Z_i,{_j\phi_i} \right) \rightarrow \left( X_i,{_j\varphi_i} \right)$ in ${\rm rep}({\mathscr{M}})$ with $\pi \alpha$ epic is epic.\\

The following proposition gives the necessary and sufficient conditions for a representation $\left( Y_i,{_j{\psi}_i} \right)$ to be a small subrepresentation of $\left( X_i,{_j\varphi_i} \right)$.

\begin{Pro}\label{3}
Let ${\mathscr{M}}= \left( F_i,{_iM_j} \right)_{i, j \in I} $ be a species. Suppose that the valued quiver $(\Gamma, \mathbf{d})$ of ${\mathscr{M}}$ is connected and acyclic and let $\left( Y_i,{_j{\psi}_i} \right)$ be a subrepresentation of a  representation $\left( X_i,{_j\varphi_i} \right)$ of $\mathscr{M}$. Then
$\left( Y_i,{_j{\psi}_i} \right)$  is a small subrepresentation of $\left( X_i,{_j\varphi_i} \right)$ if and only if the following conditions hold:
\begin{itemize}
\item[$(a)$] If $k$ is a source vertex of $\Gamma$, then $Y_k = 0$.
\item[$(b)$] If $k$ is not a source vertex of $\Gamma$, then $Y_k \subseteq {\sum_{j \rightarrow k}{\rm Im}(_k\varphi_j)}$, where the sum is over all arrows with the target $k$.
\end{itemize}
\end{Pro}
\begin{proof}
$(\Rightarrow).$ Assume that $k$ is a source vertex of  $\Gamma$ and $X_k \neq 0$. Assume that $Y_k = X_k$. We define a representation $ \left( Z_i,{_j{\theta}_i} \right)$ of ${\mathscr{M}}$ by taking $Z_k=0$, $Z_i = X_i$ for each $i \neq k$ and all ${_j{\theta}_i} = {_j\varphi_i}|_{Z_i \otimes_{F_i}{_iM_j}}$. Then $\pi \alpha$ is an epimorphism, where $\alpha= (\alpha_i)_i: \left( Z_i,{_j{\theta}_i} \right) \rightarrow \left( X_i,{_j\varphi_i} \right)$ and each $\alpha_i:Z_i \rightarrow X_i$ is the inclusion map. Since $\alpha$ is not an epimorphism, $\left( Y_i,{_j{\psi}_i} \right)$ is not a small subrepresentation of $\left( X_i,{_j\varphi_i} \right)$, which gives a contradiction. Hence $Y_k \neq X_k$. Assume that $H_k$ is an $F_k$-submodule of $X_k$ such that $H_k + Y_k = X_k$. We define a subrepresentation $( X^{'}_i, {_j\varphi^{'}_i} )$ of $\left( X_i,{_j\varphi_i} \right)$ by taking $X^{'}_k = H_k$, $X^{'}_i = X_i$ for each $i \neq k$ and all ${_j\varphi^{'}_i} = {_j\varphi_i}|_{X^{'}_i \otimes_{F_i}{_iM_j}}$. Then  $\pi \ell$ is an epimorphism, where $\ell= (\ell_i)_i: ( X^{'}_i, {_j\varphi^{'}_i} ) \rightarrow \left( X_i,{_j\varphi_i} \right)$ and each $\ell_i: X^{'}_i \rightarrow X_i$ is the inclusion map. The assumption $\left( Y_i,{_i{\psi}_j} \right)$ is a small subrepresentation of $\left( X_i,{_i\varphi_j} \right)$ yields the map $\ell$ is an epimorphism and so $H_k = X_k$. This proves that $Y_k$ is a small submodule of $ X_k$. Since $X_k$ is a right $F_k$-module and $F_k$ is a division ring, $Y_k$ is a direct summand of $X_k$. Therefore $Y_k=0$. Assume that $k$ is not a source vertex of $\Gamma$ and  $\eta_k: X_k \rightarrow X_k/{\sum_{j \rightarrow k}{\rm Im}(_k\varphi_j)}$ is the canonical quotient map, where the sum is over all arrows with the target $k$. Assume that $\eta_k(X_k) \neq 0$ and $\eta_k(Y_k) = \eta_k(X_k)$. Then $X_k = Y_k + {\sum_{j \rightarrow k}{\rm Im}(_k\varphi_j)}$. We define a subrepresentation $\left( V_i,{_j\varepsilon_i} \right)$ of $\left( X_i,{_j\varphi_i} \right)$ by taking  $V_k = {\sum_{j \rightarrow k}{\rm Im}(_k\varphi_j)}$, $V_i = X_i$ for each $i \neq k$ and ${_j\varepsilon_i} = {_j\varphi_i}|_{V_i \otimes_{F_i}{_iM_j}}$ for each $i$ and $j$. Therefore $\pi \beta$ is an epimorphism, where $\beta=(\beta_i)_i: \left( V_i,{_j\varepsilon_i} \right) \rightarrow \left( X_i,{_j\varphi_i} \right)$ and each $\beta_i: V_i \rightarrow X_i$ is the inclusion map. Thus $\left( Y_i,{_j{\psi}_i} \right)$ is not a small subrepresentation of $\left( X_i,{_j\varphi_i} \right)$, which gives a contradiction.  It follows that $\eta_k(Y_k) \neq \eta_k(X_k)$.  Assume that  $E_k$ is an $F_k$-submodule of $\eta_k(X_k)$ such that $E_k + \eta_k(Y_k) = \eta_k(X_k)$. Set $D_k=\eta^{-1}_k(E_k)$. Then $_k\varphi_j(X_j\otimes_{F_j}{_jM_k}) \subset D_k$ for each arrow $ j \rightarrow k$ of $\Gamma$.  Hence $D_k + Y_k = X_k$ and $({X^{''}_i}, {_j\varphi^{''}_i})$  is a subrepresentation of $\left( X_i,{_j{\varphi}_i}\right)$, where ${X^{''}_k} = D_k$, ${X^{''}_i} = X_i$ for each $i \neq k$ and all ${_j\varphi^{''}_i} = {_j\varphi_i}|_{{X^{''}_i} \otimes_{F_i}{_iM_j}}$. Therefore $\pi \xi$ is an epimorphism, where $\xi= (\xi_i)_i: ( X^{''}_i, {_j\varphi^{''}_i} ) \rightarrow \left( X_i,{_j\varphi_i} \right)$ and each $\xi_i: X^{''}_i \rightarrow X_i$ is the inclusion map.
The assumption $\left( Y_i,{_i{\psi}_j} \right)$ is a small subrepresentation of $\left( X_i,{_i\varphi_j} \right)$ yields  $D_k = X_k$. It follows that $E_k=\eta_k(X_k)$. Consequently, $\eta_k(Y_k)$  is a small submodule of $\eta_k(X_k)$. Therefore $Y_k \subseteq {\sum_{j \rightarrow k}{\rm Im}(_k\varphi_j)}$.\\
$(\Leftarrow)$. Let $\left( Z_i,{_j{\phi}_i} \right)$ be a representation of ${\mathscr{M}}$ and $f=(f_i): \left( Z_i,{_j{\phi}_i} \right) \rightarrow \left( X_i,{_j\varphi_i} \right)$ be a morphism in ${\rm Rep}({\mathscr{M}})$ such that $\pi f$ is an epimorphism. Then for each  vertex $i$ of $\Gamma$,  ${\rm Im}f_i + Y_i = X_i$.  If $k$ is a source vertex of $\Gamma$, then the assumption $(a)$ yields $X_k={\rm Im}f_k$. Assume that $k$ is not a source vertex of $\Gamma$. Let $\eta_k: X_k \rightarrow X_k/{\sum_{j \rightarrow k}{\rm Im}(_k\varphi_j)}$ be the canonical quotient map. Since ${\rm Im}f_k + Y_k = X_k$, $\eta_k({\rm Im}f_k) + \eta_k(Y_k) = \eta_k(X_k)$. Hence the assumption $(b)$ yields $\eta_k({\rm Im}f_k) = \eta_k(X_k)$. Therefore it is sufficient to show that for each arrow $j \rightarrow k$ of $\Gamma$, $_k\varphi_j(X_j\otimes_{F_j}{_jM_k}) \subseteq {\rm Im}f_k$. Assume that there exists an arrow $j_1 \rightarrow k$ of $\Gamma$ such that $_k\varphi_{j_1}(X_{j_1}\otimes_{F_{j_1}}{_{j_1}M_k}) \nsubseteq {\rm Im}f_k$. Since ${_j\varphi_i}({f_i}\otimes id_{_iM_j}) = f_j {_j\phi_i}$ for each arrow $i \rightarrow j$ of $\Gamma$, we have a representation $({\rm Im}{f_i}, {_j{\delta}_i})$, where ${_j{\delta}_i} = {_j\varphi_i}|_{{\rm Im}{f_i}\otimes_{F_i}{_iM_j}}$. Thus for each arrow $i \rightarrow j$ of $\Gamma$, we have the following commutative diagram
$$\xymatrix@1{
\ar@{}[dr]
{\rm Im}{f_i}\otimes_{F_i}{_iM_j}  \ar@{^{(}->}[d]_{\ell_i \otimes 1} \ar[r]^<(.3){_j{\delta}_i} & {\rm Im}f_j \ar@{^{(}->}[d]^{\ell_j} \\
X_i \otimes_{F_i}{_iM_j} \ar[r]^<(.3){_j\varphi_i} & {X_j}}$$
where each $\ell_i: {\rm Im}{f_i} \rightarrow X_i$ is the inclusion map. Therefore  ${\rm Im}{f_{j_1}}\neq X_{j_1}$. Since ${\rm Im}{f_{j_1}} + Y_{j_1} = X_{j_1}$, $Y_{j_1}$ is not a small submodule of $X_{j_1}$. By the above argument $j_1$ is not a source vertex of $\Gamma$ and by (b), $\eta_{j_1}(Y_{j_1})=0$. Since ${{\rm Im}{f_{j_1}} + Y_{j_1} = X_{j_1}}$, ${\eta_{j_1}({\rm Im}{f_{j_1}}) + \eta_{j_1}(Y_{j_1}) = \eta_{j_1}(X_{j_1})}$. It follows that $\eta_{j_1}({\rm Im}{f_{j_1}}) = \eta_{j_1}(X_{j_1})$. Therefore there exists an arrow $j_2 \rightarrow j_1$ of $\Gamma$ such that $_{j_1}\varphi_{j_2}(X_{j_2}\otimes_{F_{j_2}}{_{j_2}M_{j_1}}) \nsubseteq {\rm Im}{f_{j_1}}$. It follows that $j_2$ is not a source vertex of $\Gamma$ and the same argument shows that there exists an arrow $j_3 \rightarrow j_2$ of $\Gamma$ such that $_{j_2}\varphi_{j_3}(X_{j_3}\otimes_{F_{j_3}}{_{j_3}M_{j_2}}) \nsubseteq {\rm Im}{f_{j_2}}$. Continuing in this way, we get a cycle in $(\Gamma, \mathbf{d})$ which gives a contradiction. Then $_k\varphi_j(X_j\otimes_{F_j}{_jM_k}) \subseteq {\rm Im}{f_k}$ for each arrow $ j \rightarrow k$ of $\Gamma$. It follows that $f$ is an epimorphism. Consequently $\left( Y_i,{_j{\psi}_i} \right)$ is a small subrepresentation of $\left( X_i,{_j\varphi_i} \right)$.
\end{proof}

 In the following proposition we compute the radical and top of a representation of species.

\begin{Pro}\label{4}
Let ${\mathscr{M}}= \left( F_i,{_iM_j} \right)_{i, j \in I} $ be a species. Suppose that the valued quiver $(\Gamma, \mathbf{d})$ of ${\mathscr{M}}$ is connected and acyclic and let $\left( X_i,{_i\varphi_j} \right)$ be a finite-dimensional representation of ${\mathscr{M}}$. Then
\begin{itemize}
\item[$(a)$] ${\rm rad}(\left( X_i,{_j\varphi_i} \right)) = \left( Y_i,{_j\psi_i} \right)$, where $Y_k = \sum_{j \rightarrow k}{\rm Im}(_k\varphi_j)$ if $k$ is not a source vertex of $\Gamma$, $Y_k=0$ if $k$ is a source vertex of $\Gamma$ and ${_j\psi_i}= {_j\varphi_i}|_{Y_i \otimes_{F_i}{_iM_j}}$ for each $i$ and $j$.
\item[$(b)$] ${\rm top}(\left( X_i,{_i\varphi_j} \right)) = \left( Z_i,{_j\gamma_i} \right)$, where $Z_k = X_k/ \sum_{j \rightarrow k}{\rm Im}(_k\varphi_j)$ if $k$ is not a source vertex of $\Gamma$, $Z_k = X_k$ if $k$ is a source vertex of $\Gamma$ and ${_j\gamma_i}= 0$ for each $i$ and $j$.
\end{itemize}
\end{Pro}
\begin{proof}
For each vertex $k$ of $\Gamma$, define the representation ${\bf{F}}_k=(W_i, {_j\chi_i})$, where $W_i = 0$ for $i \neq k$, $W_k = F_k$ and all ${_j\chi_i}=0$.  Since the species $\mathscr{M}$ is acyclic, by \cite[Proposition 1.1]{Sim3}, the simple representations ${\bf{F}}_k$ form a complete list of all non-isomorphic simple objects in ${\rm rep}({\mathscr{M}})$, where $k$ is a vertex of $\Gamma$. Therefore $(V_i, {_j\phi_i})$ is a maximal subrepresentation of $(X_i, {_j\varphi_i})$  if and only if there exists a vertex $k$ of $\Gamma$ such that $(V_i, {_j\phi_i})$ satisfies one of the following conditions:
\begin{itemize}
\item[$(a')$] $k$ is a source vertex of $\Gamma$, $V_k$ is a maximal submodule of $X_k$,  $V_i = X_i$ for each $i \neq k$, and ${_j\phi_i} = {_j\varphi_i}|_{{V}_i \otimes_{F_i}{_iM_j}}$ for each $i$ and $j$,
\item[$(b')$] $k$ is not a source vertex of $\Gamma$, $V_k$ is a maximal submodule of $X_k$ which contains\\  $\sum_{j \rightarrow k}{\rm Im}(_k\varphi_j)$,  $V_i = X_i$ for each $i \neq k$, and ${_j\phi_i} = {_j\varphi_i}|_{{V}_i \otimes_{F_i}{_iM_j}}$ for each $i$ and $j$.
\end{itemize}
It follows that ${\rm rad}(\left( X_i,{_j\varphi_i} \right)) = \left( Y_i,{_j\psi_i} \right)$, where $Y_i = \sum_{j \rightarrow i}{\rm Im}(_i\varphi_j)$ if $i$ is not a source vertex of $\Gamma$, whereas $Y_i=0$ if $i$ is a source vertex of $\Gamma$ and ${_j\psi_i}= {_j\varphi_i}|_{Y_i \otimes_{F_i}{_iM_j}}$ for each $i$ and $j$. Consequently, ${\rm top}(\left( X_i,{_i\varphi_j} \right)) = \left( Z_i,{_j\gamma_i} \right)$, where $Z_i = X_i/ \sum_{j \rightarrow i}{\rm Im}(_i\varphi_j)$ if $i$ is not a source vertex of $\Gamma$, whereas $Z_i = X_i$ if $i$ is a source vertex of $\Gamma$ and ${_j\gamma_i}= 0$ for each $i$ and $j$.
\end{proof}

As an immediate consequence of Propositions \ref{3} and \ref{4}, we obtain the following corollary.

 \begin{Cor}\label{m77}
Let ${\mathscr{M}}= \left( F_i,{_iM_j} \right)_{i, j \in I} $ be a species. Suppose that the valued quiver $(\Gamma, \mathbf{d})$ of ${\mathscr{M}}$ is connected and acyclic and  let $\left( X_i,{_i\varphi_j} \right)$ be a finite-dimensional representation of ${\mathscr{M}}$. Then ${\rm rad}(\left( X_i,{_j\varphi_i} \right))$ is the biggest small subrepresentation of $\left( X_i,{_i\varphi_j} \right)$.
 \end{Cor}

It is well known that, for a quiver $Q=(Q_0, Q_1)$, the category of representations of $Q$ over a field $K$ is equivalent to Mod-$KQ$, where $KQ$ is the path $K$-algebra of $Q$, see \cite{assem}. This fact was generalized nicely for species. For a species ${\mathscr{M}}= \left( F_i,{_iM_j} \right)_{i, j \in I}$, one can form a tensor ring $R_{\mathscr{M}}= \bigoplus_{t\geq 0}N^{(t)}$, where  $A = N^{(0)} = \bigoplus_{i \in I}F_i$, $N^{(1)}= \bigoplus_{i, j \in I}{_iM_j}$ and $N^{(t)}= N^{(t-1)} \otimes_{A} N^{(1)} $ for $t \geq 2$, with the component-wise addition and the multiplication induced by taking tensor products. The ring $R_{\mathscr{M}}$ is called the {\it tensor ring of} ${\mathscr{M}}$ (see \cite{3}). Following \cite{DS1, Sim3}, a species ${\mathscr{M}}$ is called {\it right {\rm (}resp., left{\rm)} finite-dimensional} if the dimensions $d_{ij}$ (resp., $d_{ji}$) are finite for all $i \neq j$. The species ${\mathscr{M}}$ is called {\it finite-dimensional} if ${\mathscr{M}}$ is left and right finite-dimensional. Assume that ${\mathscr{M}}= \left( F_i,{_iM_j} \right)_{i, j \in I}$ is a right finite-dimensional species and the valued quiver $(\Gamma, \mathbf{d})$ of ${\mathscr{M}}$ is acyclic. We now define as in \cite{2} a functor
\begin{equation}\label{f}
\mathrm{F}: {\rm rep}({\mathscr{M}}) \longrightarrow {\rm mod-}R_{\mathscr{M}}
\end{equation}
as follows. For each object $X=\left( X_i,{_j\varphi_i} \right)$ in ${\rm rep}({\mathscr{M}})$, we set $\mathrm{F}(X) = \bigoplus_{i \in I}X_i$. The reader may easily verify that $\mathrm{F}(X)$ is a finitely generated right $R_{\mathscr{M}}$-module. If $(\alpha_i)_{i \in I}: \left( X_i,{_j\varphi_i} \right) \rightarrow \left( Y_i,{_j{\psi}_i} \right)$ is a morphism in ${\rm rep}({\mathscr{M}})$, we define $\mathrm{F}((\alpha_i)_{i \in I})$ to be $\bigoplus_{i \in I} \alpha_i: \bigoplus_{i \in I}X_i \rightarrow \bigoplus_{i \in I}Y_i$. It is easy to verify that $\mathrm{F}((\alpha_i)_{i \in I})$ is an $R_{{\mathscr{M}}}$-homomorphism.  It is well known that the functor $\mathrm{F}$ is an equivalence (see \cite[Proposition 1.1]{Sim3}). From now on, we fix  the functor $\mathrm{F}$. Note that when the species ${\mathscr{M}}$ is right finite-dimensional, then the tensor ring $R_{\mathscr{M}}$ is semiprimary. In this case,  Propositions \ref{3} and \ref{4} and Corollary \ref{m77} are an automatic translation of categorical properties from modules to representations. In particular, we have the following fact.

\begin{Pro}\label{5}
Let ${\mathscr{M}}= \left( F_i,{_iM_j} \right)_{i, j \in I} $ be a right finite-dimensional species. Suppose that the valued quiver $(\Gamma, \mathbf{d})$ of ${\mathscr{M}}$ is connected and acyclic and let $\left( X_i,{_j\varphi_i} \right)$ be a finite-dimensional representation of ${\mathscr{M}}$. Then ${\rm top}(\mathrm{F}\left( X_i,{_j\varphi_i} \right)) \cong \mathrm{F}({\rm top}(\left( X_i,{_j\varphi_i} \right)))$.
\end{Pro}

\section{Representations of species and modules over hereditary artinian rings}

Let ${\mathscr{M}} = \left( F_i,{_iM_j} \right)_{i, j \in I} $ be a finite-dimensional species and suppose that the valued quiver $(\Gamma, \mathbf{d})$ of ${\mathscr{M}}$  is acyclic and connected. Let $k$ be a sink (resp., source) vertex of $\Gamma$ and ${\mathscr{M}}^k = \left( F_i,{_iN_j} \right)_{i, j \in I}$ be the new species, where
$$_iN_j=\begin{cases}
{_jM_k}^L & \mbox{if $i=k$},\\
{_iM_j} & \mbox{if $i\neq k$ and $j \neq k$},\\
0 & \mbox{if $j = k$}.
\end{cases}$$
(resp., $_iN_k = {_kM_i}^R$, $_kN_i = 0$ for each $i$ and $_iN_j = {_iM_j}$ for each $i\neq k$, $j \neq k$) (see \cite{DS1}). We recall from \cite{DR1} and \cite{ringel} (see also \cite{APR} and \cite{DRS, DS1, Sim3}) that a pair of {\it partial Coxeter functors} (or {\it reflection functors}, see \cite[Ch. VII]{assem}) $$\xymatrix{{\rm rep}({\mathscr{M}}) \ar@<1ex>[r]^{S_k^+} \ar@{<-}[r]_{S_k^-} & {\rm rep}({\mathscr{M}}^k)}$$ is defined as follows. Let $X= \left( X_i,{_j\varphi_i} \right)$ be a finite-dimensional representation of ${\mathscr{M}}$. Define ${{S_k^+}X=\left( Y_i,{_j\psi_i} \right)}$, where $Y_i = X_i$ for each $i \neq k$ and $Y_k$ is the kernel of the morphism $({_k\varphi_j})_j$  $$0 \rightarrow Y_k \overset{({_j\sigma_k})_j}{\rightarrow } \bigoplus_{j \in \Gamma_0} X_j \otimes_{F_j}{{_jM_k}} \overset{({_k\varphi_j})_j}{\rightarrow} X_k.$$ By using the natural isomorphism ${\rm Hom}_{F_j}(Y_k \otimes_{F_k}{{_jM_k}^L}, X_j) \cong {\rm Hom}_{F_k}(Y_k, X_j \otimes_{F_j}{{_jM_k}})$, we get ${_j\psi_k}: Y_k \otimes_{F_k}{{_jM_k}^L} \rightarrow Y_j$  and ${_j\psi_i} = {_j\varphi_i}$ for $i \neq k$.  Also, if $\alpha=(\alpha_i)_i: X \rightarrow X^{'} $ is a morphism in ${\rm rep}({\mathscr{M}})$, then $S_k^+\alpha= (\beta_i)_i$ is defined by $\beta_i=\alpha_i$ for $i \neq k$ and $\beta_k: Y_k \rightarrow Y_k^{'}$ as the restriction of $$\bigoplus_{j \in \Gamma_0}(\alpha_j \otimes 1): \bigoplus_{j \in \Gamma_0} X_j \otimes {_jM_k}\rightarrow \bigoplus_{j \in \Gamma_0} X_j^{'}\otimes{_jM_k}$$
Also, for each sink vertex $k$ of $\Gamma$, define the linear transformation $s_k^+: {\Bbb{Z}}^n \rightarrow {\Bbb{Z}}^n$ by $s_k^+x=y$, where $|\Gamma_0|= n$, $y_i=x_i$ for $i \neq k$ and $y_k = -x_k + \sum_{i \in \Gamma_0}d_{ik}x_i$. For each finite-dimensional indecomposable representation $X= \left( X_i,{_j\varphi_i} \right)$ of ${\mathscr{M}}$, we can see that each ${_k\varphi_i}$ is an epimorphism. It follows that $\textbf{dim}~S_k^+X = s_k^+(\textbf{dim}~X)$, where $\textbf{dim}~X = \left( {\rm dim}{(X_i)}_{F_i}\right) _{i \in I}$. The functor ${S_k^-}$ is defined analogously. The functors ${S_k^+}$ and ${S_k^-}$ induce quasi-inverse equivalences between
the full subcategory of ${\rm rep}({\mathscr{M}})$ of the representations having no direct summand isomorphic to the simple projective representation $\mathbf{F}_k$, and the full subcategory of ${\rm rep}({\mathscr{M}}^k)$ of the representations having no direct summand isomorphic to the simple injective representation $\mathbf{F}_k$ (see \cite{DR1} and \cite{DS1}).

Let $k$ be a vertex of $\Gamma$ and let $s_k\Gamma$ be the quiver obtained from $\Gamma$ by reversing the direction of all arrows starting or end in $k$, see \cite{DR1} and \cite[Ch. VII]{assem}. An admissible sequence of sinks in a quiver $\Gamma$ is a sequence $k_1, \ldots, k_n$ of all vertices in $\Gamma$  such that each vertex $k_t$ is a sink in $s_{k_{t-1}} \cdots s_{k_1}\Gamma$ for all $1 \leq t \leq n$ \cite{assem, DR1}. Let $k_1, \ldots, k_n$ be an admissible sequence in $\Gamma$ and $({k'_j})_j$ be a sequence of vertices of $\Gamma$, where $j \in {\Bbb{Z}}$,
 $$k'_j = \begin{cases}
k_{r+1}& \mbox{if $j\geq 0, ~j=tn+r, ~0 \leq r <n$};\\
k_{n-r+1}& \mbox{if $j < 0, ~-j=tn+r, ~0 \leq r <n$},
\end{cases}$$
and $t$ is a positive integer. For any $m \in {\Bbb{Z}}$, the species $\mathscr{M}^{(m)}$ is defined in \cite{DS1} as follows.
  $$\mathscr{M}^{(m)}  = \begin{cases}
(\mathscr{M}^{(m-1)})^{k^{'}_{m-1}} & \mbox{if $m \geq 1$};\\
(\mathscr{M}^{(m+1)})^{k^{'}_m} & \mbox{if $m \leq {-1}$},
\end{cases}$$
where $\mathscr{M}^{(0)} = \mathscr{M}$. The species $\mathscr{M}$ has the {\it right {\rm (resp.}, left{\rm)} finite-dimensional property} if the species $\mathscr{M}^{(m)}$ are finite-dimensional for all $m \geq 0$ (resp., $m \leq 0$). $\mathscr{M}$ has the {\it finite-dimensional property} if it has both the left and the right finite-dimensional property. If $\mathscr{M}$ has the right (resp., left) finite-dimensional property, then the following sequence is a right (resp., left) sequence of partial Coxeter functors of $\mathscr{M}$
 $$ \xymatrix{{\rm rep}({\mathscr{M}}) \ar@<1ex>[r]^{S_{1}^+} \ar@{<-}[r]_{S_{1}^-} & {\rm rep}({\mathscr{M}}^{(1)}) \ar@<1ex>[r] \ar@{<-}[r] & \cdots \ar@<1ex>[r] \ar@{<-}[r] & {\rm rep}({\mathscr{M}}^{(m-1)})\ar@<1ex>[r]^{S_{m}^+} \ar@{<-}[r]_{S_{m}^-} & {\rm rep}({\mathscr{M}}^{(m)})\ar@<1ex>[r] \ar@{<-}[r]& \cdots}$$
$${\rm (resp.,}  \xymatrix{ \cdots \ar@<1ex>[r] \ar@{<-}[r]& {\rm rep}({\mathscr{M}}^{(-m-1)}) \ar@<1ex>[r]^{S_{-m}^+} \ar@{<-}[r]_{S_{-m}^-} & {\rm rep}({\mathscr{M}}^{(-m)}) \ar@<1ex>[r] \ar@{<-}[r] & \cdots \ar@<1ex>[r] \ar@{<-}[r] & {\rm rep}({\mathscr{M}}^{(-1)})\ar@<1ex>[r]^{S_{0}^+} \ar@{<-}[r]_{S_{0}^-} & {\rm rep}({\mathscr{M}})})$$
 where $S_j^+$ and $S_j^-$ are the pair of partial Coxeter functors corresponding to the sink $k'_{j-1}$ in the valued quiver of ${\mathscr{M}}^{(j-1)}$. We denote by ${\mathbf{F}}_{k'_j}^{(j)}$ the simple projective representation in ${\mathscr{M}}^{(j)}$ corresponding to the sink $k'_j$ \cite{DS1}.\\

We need the following result proved in \cite[Theorem 1]{DS1}.

\begin{The}\label{DS1}
Let ${\mathscr{M}}$ be a finite-dimensional species and suppose that its valued quiver $(\Gamma, \mathbf{d})$ is connected and acyclic. Then ${\mathscr{M}}$ is  representation-finite if and only if ${\mathscr{M}}$ has the finite-dimensional property and there exists $m>0$ such that $s_m^+\cdots s_1^+(e_i)\ngeqq 0$ for any source $i$ of $\Gamma$. Moreover, if $m$ is minimal with the above property and $|\Gamma_0|=n$, then the mapping $\textbf{dim}: {\rm rep}({\mathscr{M}}) \rightarrow {\Bbb{Z}}^n$ is a one-one correspondence between
isomorphism classes of finite-dimensional indecomposable representations of ${\mathscr{M}}$ and vectors in ${\Bbb{Z}}^n$ of the form $s_1^-\cdots s_{t}^-(e_{k_t^{'}})$, where $t < m$ and $k_t^{'}$ is a sink in the valued quiver of ${\mathscr{M}}^{(t)}$.  In other words any finite-dimensional indecomposable representation $X$ of ${\mathscr{M}}$ has the form $X \cong \textbf{P}_i$ for some $0 \leq i < m$, where $\textbf{P}_0 = {\mathbf{F}}_{k_1}$ and $\textbf{P}_i = {S_1^-} \cdots {S_{i}^-}{\mathbf{F}}_{k'_i}^{(i)}$.
\end{The}

We recall from \cite{DRS} that a sequence $a = (a_1, \ldots, a_m)$ of length $m \geq 2$ with $a_i \in {\Bbb{N}}$ is called a {\it dimension sequence} provided there exist $x_i, y_i \in {\Bbb{N}}$ ($1 \leq i \leq m$), with $$ a_ix_i = x_{i-1} + x_{i+1} \hspace{5mm} {\rm and} \hspace{5mm} a_iy_i = y_{i-1} + y_{i+1} \hspace{1cm} (1 \leq i < m),$$
where $x_0= -1$, $y_0 = y_m = x_1 = 0$ and $x_m = y_1 = 1$.

\begin{Lem}\label{2.12.1}
Let $F$ and $G$ be division rings, $M$ be an $F$-$G$-bimodule and $\Lambda= \tiny{\begin{pmatrix}
F & M\\
0& G
\end{pmatrix}}$ be an artinian ring. Then there exist precisely $3$ pairwise non-isomorphic finitely generated indecomposable right $\Lambda$-modules if and only if  ${\rm r.dim}~M = {\rm l.dim}~M =1$. Moreover in this case  $F \cong G$ as division rings and $M^L \cong M^R$ as $G$-$F$-bimodules.
\end{Lem}
\begin{proof}
Assume that $\Lambda= \tiny{\begin{pmatrix}
F & M\\
0& G
\end{pmatrix}}$ is an artinian ring such that there exist precisely $3$ pairwise non-isomorphic finitely generated indecomposable right $\Lambda$-modules.
Then by \cite[Corollary 3.5]{Sim5}, $d_3(M)$ is a dimension sequence. It follows that by \cite[Lemma 3.1]{Sim5} and \cite[Proposition $2^{'}$]{DRS},  ${\rm r.dim}~M = {\rm l.dim}~M =1$. \\
Conversely, assume that ${\rm r.dim}~M = {\rm l.dim}~M =1$. Let $M = xG$ for some $0 \neq x \in M$. Define a division ring embedding $\alpha: F \rightarrow G$ by the formula $fx=x\alpha(f)$ for any $f \in F$. So we have an $F$-$G$-isomorphism $\varphi: M \rightarrow {_{\alpha(F)}G_G}$ by the formula $\varphi(xg)=g$. Since ${\rm l.dim}~M =1$, we can assume that $M=Fx$. Define a division ring embedding $\beta: G \rightarrow F$ by the formula $xg=\beta(g)x$. Thus we have an $F$-$G$-isomorphism $\psi: M \rightarrow {_FF_{\beta(G)}}$ by the formula $\psi(fx)=f$. Therefore ${_{\alpha(F)}G_G} \cong M \cong {_FF_{\beta(G)}}$ as $F$-$G$-bimodules. It follows that $M^R \cong {_GG_{\alpha(F)}}$ and $M^L \cong {_{\beta(G)}F_F}$ as $G$-$F$-bimodules. Clearly  $\beta = \alpha^{-1}$ and the isomorphism $\psi \varphi^{-1}$ is just equal $\beta$. Notice that the same $\beta$ yields also the isomorphism ${_GG_{\alpha(F)}} \cong {_{\beta(G)}F_F}$, so $M^{L} \cong M^{R}$ as $G$-$F$-bimodules. Therefore by \cite[Proposition 2.6]{DR1}, there exist precisely $3$ pairwise non-isomorphic finitely generated indecomposable right $\Lambda$-modules.
\end{proof}

Let $F$ and $G$ be division rings and $M$ be an $F$-$G$-bimodule. We say that the $F$-$G$-bimodule $M$ is {\it trivial} if there exist precisely $3$ pairwise non-isomorphic finitely generated indecomposable right $\tiny{\begin{pmatrix}
F & M\\
0& G
\end{pmatrix}}$-modules.

Let $\Lambda$ be a basic hereditary ring and suppose that the number of vertices of the valued quiver $(\Gamma_{\Lambda}, \mathbf{d})$ of $\Lambda$ is the natural number $n$. A vector $x \in {\Bbb{N}}^n$ is called {\it branch system} of the Coxeter valued quiver $(\mathscr{C}_{\Lambda}, \mathbf{m})$ of $\Lambda$ if there exists an admissible sequence $k_1, \ldots, k_t$ in the valued quiver $(\Gamma_{\Lambda}, \mathbf{d})$ of $\Lambda$ such that ${s_t^+} \cdots {s_{1}^+}x=e_j$ for some $1 \leq j \leq n$. This generalizes the positive part of the usual rank 2 root systems (see \cite{DRS}). The following well-known result is playing a key role in our study in this paper.

\begin{Pro}\label{r2}
Let $\Lambda$ be a basic hereditary artinian ring. Then the underlying Coxeter valued  graph of the Coxeter valued quiver $(\mathscr{C}_{\Lambda}, \mathbf{m})$ of $\Lambda$ is a disjoint union of the Coxeter diagrams presented in Table B if and only if $\Lambda$ is a representation-finite ring. Moreover in this case
\begin{itemize}
\item[$(i)$] There exists a bijection between the indecomposable finite-dimensional representations of $\mathscr{M}_{\Lambda}$ and the branch system of the Coxeter valued quiver $(\mathscr{C}_{\Lambda}, \mathbf{m})$ of $\Lambda$.
\item[$(ii)$] The ring $\Lambda$ is isomorphic with the tensor ring $R_{{\mathscr{M}}_{\Lambda}}$ of the species ${\mathscr{M}}_{\Lambda}$ of $\Lambda$.
\item[$(iii)$] The functor $\mathrm{F}: {\rm rep}({\mathscr{M}}_{\Lambda}) \longrightarrow$ mod-$\Lambda$ in \eqref{f} with $\Lambda=R_{{\mathscr{M}}_{\Lambda}}$ is an equivalence of categories.
\item[$(iv)$] If the underlying Coxeter valued graph of the Coxeter valued quiver $(\mathscr{C}_{\Lambda}, \mathbf{m})$ of $\Lambda$ is  one of the Coxeter diagrams $\mathbb{A}_n$, $\mathbb{B}_n$, $\mathbb{D}_n$, $\mathbb{E}_6$, $\mathbb{E}_7$, $\mathbb{E}_8$ and $\mathbb{F}_4$ presented in Table B, then the underlying valued graph of the valued quiver $(\Gamma_{\Lambda}, \mathbf{d})$ of $\Lambda$ is one of the valued Dynkin diagrams presented in Table C.
\end{itemize}
\end{Pro}
\begin{center}
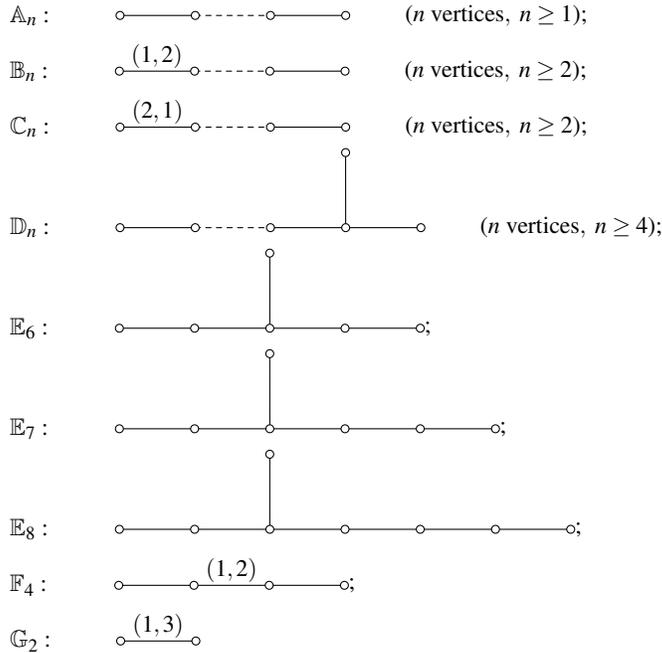
\textbf{Table C.}  Valued Dynkin diagrams (consult also \cite{simson})
\end{center}

\definecolor{qqqqff}{rgb}{0.,0.,1.}
	\begin{tikzpicture}[line cap=round,line join=round,>=triangle 45,x=1.0cm,y=1.0cm]
	\draw (2.,9.)-- (3.,9.);
	\draw [dash pattern=on 2pt off 2pt](3.,9.)--(4.,9.);
	\draw (4.,9.)-- (5.,9.);
		\begin{scriptsize}
	\draw[color=black] (.8,9.) node {$\mathbb{A}_n:$};
	\draw [fill=white] (2.,9.) circle (1.5pt);
	\draw [fill=white] (3.,9.) circle (1.5pt);
	\draw [fill=white] (4.,9.) circle (1.5pt);
	\draw [fill=white] (5.,9.) circle (1.5pt);
	\draw[color=black] (7.,9.) node {($n ~~{\rm vertices},~ n\geq 1$);};
	\end{scriptsize}
	\end{tikzpicture}

\definecolor{qqqqff}{rgb}{0.,0.,1.}
	\begin{tikzpicture}[line cap=round,line join=round,>=triangle 45,x=1.0cm,y=1.0cm]
	\draw (2.,9.)-- (3.,9.);
	\draw [dash pattern=on 2pt off 2pt](3.,9.)--(4.,9.);
	\draw (4.,9.)-- (5.,9.);
		\begin{scriptsize}
	\draw[color=black] (.8,9.) node {$\mathbb{B}_n:$};
	\draw [fill=white] (2.,9.) circle (1.5pt);
	\draw[color=black] (2.5,9.2) node {$(1,2)$};
	\draw [fill=white] (3.,9.) circle (1.5pt);
	\draw [fill=white] (4.,9.) circle (1.5pt);
	\draw [fill=white] (5.,9.) circle (1.5pt);
	\draw[color=black] (7.,9.) node {($n ~~{\rm vertices},~ n\geq 2$);};
	\end{scriptsize}
	\end{tikzpicture}

	\definecolor{qqqqff}{rgb}{0.,0.,1.}
	\begin{tikzpicture}[line cap=round,line join=round,>=triangle 45,x=1.0cm,y=1.0cm]
	\draw (2.,9.)-- (3.,9.);
	\draw [dash pattern=on 2pt off 2pt](3.,9.)--(4.,9.);
	\draw (4.,9.)-- (5.,9.);
		\begin{scriptsize}
	\draw[color=black] (.8,9.) node {$\mathbb{C}_n:$};
	\draw [fill=white] (2.,9.) circle (1.5pt);
	\draw[color=black] (2.5,9.2) node {$(2,1)$};
	\draw [fill=white] (3.,9.) circle (1.5pt);
	\draw [fill=white] (4.,9.) circle (1.5pt);
	\draw [fill=white] (5.,9.) circle (1.5pt);
	\draw[color=black] (7.,9.) node {($n ~~{\rm vertices},~ n\geq 2$);};
	\end{scriptsize}
	\end{tikzpicture}

	\definecolor{qqqqff}{rgb}{0.,0.,1.}
	\begin{tikzpicture}[line cap=round,line join=round,>=triangle 45,x=1.0cm,y=1.0cm]
	\draw (2.,9.)-- (3.,9.);
	\draw [dash pattern=on 2pt off 2pt](3.,9.)--(4.,9.);
	\draw (4.,9.)-- (5.,9.);
	\draw (5.,9.)-- (6.,9.);
	\draw (5.,9.)-- (6.,10.9999999999987,2.6);
		\begin{scriptsize}
	\draw[color=black] (.8,9.) node {$\mathbb{D}_n:$};
	\draw [fill=white] (2.,9.) circle (1.5pt);
	\draw [fill=white] (3.,9.) circle (1.5pt);
	\draw [fill=white] (4.,9.) circle (1.5pt);
	\draw [fill=white] (5.,9.) circle (1.5pt);
	\draw [fill=white] (6.,9.) circle (1.5pt);
	\draw [fill=white] (6.,10.9999999999987,2.6) circle (1.5pt);
	\draw[color=black] (8.,9.) node {($n ~~{\rm vertices},~ n\geq 4$);};
	\end{scriptsize}
	\end{tikzpicture}

	\definecolor{qqqqff}{rgb}{0.,0.,1.}
	\begin{tikzpicture}[line cap=round,line join=round,>=triangle 45,x=1.0cm,y=1.0cm]
	\draw (2.,9.)-- (3.,9.);
	\draw (3.,9.)-- (4.,9.);
	\draw (4.,9.)-- (5.,9.);
	\draw (4.,9.)-- (5.,10.9999999999987,2.6);
	\draw (5.,9.)-- (6.,9.);
		\begin{scriptsize}
	\draw[color=black] (.8,9.) node {$\mathbb{E}_6:$};
	\draw [fill=white] (2.,9.) circle (1.5pt);
	\draw [fill=white] (3.,9.) circle (1.5pt);
	\draw [fill=white] (4.,9.) circle (1.5pt);
	\draw [fill=white] (5.,9.) circle (1.5pt);
	\draw [fill=white] (5.,10.9999999999987,2.6) circle (1.5pt);
	\draw [fill=white] (6.,9.) circle (1.5pt);
	\draw [color=black] (6.1,9.) node {;};
	\end{scriptsize}
	\end{tikzpicture}

	\definecolor{qqqqff}{rgb}{0.,0.,1.}
	\begin{tikzpicture}[line cap=round,line join=round,>=triangle 45,x=1.0cm,y=1.0cm]
	\draw (2.,9.)-- (3.,9.);
	\draw (3.,9.)-- (4.,9.);
	\draw (4.,9.)-- (5.,9.);
	\draw (4.,9.)-- (5.,10.9999999999987,2.6);
	\draw (5.,9.)-- (6.,9.);
	\draw (6.,9.)-- (7.,9.);
		\begin{scriptsize}
	\draw[color=black] (.8,9.) node {$\mathbb{E}_7:$};
	\draw [fill=white] (2.,9.) circle (1.5pt);
	\draw [fill=white] (3.,9.) circle (1.5pt);
	\draw [fill=white] (4.,9.) circle (1.5pt);
	\draw [fill=white] (5.,9.) circle (1.5pt);
	\draw [fill=white] (5.,10.9999999999987,2.6) circle (1.5pt);
	\draw [fill=white] (6.,9.) circle (1.5pt);
	\draw [fill=white] (7.,9.) circle (1.5pt);
	\draw [color=black] (7.1,9.) node {;};
	\end{scriptsize}
	\end{tikzpicture}

	\definecolor{qqqqff}{rgb}{0.,0.,1.}
	\begin{tikzpicture}[line cap=round,line join=round,>=triangle 45,x=1.0cm,y=1.0cm]
	\draw (2.,9.)-- (3.,9.);
	\draw (3.,9.)-- (4.,9.);
	\draw (4.,9.)-- (5.,9.);
	\draw (4.,9.)-- (5.,10.9999999999987,2.6);
	\draw (5.,9.)-- (6.,9.);
	\draw (6.,9.)-- (7.,9.);
	\draw (7.,9.)-- (8.,9.);
		\begin{scriptsize}
	\draw[color=black] (.8,9.) node {$\mathbb{E}_8:$};
	\draw [fill=white] (2.,9.) circle (1.5pt);
	\draw [fill=white] (3.,9.) circle (1.5pt);
	\draw [fill=white] (4.,9.) circle (1.5pt);
	\draw [fill=white] (5.,9.) circle (1.5pt);
	\draw [fill=white] (5.,10.9999999999987,2.6) circle (1.5pt);
	\draw [fill=white] (6.,9.) circle (1.5pt);
	\draw [fill=white] (7.,9.) circle (1.5pt);
	\draw [fill=white] (8.,9.) circle (1.5pt);
	\draw [color=black] (8.1,9.) node {;};
	\end{scriptsize}
	\end{tikzpicture}
	
		\definecolor{qqqqff}{rgb}{0.,0.,1.}
	\begin{tikzpicture}[line cap=round,line join=round,>=triangle 45,x=1.0cm,y=1.0cm]
	\draw (2.,9.)-- (3.,9.);
	\draw (3.,9.)-- (4.,9.);
	\draw (4.,9.)-- (5.,9.);
		\begin{scriptsize}
	\draw[color=black] (.8,9.) node {$\mathbb{F}_4:$};
	\draw [fill=white] (2.,9.) circle (1.5pt);
	\draw [fill=white] (3.,9.) circle (1.5pt);
	\draw[color=black] (3.5,9.2) node {$(1,2)$};
	\draw [fill=white] (4.,9.) circle (1.5pt);
	\draw [fill=white] (5.,9.) circle (1.5pt);
	\draw [color=black] (5.1,9.) node {;};
	\end{scriptsize}
	\end{tikzpicture}
	
		\definecolor{qqqqff}{rgb}{0.,0.,1.}
	\begin{tikzpicture}[line cap=round,line join=round,>=triangle 45,x=1.0cm,y=1.0cm]
	\draw (2.,9.)-- (3.,9.);
		\begin{scriptsize}
	\draw[color=black] (.8,9.) node {$\mathbb{G}_2:$};
	\draw [fill=white] (2.,9.) circle (1.5pt);
	\draw[color=black] (2.5,9.2) node {$(1,3)$};
	\draw [fill=white] (3.,9.) circle (1.5pt);
	
	\end{scriptsize}
	\end{tikzpicture}
\begin{proof}
It follows from \cite[Theorem 2]{DRS}.\\
 $(i)$. It follows from Theorem \ref{DS1}, see also \cite[Theorem 2]{DRS}.\\
$(ii)$. The proof in case when $\Lambda$ is a basic hereditary finite-dimensional algebra over a field is given in \cite{2}. The proof in the general case follows from Theorem 4.5 and Corollary 4.6 in \cite{sims} together with Lemma 3.3 in \cite{simson}, see also \cite[Theorem 3]{Drozd}.\\
$(iii)$. By Lemma 3.1  and Corollary 3.5 in \cite{Sim5} and Proposition $2^{'}$ in \cite{DRS}, the species ${\mathscr{M}}_{\Lambda}$ of $\Lambda$ is finite-dimensional. Since the valued quiver $(\Gamma_{\Lambda}, \mathbf{d})$  of $\Lambda$  is acyclic and connected, by \cite[Proposition 1.1]{Sim3}, the functor $\mathrm{F}: {\rm rep}({\mathscr{M}}_{\Lambda}) \longrightarrow$ mod-$\Lambda$ in (\ref{f}) with $\Lambda=R_{{\mathscr{M}}_{\Lambda}}$ is an equivalence of categories.\\
$(iv)$. Assume that the underlying Coxeter valued graph of the Coxeter valued quiver $(\mathscr{C}_{\Lambda}, \mathbf{m})$ of $\Lambda$ is  one of the Coxeter diagrams $\mathbb{A}_n$, $\mathbb{D}_n$, $\mathbb{E}_6$, $\mathbb{E}_7$ and $\mathbb{E}_8$ presented in Table B. Then by Lemma \ref{2.12.1}, the underlying valued graph of the valued quiver $(\Gamma_{\Lambda}, \mathbf{d})$ of $\Lambda$ is one of the valued Dynkin diagrams $\mathbb{A}_n$, $\mathbb{D}_n$, $\mathbb{E}_6$, $\mathbb{E}_7$ and $\mathbb{E}_8$ presented in Table C, respectively. If the underlying Coxeter valued graph of the Coxeter valued quiver $(\mathscr{C}_{\Lambda}, \mathbf{m})$ of $\Lambda$ is the Coxeter diagram $\mathbb{B}_n$ presented in Table B, then by Lemma 3.1  and Corollary 3.5 in \cite{Sim5} and Proposition $2^{'}$ in \cite{DRS}, the underlying valued graph of the valued quiver $(\Gamma_{\Lambda}, \mathbf{d})$ of $\Lambda$ is one of the valued Dynkin diagrams $\mathbb{B}_n$ and $\mathbb{C}_n$ presented in Table C. If the underlying Coxeter valued graph of the Coxeter valued quiver $(\mathscr{C}_{\Lambda}, \mathbf{m})$ of $\Lambda$ is the Coxeter diagram $\mathbb{F}_4$ presented in Table B, by the same arguments the underlying  valued graph of the valued quiver $(\Gamma_{\Lambda}, \mathbf{d})$ of $\Lambda$ is the valued Dynkin diagram $\mathbb{F}_4$ presented in Table C.
\end{proof}

Let $\Lambda$ be a basic hereditary ring such that the underlying valued graph of the valued quiver $(\Gamma_{\Lambda}, \mathbf{d})$ of $\Lambda$ is one of the valued Dynkin diagrams presented in Table C.  Therefore there exist non-zero natural numbers $f_i$ satisfying $d_{ij}f_j = d_{ji}f_i$ for each vertex $i , j$ of $\Gamma_{\Lambda}$. Assume that the number of vertices of the valued quiver $(\Gamma_{\Lambda}, \mathbf{d})$ of $\Lambda$ is the natural number $n$. We denote by ${\Bbb{Q}}^n$ the vector space of all $x=(x_1, \ldots, x_n)$ over the field of rational numbers. Then we define a symmetric positive definite bilinear form $B: {\Bbb{Q}}^n \times {\Bbb{Q}}^n \rightarrow {\Bbb{Q}}$ as follows. For each $x, y \in {\Bbb{Q}}^n$, $B(x, y) = \sum_i f_ix_iy_i - \frac{1}{2} \sum_{i,j}d_{ij}f_jx_iy_j$.  For each vertex $k$ of $\Gamma_{\Lambda}$, we have a reflection $s_k$, where $s_k: {\Bbb{Q}}^n \rightarrow {\Bbb{Q}}^n$ is a linear transformation given by $s_kx= x - (2 B(x,e_k)/B(e_k,e_k))e_k$. A group of all linear transformations of ${\Bbb{Q}}^n$ generated by the reflections $s_k$, $k$ is a vertex of $\Gamma_{\Lambda}$, is called {\it Weyl group} and is denoted by ${\mathcal{W}}$. It is well known that the set $R=\lbrace x \in {\Bbb{Q}}^n ~ | ~ x=we_k ~{\rm for~ some}~ w \in {\mathcal{W}}~ {\rm and}~ k~{\rm is ~a~vertex~of~}\Gamma_{\Lambda} \rbrace$ is a reduced root system such that the set $\lbrace e_k ~|~ k~{\rm is ~a~vertex~of~} \Gamma_{\Lambda} \rbrace$ is a base for $R$ see \cite{DR1} and \cite{Bourbaki, hum}.  If the underlying Coxeter valued graph of the Coxeter valued quiver $(\mathscr{C}_{\Lambda}, \mathbf{m})$ of $\Lambda$ is one of the Coxeter diagrams $\mathbb{A}_n$, $\mathbb{B}_n$, $\mathbb{D}_n$, $\mathbb{E}_6$, $\mathbb{E}_7$, $\mathbb{E}_8$ and $\mathbb{F}_4$  presented in Table B, then by Theorem \ref{DS1}, \cite[Theorem 6.5]{ringel} and \cite[Theorem 1]{DRS} (see also \cite[Proposition 1.9]{DR1}), there exists a bijection between the isomorphism classes of finite-dimensional indecomposable representations of ${\mathscr{M}}_{\Lambda}$ and the positive roots of $(\Gamma_{\Lambda}, \mathbf{d})$. Now we assume that the underlying Coxeter valued graph of the Coxeter valued quiver $(\mathscr{C}_{\Lambda}, \mathbf{m})$ of $\Lambda$ is the Coxeter diagram $\mathbb{G}_2$ presented in Table B, then the category mod-$\Lambda$ has exactly 6 non-isomorphic indecomposable modules. By \cite[Corollary 3.5]{Sim5}, $d_6(M)$ is a dimension sequence. Therefore by \cite[Proposition $2^{'}$]{DRS}, $d_6(M)$ is one of the following sequences up to cyclic permutation and reversion:
$$(1,2,2,2,1,4), (1,2,3,1,2,3), (1,3,1,3,1,3)$$
If $d_6(M)=(1,3,1,3,1,3)$, then by Lemma 3.1 and Corollary 3.5 in \cite{Sim5}, the underlying valued graph of the valued quiver $(\Gamma_{\Lambda}, \mathbf{d})$ of $\Lambda$ is $\mathbb{G}_2$ presented in Table C. Thus by Proposition \ref{r2} and \cite[Proposition 1.9]{DR1}, the branch system of the $(\mathscr{C}_{\Lambda}, \mathbf{m})$ of $\Lambda$ is exactly the positive roots of $(\Gamma_{\Lambda}, \mathbf{d})$. But if $d_6(M)$ is one of the sequences (up to cyclic permutation and reversion) $(1,2,2,2,1,4)$ and $(1,2,3,1,2,3)$, then the underlying valued graph of the valued quiver $(\Gamma_{\Lambda}, \mathbf{d})$ of $\Lambda$ is one of the valued diagrams
\begin{center}
\dynkin[Coxeter, edge length=1.0cm, gonality={(1,2)}]{G}{2}, \dynkin[Coxeter, edge length=1.0cm, gonality={(2,2)}]{G}{2}, \dynkin[Coxeter, edge length=1.0cm, gonality={(1,4)}]{G}{2}, \dynkin[Coxeter, edge length=1.0cm, gonality={(2,3)}]{G}{2}, \dynkin[Coxeter, edge length=1.0cm, gonality={(2,1)}]{G}{2}, \dynkin[Coxeter, edge length=1.0cm, gonality={(4,1)}]{G}{2}, \dynkin[Coxeter, edge length=1.0cm, gonality={(3,1)}]{G}{2}
\end{center}
It follows that by Proposition \ref{r2}, the branch system of its Coxeter valued quiver is different of the positive roots of the corresponding Dynkin diagram. Therefore in Section 5, by using reflection functors,  we study the K\"othe property for basic hereditary rings $\Lambda$ of the Dynkin type $\mathbb{G}_2$ and of the Coxeter types $\mathbb{H}_3$, $\mathbb{H}_4$ and $\mathbb{I}_2(p)$ with $p=5$ or $7 \leq p < \infty$.

\section{Right K\"othe rings of Dynkin type}

Following \cite{Ringel2}, we say that a finitely generated indecomposable right $\Lambda$-module $M$ is {\it multiplicity-free top} if composition factors of top$(M)$ are pairwise non-isomorphic. The species ${\mathscr{M}}_{\Lambda}$ has {\it the multiplicity-free top property} if every finite-dimensional indecomposable representation of ${\mathscr{M}}_{\Lambda}$ is multiplicity-free top. \\

We start this section with the following fact that is frequently used in our study of right K\"othe rings of Dynkin type.

\begin{Pro}\label{r3}
The following two conditions are equivalent for a basic ring $\Lambda$.
\begin{itemize}
\item[$(i)$] $\Lambda$ is a right K\"othe ring.
\item[$(ii)$] $\Lambda$ is artinian and every indecomposable right $\Lambda$-module of finite length is multiplicity-free top.
\end{itemize}
If, in addition, $\Lambda$ is hereditary then $(i)$ is equivalent with the following statement:
\begin{itemize}
\item[$(iii)$] $\Lambda$ is a representation-finite ring and  the species ${\mathscr{M}}_{\Lambda}$ has the multiplicity-free top property.
\end{itemize}
\end{Pro}
\begin{proof}
The equivalence of $(i) \Leftrightarrow (ii)$ is a consequence of \cite[Corollary 3.3]{FN}.\\
$(ii) \Leftrightarrow (iii)$ Since $\Lambda$ is assumed to be hereditary, the equivalence of $(ii)$ and $(iii)$ follows from Propositions \ref{5} and \ref{r2}, and the equivalence $(i) \Leftrightarrow (ii)$ proved earlier.
\end{proof}

\begin{Pro}\label{2.12}
If $\Lambda$ is a basic hereditary artinian ring of the Dynkin type $\mathbb{A}_n$, then $\Lambda$ is a right K\"othe ring.
\end{Pro}
\begin{proof}
Assume that $\Lambda$ is a basic hereditary artinian ring of the Dynkin type $\mathbb{A}_n$. Then there exists a bijection between the finite-dimensional indecomposable representations of ${\mathscr{M}}_{\Lambda}$ and the positive roots of $(\Gamma_{\Lambda}, \mathbf{d})$.  Since by \cite[Page 265]{Bourbaki}, the positive roots of $(\Gamma_{\Lambda}, \mathbf{d})$ are $\sum_{i \leq k < j}e_k $, where $1 \leq i < j \leq n+1$,   by Proposition \ref{4}, every finite-dimensional indecomposable representation of the species ${\mathscr{M}}_{\Lambda}$ is multiplicity-free top. Therefore by Proposition \ref{r3}, $\Lambda$ is a right K\"othe ring.
\end{proof}

\begin{Pro}\label{2.1}
Let $\Lambda$ be a basic hereditary ring of the Dynkin type $\mathbb{D}_n$. Then $\Lambda$ is a right K\"othe ring if and only if $\Lambda$ is an artinian ring such that the following conditions hold:
 \begin{center}
 \definecolor{qqqqff}{rgb}{0.,0.,1.}
	\begin{tikzpicture}[line cap=round,line join=round,>=triangle 45,x=1.0cm,y=1.0cm]
	\draw (2.,9.)-- (3.,9.);
	\draw [dash pattern=on 2pt off 2pt](3.,9.)--(4.,9.);
	\draw (4.,9.)-- (5.,9.);
	\draw (5.,9.)-- (6.,9.);
	\draw (5.,9.)-- (6.,10.9999999999987,2.6);
		\begin{scriptsize}
	\draw[color=black] (.8,9.) node {$\mathbb{D}_n:$};
	\draw [fill=white] (2.,9.) circle (1.5pt);
	\draw[color=black] (2.,8.7) node {$1$};
	\draw [fill=white] (3.,9.) circle (1.5pt);
	\draw[color=black] (3.,8.7) node {$2$};
	\draw [fill=white] (4.,9.) circle (1.5pt);
	\draw[color=black] (4.,8.7) node {$n-3$};
	\draw [fill=white] (5.,9.) circle (1.5pt);
	\draw[color=black] (5.,8.7) node {$n-2$};
	\draw [fill=white] (6.,9.) circle (1.5pt);
	\draw[color=black] (6.,8.7) node {$n-1$};
	\draw [fill=white] (6.,10.9999999999987,2.6) circle (1.5pt);
	\draw[color=black] (6.2,10.9999999999987,2.6) node {$n$};
	
	\end{scriptsize}
	\end{tikzpicture}
 \end{center}
  \begin{itemize}
\item[$(a)$] $|{(n-2)}^{+}| \leq 2$;
\item[$(b)$] For each $i\leq n-3$, there exists at most one arrow with the source $i$.
\end{itemize}
 \end{Pro}
 \begin{proof}
$(\Rightarrow).$ Assume that $|{(n-2)}^{+}| > 2$. Since by Proposition \ref{r2}, the underlying valued graph of the valued quiver $(\Gamma_{\Lambda}, \mathbf{d})$ of $\Lambda$ is the valued Dynkin diagram $\mathbb{D}_n$ presented in Table C, by using \cite[Page 271]{Bourbaki}, there exists a finite-dimensional indecomposable representation $X$ of ${\mathscr{M}}_{\Lambda}$ with the dimension vector $\textbf{dim}~X = {\bsm\hspace{0.6cm}1\\0\cdots 0121\esm}$. Therefore by Proposition \ref{4}, there exists a finite-dimensional indecomposable representation $X$ of ${\mathscr{M}}_{\Lambda}$ such that ${\rm top}(X) \cong {\bf{F}}_{n-2} \oplus {\bf{F}}_{n-2}$, which is a contradiction by Proposition \ref{r3}. Now, we show that for each $i \leq n-3$, there exists at most one arrow with the source $i$. Assume that there exists $i \leq n-3$ such that $|i^{+}|=2$. Since by using \cite[Page 271]{Bourbaki}, there exists a finite-dimensional indecomposable representation $Y$ of ${\mathscr{M}}_{\Lambda}$ with the dimension vector $\textbf{dim}~Y = {\bsm\hspace{1.3cm}1\\0\cdots012\cdots 2221\esm}$, where the first $2$ is in the i-th coordinate, by Proposition \ref{4}, there exists a finite-dimensional indecomposable representation $Y$ of ${\mathscr{M}}_{\Lambda}$ which is not multiplicity-free top. It follows that by Proposition \ref{r3}, $\Lambda$ is not right K\"othe, which is a contradiction.\\
$(\Leftarrow).$ Assume that $X=(X_i, {_j{\varphi}_i})$ is a finite-dimensional indecomposable representation of ${\mathscr{M}}_{\Lambda}$ and  there exists  $1 \leq t \leq n$ such that ${{\bf{F}}_{t}} \oplus {{\bf{F}}_{t}} \subseteq {\rm top}(X)$. Since by Proposition \ref{r2}, the underlying valued graph of the valued quiver $(\Gamma_{\Lambda}, \mathbf{d})$ of $\Lambda$ is the valued Dynkin diagram $\mathbb{D}_n$ presented in Table C, by \cite[Page 271]{Bourbaki}, the positive roots of $(\Gamma_{\Lambda}, \mathbf{d})$ can be expressed as combinations of simple roots as follows: \begin{itemize}
\item[$(1)$] $e_{i} + \cdots + e_{k}$  \hspace{1cm} $ 1\leq i \leq k \leq n-1$;
\item[$(2)$] $e_{i} + \cdots + e_{k} + 2e_{k+1} + \cdots + 2e_{n-2} + e_{n-1} + e_{n}$ \hspace{1cm} $1\leq i \leq k \leq n-3$;
\item[$(3)$] $e_{i} + \cdots + e_{n-2} + e_{n}$  \hspace{1cm}  $1\leq i \leq n-2$;
\item[$(4)$] $e_{i} + \cdots + e_{n-2} + e_{n-1} + e_{n}$  \hspace{1cm} $1\leq i \leq n-2$.
\end{itemize}
By Proposition \ref{4}, ${\rm dim}{(X_t)}_{F_t} = 2$ and the vertex $t$ is source. It follows that by $(b)$, $t = n-2$. Therefore $|{(n-2)}^+| =3$, which is a contradiction. Thus every finite-dimensional indecomposable representation of ${\mathscr{M}}_{\Lambda}$ is multiplicity-free top. Therefore by Propositions \ref{r2} and \ref{r3}, $\Lambda$ is a right K\"othe ring.
 \end{proof}
An {\it arm of length} $t$ is a pair $(Q', k)$ consisting of a quiver $Q'$ of type $\mathbb{A}_n$ presented in Table A and the vertex $k$ of $Q'$ which has at most one neighbor in $Q'$. We say that a quiver $Q$ has an arm $(Q', k)$ if $Q'$ is a full subquiver of $Q$ and there are no arrows between the vertices of $Q$ outside of $Q'$ and the vertices of $Q'$ different from $k$. Let $c$ be a vertex of a quiver $\Delta$ of tree type. We say that an arrow $\alpha: x \rightarrow y$ points to $c$ provided $y$ and $c$ belong to the same connected component of the quiver obtained from $\Delta$ by deleting $\alpha$. Let $\mathscr{M} = (F_i, {_iM_j})_{i, j \in I}$ be a species such that for any $i, j \in I$, $F_i \cong F_j$ as division rings and ${\rm r.dim}~{_iM_j} = {\rm l.dim}~{_iM_j} = 1$ and suppose that the valued quiver $(\Gamma, \mathbf{d})$ of $\mathscr{M}$ is of tree type. Let $X= \left( X_i,{_j\varphi_i} \right)$ be a finite-dimensional representation of ${\mathscr{M}}$ and $(Q', k)$ be an arm of $(\Gamma, \mathbf{d})$. We say that $X$ is {\it conical} on  $(Q', k)$ provided ${_j\varphi_i}$ is injective for any arrow $i \rightarrow j$ of $Q'$ which points to $k$ and for the remaining arrows $l \rightarrow t$ of $Q'$, ${_t\varphi_l}$ is surjective. The representation $X$ is said to be {\it thin}, provided ${\rm dim}{(X_i)}_{F_i} \leq 1$ for all vertices $i$. The support of $X$ is the set of vertices $i$ with $X_i \neq 0$ \cite{R}.\\

The following lemma is a generalization of \cite[Corollary 1.4]{R}.

\begin{Lem}\label{RH}
Let $F$ be a division ring and $\mathscr{M} = (F_i, {_iM_j})_{i, j \in I}$ be a species such that for each $i$ and $j$, $F_i \cong F$ as division rings and  ${_iM_j} \cong {_FF_F}$ as bimodule and the valued quiver $(\Gamma, \mathbf{d})$ of $\mathscr{M}$ is of tree type.  Let $(Q^{'}, k)$ be an arm of $(\Gamma, \mathbf{d})$. Then any finite-dimensional representation $X$ of $\mathscr{M}$ has a decomposition as $X = X^{'} \oplus X^{''}$, where $X^{'}$ is conical on $(Q^{'}, k)$ and the support of $X^{''}$ is contained in $Q^{'} \setminus \lbrace k\rbrace$. In particular, if $X$ is a finite-dimensional indecomposable representation of $\mathscr{M}$ such that $X_k \neq 0$, then $X$ is conical on $(Q^{'}, k)$.
\end{Lem}
\begin{proof}
Let $X= \left( X_i,{_j\varphi_i} \right)$ be a finite-dimensional representation of $\mathscr{M}$. Let $\mathscr{M}^{'}$ be a subspecies of $\mathscr{M}$ with the valued quiver $Q^{'}$. Then by \cite[Theorem]{DR1}, there exists a bijection between the finite-dimensional indecomposable representations of $\mathscr{M}^{'}$ and the positive roots of $Q^{'}$. Moreover $\mathscr{M}^{'}$ is representation-finite. Therefore by using \cite[Page 265]{Bourbaki}, every finite-dimensional representation of $\mathscr{M}^{'}$ is a direct sum of thin indecomposable representations. Thus the restriction $X|Q^{'}$ of $X$ to $Q^{'}$ is a direct sum of thin indecomposable representations $X(j)$ with $j \in J$. Let $J^{'}$ be the set of indices $j \in J$ such that $X(j)_k \neq 0$ and let $J^{''}$ be the set of indices $j \in J$ such that $X(j)_k = 0$. For each vertex $t$ of $Q^{'}$, we set $X^{'}_t = \bigoplus_{j \in J^{'}} X(j)_t$ and $X^{''}_t = \bigoplus_{j \in J^{''}} X(j)_t$. Moreover if $t$ is a vertex in $\Gamma \setminus Q^{'}$, then we set  $X^{'}_t = X_t$ and $X^{''}_t = 0$. Thus $X = X^{'} \oplus X^{''}$, where $X^{'}= \left( X^{'}_i,{_j\varphi_i} \right)$ and $X^{''}= \left( X^{''}_i,{_j\varphi_i} \right)$. Since the representations $X(j)$ with $j \in J^{'}$ are conical on $(Q^{'}, k)$,  $X^{'}$ is conical on $(Q^{'}, k)$. Also, since the representations $X(j)$ with $j \in J^{''}$ satisfy $X(j)_k = 0$,  the support of $X^{''}$ is contained in $Q^{'} \setminus \lbrace k \rbrace$. Therefore the proof complete.
\end{proof}

Let ${\mathscr{M}}$ be a species and assume that ${\mathscr{M}}^{'}$ is a subspecies of ${\mathscr{M}}$. If the species ${\mathscr{M}}$ has the multiplicity-free top property, then clearly the species ${\mathscr{M}}^{'}$ has the multiplicity-free top property. This fact is frequently used in the rest of the paper.

\begin{Pro}\label{3.4}
Let $\Lambda$ be a basic hereditary ring of the Dynkin type $\mathbb{E}_6$. Then $\Lambda$ is a right K\"othe ring if and only if $\Lambda$ is an artinian ring such that the following conditions hold:
 \begin{center}
 \definecolor{qqqqff}{rgb}{0.,0.,1.}
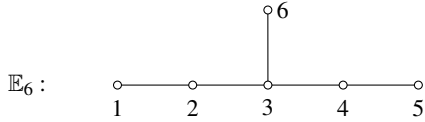
	\begin{tikzpicture}[line cap=round,line join=round,>=triangle 45,x=1.0cm,y=1.0cm]
	\draw (2.,9.)-- (3.,9.);
	\draw (3.,9.)-- (4.,9.);
	\draw (4.,9.)-- (5.,9.);
	\draw (4.,9.)-- (5.,10.9999999999987,2.6);
	\draw (5.,9.)-- (6.,9.);
		\begin{scriptsize}
	\draw[color=black] (.8,9.) node {$\mathbb{E}_6:$};
	\draw [fill=white] (2.,9.) circle (1.5pt);
	\draw[color=black] (2.,8.7) node {$1$};
	\draw [fill=white] (3.,9.) circle (1.5pt);
	\draw[color=black] (3.,8.7) node {$2$};
	\draw [fill=white] (4.,9.) circle (1.5pt);
	\draw[color=black] (4.,8.7) node {$3$};
	\draw [fill=white] (5.,9.) circle (1.5pt);
	\draw[color=black] (5.,8.7) node {$4$};
	\draw [fill=white] (5.,10.9999999999987,2.6) circle (1.5pt);
	\draw[color=black] (5.2,10.9999999999987,2.6) node {$6$};
	\draw [fill=white] (6.,9.) circle (1.5pt);
	\draw[color=black] (6.,8.7) node {$5$};
	\end{scriptsize}
	\end{tikzpicture}
 \end{center}
\begin{itemize}
\item[$(a)$]  $1 \leq |3^{+}| \leq 2$, $|4^{+}| \leq 1$ and $|2^{+}| \leq 1$;
\item[$(b)$] For each $y \in 3^{-}$, there exists at least one arrow with the target $y$.
\end{itemize}
 \end{Pro}
 \begin{proof}
 $(\Rightarrow).$ Assume that $\Lambda$ is a right K\"othe ring. Let ${\mathscr{M}}^{'}$ be a subspecies of ${\mathscr{M}_{\Lambda}}$ such that the underlying valued graph of the valued quiver $(\Gamma^{'}, \mathbf{d})$ of ${\mathscr{M}}^{'}$ is the valued Dynkin diagram $\mathbb{D}_5$ presented in Table C. Then by using Proposition \ref{2.1},  $|{3}^{+}| \leq 2$ and there exists at most one arrow with the source $2$ and one arrow with the source $4$. Consequently $|4^{+}| \leq 1$ and $|2^{+}| \leq 1$. If $|3^{+}|=0$, then by Proposition \ref{r2}, the Coxeter valued quiver $(\mathscr{C}_{\Lambda}, \mathbf{m})$ of $\Lambda$ has the orientation
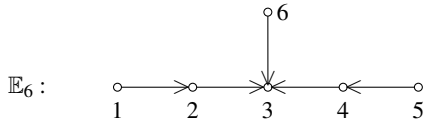
\begin{center}
	\definecolor{qqqqff}{rgb}{0.,0.,1.}
	\begin{tikzpicture}[line cap=round,line join=round,>=triangle 45,x=1.0cm,y=1.0cm]
	\draw (2.,9.)-- (3.,9.);
	\draw (3.,9.)-- (4.,9.);
	\draw (4.,9.)-- (5.,9.);
	\draw (4.,9.)-- (5.,10.9999999999987,2.6);
	\draw (5.,9.)-- (6.,9.);
		\begin{scriptsize}
	\draw[color=black] (.8,9.) node {$\mathbb{E}_6:$};
	\draw[color=black] (.8,9.) node {};
	\draw [fill=white] (2.,9.) circle (1.5pt);
	\draw[color=black] (2.,8.7) node {$1$};
	\draw[color=black] (2.85,9.) node {$>$};
	\draw [fill=white] (3.,9.) circle (1.5pt);
	\draw[color=black] (3.,8.7) node {$2$};
	\draw[color=black] (3.85,9.) node {$>$};
	\draw [fill=white] (4.,9.) circle (1.5pt);
	\draw[color=black] (4.,8.7) node {$3$};
	\draw[color=black] (4.14,9.) node {$<$};
	\draw [fill=white] (5.,9.) circle (1.5pt);
	\draw[color=black] (5.,8.7) node {$4$};
    \draw[color=black] (5.14,9.) node {$<$};
	\draw [fill=white] (5.,10.9999999999987,2.6) circle (1.5pt);
	\draw[color=black] (5.2,10.9999999999987,2.6) node {$6$};
	\draw[color=black] (4.,9.12) node {$\vee$};
	\draw [fill=white] (6.,9.) circle (1.5pt);
	\draw[color=black] (6.,8.7) node {$5$};
	\end{scriptsize}
	\end{tikzpicture}
\end{center}
 Hence by using \cite[Page 275]{Bourbaki}, there exists a finite-dimensional indecomposable representation $X$ of ${\mathscr{M}}_{\Lambda}$ with the dimension vector $\textbf{dim}~X = {\bsm\hspace{-.5mm}2\\12321\esm}$. It follows that by Proposition \ref{4}, there exists a finite-dimensional indecomposable representation $X$ of ${\mathscr{M}}_{\Lambda}$ such that ${{\bf{F}}_{6}} \oplus {{\bf{F}}_{6}} \subseteq {\rm top}(X)$ which is a contradiction, by Proposition \ref{r3}. Therefore $1 \leq |3^{+}| \leq 2$. Now we show that for each $y \in 3^{-}$, there exists at least one arrow with the target $y$. Assume that there exists $y \in 3^{-}$ such that there is no arrow with the target $y$. Since $|4^{+}| \leq 1$ and $|2^{+}| \leq 1$, $y=6$. Therefore by Proposition \ref{4}, ${{\bf{F}}_{6}} \oplus {{\bf{F}}_{6}} \subseteq {\rm top}(X)$, where $X$ is a finite-dimensional indecomposable representation of ${\mathscr{M}}_{\Lambda}$ with the dimension vector $\textbf{dim}~X = {\bsm\hspace{-.5mm}2\\12321\esm}$. Hence $X$ is not multiplicity-free top, which is a contradiction, by Proposition \ref{r3}. \\
$(\Leftarrow)$. Let ${\mathscr{M}_{\Lambda}}= (F_i, {_iM_j})_{i, j \in I}$ be the species of $\Lambda$. Assume that $X =(X_i, {_j{\varphi}_i})$ is a finite-dimensional indecomposable representation of ${\mathscr{M}_{\Lambda}}$ such that ${{\bf{F}}_{i}} \oplus {{\bf{F}}_{i}} \subseteq {\rm top}(X)$ for some vertex $i$ of $\Gamma_{\Lambda}$. Since by Lemma \ref{2.12.1}, there exists a division ring $F$ such that each $F_i\cong F$ as division rings, by using \cite[Page 275]{Bourbaki},  ${\rm dim}{(X_1)}_F \leq 1$,  ${\rm dim}{(X_2)}_F \leq 2$, ${\rm dim}{(X_3)}_F \leq 3$, ${\rm dim}{(X_4)}_F \leq 2$, ${\rm dim}{(X_5)}_F \leq 1$ and ${\rm dim}{(X_6)}_F \leq 2$. Thus by assumptions $(a)$ and $(b)$ and Proposition \ref{4},  $i=3$ and ${\rm dim}{(X_3)}_F=3$. It follows that by Lemma \ref{RH}, $X$ is one of the following representations
$$\hskip .5cm \xymatrix{
&& X_6 \ar @{<<-}[d]\\
X_1 \ar@{^{(}->}[r]& X_2 \ar@{^{(}->}[r]^{_3{\varphi}_2}& X_3\ar @{->>}[r]& X_4 \ar @{->>}[r]& X_5} \hskip .5cm ~ \hskip .5cm \xymatrix{
&& X_6 \ar @{<<-}[d]\\
X_1 \ar@{^{(}->}[r] & X_2 \ar@{^{(}->}[r]^{_3{\varphi}_2}& X_3 \ar@{->>}[r]& X_4 \ar@^{<-)}[r]& X_5}\hskip .5cm $$~  $$\hskip .5cm \xymatrix{
&& X_6 \ar @{<<-}[d]\\
X_1 \ar@{^{(}->}[r]& X_2 \ar@{<<-}[r]& X_3 \ar@^{<-)}[r]^{_3{\varphi}_4}& X_4 \ar@^{<-)}[r] & X_5} \hskip .5cm ~ \hskip .5cm \xymatrix{
&& X_6 \ar @{<<-}[d]\\
X_1 \ar @{<<-}[r]& X_2 \ar @{<<-}[r]& X_3\ar@^{<-)}[r]^{_3{\varphi}_4}& X_4 \ar@^{<-)}[r]& X_5}\hskip .5cm$$
where ${\rm dim}{({\rm Im}({_3{\varphi}_2}))}_F = 1$ and ${\rm dim}{({\rm Im}({_3{\varphi}_4}))}_F = 1$. Since there exists a bijection between the finite-dimensional indecomposable representations of $\mathscr{M}_{\Lambda}$ and the positive roots of $(\Gamma_{\Lambda}, \mathbf{d})$,
by using \cite[Page 275]{Bourbaki}, there is no any indecomposable representations of ${\mathscr{M}}_{\Lambda}$ with the dimension vectors ${\bsm\hspace{-.5mm}*\\ *13**\esm}$ or ${\bsm\hspace{-.5mm}*\\ **31*\esm}$. It follows that $X$ is not indecomposable, which is a contradiction. Thus every finite-dimensional indecomposable representation of ${\mathscr{M}}_{\Lambda}$ is multiplicity-free top. Therefore by Propositions \ref{r2} and  \ref{r3}, $\Lambda$ is a right K\"othe ring.
 \end{proof}

 \begin{Pro}\label{3.5}
Let $\Lambda$ be a basic hereditary ring of the Dynkin type $\mathbb{E}_7$. Then $\Lambda$ is a right K\"othe ring if and only if $\Lambda$ is an artinian ring such that the Coxeter valued quiver $(\mathscr{C}_{\Lambda}, \mathbf{m})$ of $\Lambda$ has the orientation
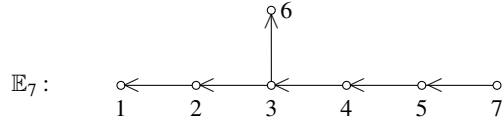
\begin{center}
	\definecolor{qqqqff}{rgb}{0.,0.,1.}
	\begin{tikzpicture}[line cap=round,line join=round,>=triangle 45,x=1.0cm,y=1.0cm]
	\draw (2.,9.)-- (3.,9.);
	\draw (3.,9.)-- (4.,9.);
	\draw (4.,9.)-- (5.,9.);
	\draw (4.,9.)-- (5.,10.9999999999987,2.6);
	\draw (5.,9.)-- (6.,9.);
	\draw (6.,9.)-- (7.,9.);
		\begin{scriptsize}
	\draw[color=black] (.8,9.) node {$\mathbb{E}_7:$};
	\draw[color=black] (.8,9.) node {};
	\draw [fill=white] (2.,9.) circle (1.5pt);
	\draw[color=black] (2.,8.7) node {$1$};
	\draw[color=black] (2.14,9.) node {$<$};
	\draw [fill=white] (3.,9.) circle (1.5pt);
	\draw[color=black] (3.,8.7) node {$2$};
	\draw[color=black] (3.14,9.) node {$<$};
	\draw [fill=white] (4.,9.) circle (1.5pt);
	\draw[color=black] (4.,8.7) node {$3$};
	\draw[color=black] (4.14,9.) node {$<$};
	\draw [fill=white] (5.,9.) circle (1.5pt);
	\draw[color=black] (5.,8.7) node {$4$};
    \draw[color=black] (5.14,9.) node {$<$};
	\draw [fill=white] (5.,10.9999999999987,2.6) circle (1.5pt);
	\draw[color=black] (5.2,10.9999999999987,2.6) node {$6$};
	\draw[color=black] (4.,9.85) node {$\wedge$};
	\draw [fill=white] (6.,9.) circle (1.5pt);
	\draw[color=black] (6.,8.7) node {$5$};
	\draw [fill=white] (7.,9.) circle (1.5pt);
	\draw[color=black] (7.,8.7) node {$7$};
    \draw[color=black] (6.17,9.) node {$<$};
	\end{scriptsize}
	\end{tikzpicture}
\end{center}
 \end{Pro}
 \begin{proof}
 $(\Rightarrow).$ Assume that $\Lambda$ is a right K\"othe ring. Then by the same argument as in the proof of Propositions \ref{2.1} and \ref{3.4}, we can see that $1 \leq |3^{+}| \leq 2$, $|4^{+}| \leq 1$,  $|5^{+}| \leq 1$ and $|2^{+}| \leq 1$ and for each $y \in 3^{-}$, there exists at least one arrow with the target $y$. Since by using \cite[Page 279]{Bourbaki},  there exists a finite-dimensional indecomposable representation $Y$ of ${\mathscr{M}}_{\Lambda}$ with the dimension vector $\textbf{dim}~Y = {\bsm\hspace{-.8mm}2\\234321\esm}$, by Propositions \ref{4} and \ref{r3}, $(\mathscr{C}_{\Lambda}, \mathbf{m})$ has the orientation \begin{center}
	\definecolor{qqqqff}{rgb}{0.,0.,1.}
	\begin{tikzpicture}[line cap=round,line join=round,>=triangle 45,x=1.0cm,y=1.0cm]
	\draw (2.,9.)-- (3.,9.);
	\draw (3.,9.)-- (4.,9.);
	\draw (4.,9.)-- (5.,9.);
	\draw (4.,9.)-- (5.,10.9999999999987,2.6);
	\draw (5.,9.)-- (6.,9.);
	\draw (6.,9.)-- (7.,9.);
		\begin{scriptsize}
	\draw[color=black] (.8,9.) node {$\mathbb{E}_7:$};
	\draw[color=black] (.8,9.) node {};
	\draw [fill=white] (2.,9.) circle (1.5pt);
	\draw[color=black] (2.,8.7) node {$1$};
	\draw[color=black] (2.14,9.) node {$<$};
	\draw [fill=white] (3.,9.) circle (1.5pt);
	\draw[color=black] (3.,8.7) node {$2$};
	\draw[color=black] (3.14,9.) node {$<$};
	\draw [fill=white] (4.,9.) circle (1.5pt);
	\draw[color=black] (4.,8.7) node {$3$};
	\draw[color=black] (4.14,9.) node {$<$};
	\draw [fill=white] (5.,9.) circle (1.5pt);
	\draw[color=black] (5.,8.7) node {$4$};
    \draw[color=black] (5.14,9.) node {$<$};
	\draw [fill=white] (5.,10.9999999999987,2.6) circle (1.5pt);
	\draw[color=black] (5.2,10.9999999999987,2.6) node {$6$};
	\draw[color=black] (4.,9.85) node {$\wedge$};
	\draw [fill=white] (6.,9.) circle (1.5pt);
	\draw[color=black] (6.,8.7) node {$5$};
	\draw [fill=white] (7.,9.) circle (1.5pt);
	\draw[color=black] (7.,8.7) node {$7$};
    \draw[color=black] (6.17,9.) node {$<$};
	\end{scriptsize}
	\end{tikzpicture}
\end{center}
$(\Leftarrow).$ Let ${\mathscr{M}}_{\Lambda} = (F_i, {_iM_j})_{i, j \in I}$ be the species of $\Lambda$. Assume that $X =(X_i, {_j{\psi}_i})$ is a finite-dimensional indecomposable representation of ${\mathscr{M}}_{\Lambda}$ such that ${{\bf{F}}_{j}} \oplus {{\bf{F}}_{j}} \subseteq {\rm top}(X)$ for some vertex $j$ of $\Gamma_{\Lambda}$. Since $\Lambda$ is of the Dynkin type $\mathbb{E}_7$, by Lemma \ref{2.12.1}, there exists a division ring $F$ such that for each $i \in I$, $F_i\cong F$ as division rings. If $X_3 = 0$, then by Proposition \ref{2.12}, $X$ is multiplicity-free top, which is a contradiction. Hence $X_3 \neq 0$. Therefore by Lemma \ref{RH}, $X$ is the representation $$\hskip .5cm \xymatrix{
&&X_6\ar @{<<-}[d]\\
 X_1\ar @{<<-}[r]& X_2\ar @{<<-}[r] & X_3 \ar@^{<-)}[r]& X_4\ar@^{<-)}[r]& X_5 \ar@^{<-)}[r]& X_7}\hskip .5cm$$
If $X_7=0$, by the proof of Proposition \ref{3.4}, $X$ is not indecomposable, which is a contradiction. Thus $X_7 \neq 0$. If ${\rm dim}{(X_5)}_F = 2$,  then by Lemma \ref{RH}, $\overline{X}= \overline{Z^{'}} \oplus \overline{Z^{''}}$, where
$$\xymatrix{ \overline{X}:= X_1 \ar@{<-}[r]& X_2 \oplus X_6 \ar@{<-}[r]& X_3 \ar@^{<-)}[r]& X_4 \ar@^{<-)}[r]& X_5 \ar@^{<-)}[r]& X_7,}$$ $$\xymatrix{\overline{Z^{'}}:= X^{'}_1 \ar@{<-}[r]& X^{'}_2 \oplus X^{'}_6 \ar@{<-}[r]& X^{'}_3 \ar@{<-}[r]& X^{'}_4 \ar@{<-}[r]& X^{'}_5 \ar@{<-}[r]& 0}$$   $$\xymatrix{\overline{Z^{''}}:=X^{''}_1 \ar@{<-}[r]& X^{''}_2 \oplus X^{''}_6 \ar@{<-}[r]& X^{''}_3 \ar@{<-}[r]& X^{''}_4 \ar@{<-}[r]& X^{''}_5 \ar@{<-}[r]^{_5{\psi}_7}& X_7}$$ and ${_5{\psi}_7}: X_7 \rightarrow X^{''}_5$ is an isomorphism. Thus $X = Z^{'} \oplus Z^{''}$, where $$ \xymatrix{
&&X^{'}_6\ar @{<-}[d]\\
 Z^{'}:=X^{'}_1\ar @{<-}[r]& X^{'}_2\ar @{<-}[r]& X^{'}_3 \ar@{<-}[r]& X^{'}_4 \ar@{<-}[r]& X^{'}_5 \ar@{<-}[r]& 0}$$
 $$\xymatrix{
&&X^{''}_6\ar @{<-}[d]\\
 Z^{''}:= X^{''}_1\ar @{<-}[r]& X^{''}_2\ar @{<-}[r]& X^{''}_3 \ar@{<-}[r] & X^{''}_4\ar@{<-}[r]& X^{''}_5 \ar@{<-}[r]^{{_5{\psi}_7}}& X_7.}$$ Since by using \cite[Page 278]{Bourbaki}, ${\rm dim}{(X_7)}_F = 1$, $X^{'}_5 \neq 0$. It follows that $X$ is not indecomposable, which is a contradiction. Therefore ${\rm dim}{(X_5)}_F \neq 2$.
 Since by using \cite[Page 278]{Bourbaki},  ${\rm dim}{(X_1)}_F \leq 2$, ${\rm dim}{(X_2)}_F \leq 3$, ${\rm dim}{(X_3)}_F \leq 4$, ${\rm dim}{(X_4)}_F \leq 3$, ${\rm dim}{(X_5)}_F \leq 2$, ${\rm dim}{(X_6)}_F \leq 2$ and ${\rm dim}{(X_7)}_F \leq 1$,  so ${\rm dim}{(X_7)}_F = {\rm dim}{(X_5)}_F =1$. Hence by Proposition \ref{4}, $j=3$ or $j=4$.  If $j=4$, then ${\rm dim}{(X_4)}_F=3$. Set  $$ \xymatrix{
&&X_6\ar @{<<-}[d]\\
 T:= X_1\ar @{<<-}[r]& X_2\ar @{<<-}[r] & X_3 \ar@^{<-)}[r] & X_4\ar@^{<-)}[r]^{_4{\psi}_5} & X_5}\hskip .5cm$$
Thus by Lemma \ref{RH}, $\overline{T} = \overline{T^{'}} \oplus \overline{T^{''}}$, where $$\xymatrix{\overline{T}:= X_1 \ar@{<-}[r]& X_2 \oplus X_6 \ar@{<-}[r]& X_3 \ar@^{<-)}[r]& X_4 \ar@^{<-)}[r]& X_5,}$$ $$\xymatrix{\overline{T^{'}}:= X^{'}_1 \ar@{<-}[r]& X^{'}_2 \oplus X^{'}_6 \ar@{<-}[r]& X^{'}_3 \ar@{<-}[r]& X^{'}_4 \ar@{<-}[r]& 0,}$$ $$\xymatrix{\overline{T^{''}}:= X^{''}_1 \ar@{<-}[r]& X^{''}_2 \oplus X^{''}_6 \ar@{<-}[r]& X^{''}_3 \ar@{<-}[r]& X^{''}_4 \ar@{<-}[r]^{{_4\psi_5}}& X_5}$$ and  ${_4{\psi}_5}: X_5 \rightarrow X^{''}_4$ is an isomorphism. So $X = T^{'} \oplus T^{''}$, where $$ \xymatrix{
&&X^{'}_6\ar @{<-}[d]\\
 T^{'} := X^{'}_1\ar @{<-}[r]& X^{'}_2\ar @{<-}[r]& X^{'}_3 \ar@{<-}[r]& X^{'}_4\ar@{<-}[r]& 0 \ar@{<-}[r]& 0}$$ $$\xymatrix{
&&X^{''}_6\ar @{<-}[d]\\
 T^{''}:= X^{''}_1\ar @{<-}[r]& X^{''}_2\ar @{<-}[r]& X^{''}_3 \ar@{<-}[r]& X^{''}_4\ar@{<-}[r]^{_4{\psi}_5} & X_5 \ar@{<-}[r]& X_7}$$
Since ${\rm dim}{(X_5)}_F =1$,  $X^{'}_4 \neq 0$ and  $T^{''} \neq 0$.  Consequently, $X$ is not indecomposable, which is a contradiction. Now assume that $j=3$. Therefore either ${\rm dim}{(X_3)}_F= 3$ and ${\rm dim}{({\rm Im}({_3{\psi}_4}))}_F =1$ or ${\rm dim}{(X_3)}_F= 4$ and $1 \leq {\rm dim}{({\rm Im}({_3{\psi}_4}))}_F \leq 2$. Set $$ \xymatrix{
&&X_6\ar @{<<-}[d]\\
 H:=X_1\ar @{<<-}[r]& X_2\ar @{<<-}[r]& X_3 \ar@^{<-)}[r] & X_4}\hskip .5cm$$ Then by Lemma \ref{RH}, $\overline{H} = \overline{H^{'}} \oplus \overline{H^{''}}$, where $$\xymatrix{\overline{H}:= X_1 \ar @{<-}[r]& X_2 \oplus X_6 \ar @{<-}[r]& X_3 \ar@^{<-)}[r]& X_4,}$$  $$\xymatrix{\overline{H^{'}} := X^{'}_1 \ar@{<-}[r]& X^{'}_2 \oplus X^{'}_6 \ar@{<-}[r]& X^{'}_3 \ar@{<-}[r]& 0}$$ $$\xymatrix{\overline{H^{''}} := X^{''}_1 \ar@{<-}[r]& X^{''}_2 \oplus X^{''}_6 \ar@{<-}[r]& X^{''}_3 \ar@{<-}[r]^{{_3{\psi}_4}}& X_4}$$ and ${_3{\psi}_4}:X_4 \rightarrow X^{''}_3$ is an isomorphism. It follows that $X = H^{'} \oplus H^{''}$, where $$ \xymatrix{
&&X^{'}_6\ar @{<-}[d]\\
 H^{'} := X^{'}_1\ar @{<-}[r]& X^{'}_2\ar @{<-}[r]& X^{'}_3 \ar@{<-}[r]& 0 \ar@{<-}[r]& 0 \ar@{<-}[r]& 0}$$  $$\xymatrix{
&&X^{''}_6\ar @{<-}[d]\\
 H^{''}:= X^{''}_1\ar @{<-}[r]& X^{''}_2\ar @{<-}[r]& X^{''}_3 \ar@{<-}[r]^{{_3{\psi}_4}} & X_4 \ar@{<-}[r] & X_5 \ar@{<-}[r] & X_7}$$
 and $X^{'}_3 \neq 0$. Consequently,  $X$ is not indecomposable, which is a contradiction. Hence every finite-dimensional indecomposable representation of  ${\mathscr{M}}_{\Lambda}$ is multiplicity-free top. Therefore by Propositions \ref{r2} and \ref{r3}, $\Lambda$ is a right K\"othe ring.
 \end{proof}

\begin{Lem} \label{m6}
If $\Lambda$ is a basic hereditary artinian ring of the Dynkin type $\mathbb{E}_8$, then $\Lambda$ is not a right K\"othe ring.
\end{Lem}
 \begin{proof}
Assume that $\Lambda$ is a right K\"othe ring. Then by Proposition \ref{r3}, every finite-dimensional indecomposable representation of ${\mathscr{M}}_{\Lambda}$ is multiplicity-free top. Thus by using Proposition \ref{3.5}, $(\mathscr{C}_{\Lambda}, \mathbf{m})$ of $\Lambda$ is one of the following quivers:
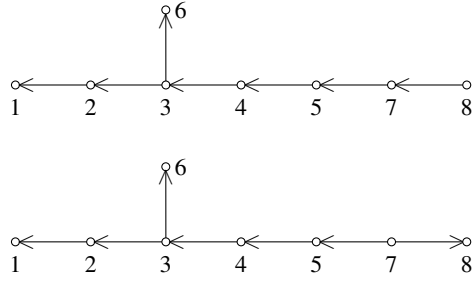
\begin{center}
	\definecolor{qqqqff}{rgb}{0.,0.,1.}
	\begin{tikzpicture}[line cap=round,line join=round,>=triangle 45,x=1.0cm,y=1.0cm]
	\draw (2.,9.)-- (3.,9.);
	\draw (3.,9.)-- (4.,9.);
	\draw (4.,9.)-- (5.,9.);
	\draw (4.,9.)-- (5.,10.9999999999987,2.6);
	\draw (5.,9.)-- (6.,9.);
	\draw (6.,9.)-- (7.,9.);
	\draw (7.,9.)-- (8.,9.);
		\begin{scriptsize}
	\draw [fill=white] (2.,9.) circle (1.5pt);
	\draw[color=black] (2.,8.7) node {$1$};
	\draw[color=black] (2.14,9.) node {$<$};
	\draw [fill=white] (3.,9.) circle (1.5pt);
	\draw[color=black] (3.,8.7) node {$2$};
	\draw[color=black] (3.14,9.) node {$<$};
	\draw [fill=white] (4.,9.) circle (1.5pt);
	\draw[color=black] (4.,8.7) node {$3$};
	\draw[color=black] (4.14,9.) node {$<$};
	\draw [fill=white] (5.,9.) circle (1.5pt);
	\draw[color=black] (5.,8.7) node {$4$};
	\draw[color=black] (5.14,9.) node {$<$};
	\draw [fill=white] (5.,10.9999999999987,2.6) circle (1.5pt);
	\draw[color=black] (5.2,10.9999999999987,2.6) node {$6$};
	\draw[color=black] (4.,9.85) node {$\wedge$};
	\draw [fill=white] (6.,9.) circle (1.5pt);
	\draw[color=black] (6.,8.7) node {$5$};
	\draw[color=black] (6.14,9.) node {$<$};
	\draw [fill=white] (7.,9.) circle (1.5pt);
	\draw[color=black] (7.,8.7) node {$7$};
	\draw[color=black] (7.14,9.) node {$<$};
	\draw [fill=white] (8.,9.) circle (1.5pt);
	\draw[color=black] (8.,8.7) node {$8$};
	\end{scriptsize}
	\end{tikzpicture}
\end{center}
\begin{center}
	\definecolor{qqqqff}{rgb}{0.,0.,1.}
	\begin{tikzpicture}[line cap=round,line join=round,>=triangle 45,x=1.0cm,y=1.0cm]
	\draw (2.,9.)-- (3.,9.);
	\draw (3.,9.)-- (4.,9.);
	\draw (4.,9.)-- (5.,9.);
	\draw (4.,9.)-- (5.,10.9999999999987,2.6);
	\draw (5.,9.)-- (6.,9.);
	\draw (6.,9.)-- (7.,9.);
	\draw (7.,9.)-- (8.,9.);
		\begin{scriptsize}
	\draw [fill=white] (2.,9.) circle (1.5pt);
	\draw[color=black] (2.,8.7) node {$1$};
	\draw[color=black] (2.14,9.) node {$<$};
	\draw [fill=white] (3.,9.) circle (1.5pt);
	\draw[color=black] (3.,8.7) node {$2$};
	\draw[color=black] (3.14,9.) node {$<$};
	\draw [fill=white] (4.,9.) circle (1.5pt);
	\draw[color=black] (4.,8.7) node {$3$};
	\draw[color=black] (4.14,9.) node {$<$};
	\draw [fill=white] (5.,9.) circle (1.5pt);
	\draw[color=black] (5.,8.7) node {$4$};
	\draw[color=black] (5.14,9.) node {$<$};
	\draw [fill=white] (5.,10.9999999999987,2.6) circle (1.5pt);
	\draw[color=black] (5.2,10.9999999999987,2.6) node {$6$};
	\draw[color=black] (4.,9.85) node {$\wedge$};
	\draw [fill=white] (6.,9.) circle (1.5pt);
	\draw[color=black] (6.,8.7) node {$5$};
	\draw[color=black] (6.14,9.) node {$<$};
	\draw [fill=white] (7.,9.) circle (1.5pt);
	\draw[color=black] (7.,8.7) node {$7$};
	\draw[color=black] (7.87,9.) node {$>$};
	\draw [fill=white] (8.,9.) circle (1.5pt);
	\draw[color=black] (8.,8.7) node {$8$};
	\end{scriptsize}
	\end{tikzpicture}
\end{center}
Since by using \cite[Page 284]{Bourbaki}, there exists a finite-dimensional indecomposable representation $Z$ of ${\mathscr{M}}_{\Lambda}$ with the dimension vector $\textbf{dim}~Z = {\bsm\hspace{-3mm}3\\2465432\esm}$, by Proposition \ref{4}, either ${\mathbf{F}}_8 \oplus {\mathbf{F}}_8 \subseteq {\rm top}(Z)$ or $\mathbf{F}_7 \oplus \mathbf{F}_7 \oplus \mathbf{F}_7 \subseteq {\rm top}(Z)$, which is a contradiction. Therefore $\Lambda$ is not right K\"othe.
 \end{proof}

\section{A characterization of right K\"othe hereditary rings}
In this section we give a characterization of right K\"othe hereditary rings.

\begin{The}\label{3.7}
Let $F$ and $G$ be division rings and $M$ be an $F$-$G$-bimodule. Then $\Lambda = \tiny{\begin{pmatrix}
F & M\\
0& G
\end{pmatrix}}$ is a right K\"othe ring if and only if there exists $m \geq 3$ such that $d_m({M})= (m-2, 1, 2, \ldots, 2, 1)$ is a dimension sequence.
 \end{The}
 \begin{proof}
 We start by some observation which is a direct consequence of Proposition \ref{4} (see also Proposition \ref{5}). For a basic hereditary ring $\Lambda$ ($\cong R_{\mathscr{M}_{\Lambda}}$) as above being right artinian, a finitely generated indecomposable right $\Lambda$-module $N=(X_F, Y_G, \varphi:X\otimes_F M \rightarrow Y)$ is multiplicity-free top if and only if either $N$ is a simple projective module and then $({\rm dim}{(X)}_F, {\rm dim}{(Y)}_G)=(0,1)$ or ${\rm dim}{(X)}_F =1$. This follows from the fact that $\varphi$ is surjective and $X_F \neq 0$ in case $N$ is not simple projective.\\
$(\Rightarrow)$. Assume that $\Lambda= \tiny{\begin{pmatrix}
F & M\\
0& G
\end{pmatrix}}$ is a right K\"othe ring. Then there exists $m\geq 3$ such that $\Lambda$ has only $m$ pairwise non-isomorphic finitely generated indecomposable right $\Lambda$-modules of finite length. Thus by \cite[Proposition 1.1]{Sim3}, $\mathscr{M}_{\Lambda}$ is a finite-dimensional species. Hence by Theorem \ref{DS1}, $\mathscr{M}_{\Lambda}$ has the finite-dimensional property and $\lbrace \textbf{P}_0, \textbf{P}_1, \ldots, \textbf{P}_{m-1} \rbrace$ is the set of all finite-dimensional indecomposable representations (up to isomorphism) of $\mathscr{M}_{\Lambda}$, where $\textbf{P}_0$ is a simple projective representation and $\textbf{P}_i=S_1^- \cdots S_i^- {\mathbf{F}}_{k^{'}_i}^{(i)}$, for each $1 \leq i \leq m-1$. By using of  \cite[Proposition 3.2]{Sim5} and \cite[Proposition 1.1]{Sim3} (see also \cite[Proposition 1]{DRS}), $M^{(m-2)L} \cong M^{2R}$ as bimodules and by using of \cite[Lemma 1.3]{Sim3} (see also \cite[Proposition 1]{DRS}), for each $1 \leq t \leq m-2$, ${\textbf{dim}}~\textbf{P}_{t+1} = d^M_{t}{\textbf{dim}}~\textbf{P}_t - {\textbf{dim}}~\textbf{P}_{t-1}$. It follows that $d_m(M)= \left( d^M_0, d^M_1, \ldots, d^M_{m-1}\right)$ is a dimension sequence. Since $\Lambda$ is right K\"othe, by Proposition \ref{r3}, every finite-dimensional indecomposable representation of $\mathscr{M}_{\Lambda}$ is multiplicity-free top. Assume that $m > 3$ and ${\textbf{dim}}~\textbf{P}_i = (x_i, y_i)$, for each $0 \leq i \leq m-1$ (comparing to the formula before Lemma \ref{2.12.1} for the sequence $a:=d_m(M)$, the indexing of $x_i$'s and $y_i$'s are shifted by 1). Then by the entrance observation $x_i = 1$, for every $i=1, \ldots, m-1$.  Since $\textbf{dim}~\textbf{P}_0 = (0, 1)$ and $\textbf{dim}~\textbf{P}_1 = (1, d^M_0)$, we have $d^M_1= x_{2} = 1$ and consequently $d^M_0 > 1$, since otherwise $m=3$ by Lemma \ref{2.12.1}. Also, we have $x_3 = d^M_2 -1$, so $d^M_2 = 2$. Proceeding in a similar way, by induction we obtain $$d_m(M)= \left( d^M_0, 1, 2, \ldots, 2, d^M_{m-1}\right).$$ Consequently, by \cite[Proposition $2^{'}$]{DRS} and Lemma \ref{2.12.1}, $\left( d^M_0 - (m-3), 1, d^M_{m-1}\right) $ is a dimension sequence. Hence, we infer that  $d^M_0= m-2$ and $d^M_{m-1} = 1$. Therefore $d_m(M)= (m-2, 1, 2, \ldots, 2, 1)$. \\
$(\Leftarrow)$. Assume that there exists $m \geq 3$ such that $d_m({M})= (m-2, 1, 2, \ldots, 2, 1)$ is a dimension sequence. Then by \cite[Corollary 3.5]{Sim5}, $\Lambda$ has only $m$ pairwise non-isomorphic finitely generated indecomposable right $\Lambda$-modules of finite length. Thus $\Lambda$ is representation-finite. It follows from \cite[Proposition 1.1]{Sim3} that $\Lambda$ is a basic hereditary artinian ring and $\mathscr{M}_{\Lambda}$ is representation-finite. Therefore by Theorem \ref{DS1}, $\lbrace \textbf{P}_0, \textbf{P}_1, \ldots, \textbf{P}_{m-1} \rbrace$ is the set of all finite-dimensional indecomposable representations (up to isomorphism) of $\mathscr{M}_{\Lambda}$, where $\textbf{P}_0$ is a simple projective representation and $\textbf{P}_i=S_1^- \cdots S_i^- {\mathbf{F}}_{k^{'}_i}^{(i)}$, for each $1 \leq i \leq m-1$. Since $d_m({M})= (m-2, 1, 2, \ldots, 2, 1)$ and by using \cite[Lemma 1.3]{Sim3} (see also \cite[Proposition 1]{DRS}), ${\textbf{dim}}~\textbf{P}_{t+1} = d^M_{t}{\textbf{dim}}~\textbf{P}_t - {\textbf{dim}}~\textbf{P}_{t-1}$ for each $1 \leq t \leq m-2$,  so computing inductively the consecutive dimension vectors, with the starting values $\textbf{dim}~\textbf{P}_0 = (0, 1)$ and $\textbf{dim}~\textbf{P}_1 = (1, m-2)$, we infer that $x_i=1$, for every $i=1, \ldots, m-1$, where $(x_i,y_i)$ are as above. Thus by entrance observation and Proposition \ref{4}, every finite-dimensional indecomposable representation of $\mathscr{M}_{\Lambda}$ is multiplicity-free top. Therefore by Proposition \ref{r3}, $\Lambda$ is a right K\"othe ring.
 \end{proof}

\begin{Examp}
{\rm Let $\Lambda=\tiny{\begin{pmatrix}
{\Bbb{C}} & {\Bbb{C}}\\
0& {\Bbb{R}}
\end{pmatrix}}$, where ${\Bbb{C}}$ is the field of complex numbers and ${\Bbb{R}}$ is the field of real numbers. Then $\Lambda$ is a basic hereditary artinian ring and  the category mod-$\Lambda$ has exactly 4 non-isomorphic indecomposable modules, since the underlying valued graph of $(\Gamma_{\Lambda}, \mathbf{d})$ is equal $\mathbb{B}_2$ from Table C and $\Lambda$ is a finite-dimensional algebra over its center ${\Bbb{R}}.{\rm I}_2$, in particular $(_{\Bbb{C}}{\Bbb{C}}_{\Bbb{R}})^L \cong (_{\Bbb{C}}{\Bbb{C}}_{\Bbb{R}})^R$, see \cite{DR1}. Thus the underlying Coxeter valued graph of the Coxeter valued quiver $(\mathscr{C}_{\Lambda}, \mathbf{m})$ of $\Lambda$ is $\mathbb{B}_2$ presented in Table B (i.e., \dynkin[Coxeter, edge length=1.0cm, gonality={4}]{B}{2}) and by the very definition  $d_4({_{{\Bbb{C}}}{\Bbb{C}}_{{\Bbb{R}}}})=(2,1,2,1)$.  Therefore by Theorem \ref{3.7}, $\Lambda$ is a right K\"othe ring. 

Let now $\Lambda=\tiny{\begin{pmatrix}
{\Bbb{R}} & {\Bbb{C}}\\
0& {\Bbb{C}}
\end{pmatrix}}$. Then $\Lambda$ has the analogous properties as in the previous case with the one exception; namely, $d_4({_{{\Bbb{R}}}{\Bbb{C}}_{{\Bbb{C}}}})=(1,2,1,2)$, hence $\Lambda$ is not a right K\"othe ring. Notice that in this case $\textbf{dim}~\textbf{P}_2 = (2, 1)$, so the representation $\textbf{P}_2$ of $\mathscr{M}_{\Lambda}$ is not multiplicity-free top.}
\end{Examp}

Let $\Lambda$ be a basic hereditary artinian ring such that the underlying Coxeter valued graph of the Coxeter valued quiver $(\mathscr{C}_{\Lambda}, \mathbf{m})$ of $\Lambda$ is one of the Coxeter diagrams $\mathbb{B}_n$, $\mathbb{F}_4$, $\mathbb{H}_3$, $\mathbb{H}_4$, $\mathbb{G}_2$ and $\mathbb{I}_2(p)$ presented in Table B. Then by Proposition \ref{r2}, $\Lambda \cong R_{{\mathscr{M}}_{\Lambda}}$. Assume that $\Lambda$ is a right K\"othe ring. Then by Proposition \ref{r3},  ${\mathscr{M}}_{\Lambda}$ has the multiplicity-free top property. Let ${\mathscr{M}}^{'}= (F, G, M)$ be a subspecies of ${\mathscr{M}}_{\Lambda}$ such that $\tiny{\begin{pmatrix}
F & M\\
0& G
\end{pmatrix}}$ has only $m \geq 3$ pairwise non-isomorphic indecomposable right modules. Therefore by Theorem \ref{3.7},  $d_m({M})= (m-2, 1, 2, \ldots, 2, 1)$ is a dimension sequence.

\begin{Pro}\label{3.1}
Let $\Lambda$ be a basic hereditary ring of the Dynkin type $\mathbb{B}_n$. Assume that ${\mathscr{M}}_{\Lambda} = \left( F_i,{_iM_j} \right)_{i, j \in I} $ is the species of $\Lambda$. Then $\Lambda$ is a right K\"othe ring if and only if $\Lambda$ is an artinian ring such that one of the following conditions holds:
\begin{itemize}
\item[$(a)$] $d_4({_{n-1}M_n}) = (2,1,2,1)$ is a dimension sequence and the Coxeter valued quiver $(\mathscr{C}_{\Lambda}, \mathbf{m})$ of $\Lambda$ is the quiver $$ 
	\definecolor{qqqqff}{rgb}{0.,0.,1.}
	\begin{tikzpicture}[line cap=round,line join=round,>=triangle 45,x=1.0cm,y=1.0cm]
	\draw (2.,9.)-- (3.,9.);
	\draw [dash pattern=on 2pt off 2pt](3.,9.)--(4.,9.);
	\draw (4.,9.)-- (5.,9.);
		\begin{scriptsize}
	\draw [fill=white] (2.,9.) circle (1.5pt);
	\draw[color=black] (2.,8.7) node {$n$};
	\draw[color=black] (2.5,9.2) node {$4$};
	\draw[color=black] (2.14,9.) node {$<$};
	\draw [fill=white] (3.,9.) circle (1.5pt);
	\draw[color=black] (3.,8.7) node {$n-1$};
	\draw [fill=white] (4.,9.) circle (1.5pt);
	\draw[color=black] (4.,8.7) node {$2$};
	\draw[color=black] (4.14,9.) node {$<$};
	\draw [fill=white] (5.,9.) circle (1.5pt);
	\draw[color=black] (5.,8.7) node {$1$};
	
	\end{scriptsize}
	\end{tikzpicture}
$$
\item[$(b)$] $d_4({_nM_{n-1}}) = (2,1,2,1)$ is a dimension sequence and the Coxeter valued quiver $(\mathscr{C}_{\Lambda}, \mathbf{m})$ of $\Lambda$ is the quiver $$ 
	\definecolor{qqqqff}{rgb}{0.,0.,1.}
	\begin{tikzpicture}[line cap=round,line join=round,>=triangle 45,x=1.0cm,y=1.0cm]
	\draw (2.,9.)-- (3.,9.);
	\draw [dash pattern=on 2pt off 2pt](3.,9.)--(4.,9.);
	\draw (4.,9.)-- (5.,9.);
	\draw (5.,9.)-- (6.,9.);
	\draw [dash pattern=on 2pt off 2pt](6.,9.)--(7.,9.);
	\draw (7.,9.)-- (8.,9.);
		\begin{scriptsize}
	\draw [fill=white] (2.,9.) circle (1.5pt);
	\draw[color=black] (2.,8.7) node {$n$};
	\draw[color=black] (2.5,9.2) node {$4$};
	\draw[color=black] (2.86,9.) node {$>$};
	\draw [fill=white] (3.,9.) circle (1.5pt);
	\draw[color=black] (3.,8.7) node {$n-1$};
	\draw [fill=white] (4.,9.) circle (1.5pt);
	\draw[color=black] (4.,8.7) node {$t-1$};
	\draw[color=black] (4.86,9.) node {$>$};
	\draw [fill=white] (5.,9.) circle (1.5pt);
	\draw[color=black] (5.,8.7) node {$t$};
	\draw [fill=white] (6.,9.) circle (1.5pt);
	\draw[color=black] (6.,8.7) node {$t+1$};
	\draw[color=black] (5.14,9.) node {$<$};
	\draw [fill=white] (7.,9.) circle (1.5pt);
	\draw[color=black] (7.,8.7) node {$2$};
	\draw[color=black] (7.14,9.) node {$<$};
	\draw [fill=white] (8.,9.) circle (1.5pt);
	\draw[color=black] (8.,8.7) node {$1$};
	
	\end{scriptsize}
	\end{tikzpicture}
$$
where $1 \leq t \leq n-1$.
\end{itemize}
 \end{Pro}
 \begin{proof}
 $(\Rightarrow)$. Assume that $\Lambda$ is a right K\"othe ring. Then by Proposition \ref{r3}, every finite-dimensional indecomposable representation of the species $\mathscr{M}_{\Lambda}$ has multiplicity-free top. The underlying Coxeter valued graph of  the Coxeter valued quiver $(\mathscr{C}_{\Lambda}, \mathbf{m})$ of $\Lambda$ is the Coxeter diagram
 \begin{center}
	\definecolor{qqqqff}{rgb}{0.,0.,1.}
	\begin{tikzpicture}[line cap=round,line join=round,>=triangle 45,x=1.0cm,y=1.0cm]
	\draw (2.,9.)-- (3.,9.);
	\draw [dash pattern=on 2pt off 2pt](3.,9.)--(4.,9.);
	\draw (4.,9.)-- (5.,9.);
		\begin{scriptsize}
	\draw[color=black] (.8,9.) node {$\mathbb{B}_n:$};
	\draw [fill=white] (2.,9.) circle (1.5pt);
	\draw[color=black] (2.,8.7) node {$n$};
	\draw[color=black] (2.5,9.2) node {$4$};
	\draw [fill=white] (3.,9.) circle (1.5pt);
	\draw[color=black] (3.,8.7) node {$n-1$};
	\draw [fill=white] (4.,9.) circle (1.5pt);
	\draw[color=black] (4.,8.7) node {$2$};
	\draw [fill=white] (5.,9.) circle (1.5pt);
	\draw[color=black] (5.,8.7) node {$1$};
	
	\end{scriptsize}
	\end{tikzpicture}
\end{center}
thus by Theorem \ref{3.7}, either $d_4({_{n-1}M_n}) = (2,1,2,1)$ is a dimension sequence or  $d_4({_{n}M_{n-1}}) = (2,1,2,1)$ is a dimension sequence. If
$d_4({_{n-1}M_n}) = (2,1,2,1)$ is a dimension sequence,  then by Proposition \ref{r2}, the underlying valued graph of the valued quiver $(\Gamma_{\Lambda}, \mathbf{d})$ is the valued Dynkin diagram $\mathbb{B}_n$ presented in Table C. Since there exists a bijection between the isomorphism classes of finite-dimensional indecomposable representations of $\mathscr{M}_{\Lambda}$ and the positive roots of $(\Gamma_{\Lambda}, \mathbf{d})$ and since by \cite[Pages 267-268]{Bourbaki}, the positive roots of $(\Gamma_{\Lambda}, \mathbf{d})$ can be expressed as combinations of simple roots as follows:
 \begin{itemize}
 \item[$(1)$] $\sum_{i \leq k \leq n}e_k$ \hspace{1cm} $(1 \leq i \leq n)$;
 \item[$(2)$] $\sum_{i \leq k \leq j}e_k$ \hspace{1cm} $(1 \leq i < j \leq n)$;
 \item[$(3)$]  $\sum_{i \leq k < j}e_k + 2 \sum_{j \leq k \leq n}e_k$ \hspace{1cm} $(1 \leq i < j \leq n)$,
 \end{itemize}
 thus by Proposition \ref{4}, the Coxeter valued quiver $(\mathscr{C}_{\Lambda}, \mathbf{m})$ of $\Lambda$ is
\begin{center}
	\definecolor{qqqqff}{rgb}{0.,0.,1.}
	\begin{tikzpicture}[line cap=round,line join=round,>=triangle 45,x=1.0cm,y=1.0cm]
	\draw (2.,9.)-- (3.,9.);
	\draw [dash pattern=on 2pt off 2pt](3.,9.)--(4.,9.);
	\draw (4.,9.)-- (5.,9.);
		\begin{scriptsize}
	\draw [fill=white] (2.,9.) circle (1.5pt);
	\draw[color=black] (2.,8.7) node {$n$};
	\draw[color=black] (2.5,9.2) node {$4$};
	\draw[color=black] (2.14,9.) node {$<$};
	\draw [fill=white] (3.,9.) circle (1.5pt);
	\draw[color=black] (3.,8.7) node {$n-1$};
	\draw [fill=white] (4.,9.) circle (1.5pt);
	\draw[color=black] (4.,8.7) node {$2$};
	\draw[color=black] (4.14,9.) node {$<$};
	\draw [fill=white] (5.,9.) circle (1.5pt);
	\draw[color=black] (5.,8.7) node {$1$};
	
	\end{scriptsize}
	\end{tikzpicture}
\end{center}
 If $d_4({_{n}M_{n-1}}) = (2,1,2,1)$ is a dimension sequence, then by Proposition \ref{r2}, the underlying valued graph of the valued quiver $(\Gamma_{\Lambda}, \mathbf{d})$ is the valued Dynkin diagram $\mathbb{C}_n$ presented in Table C. Since there exists a bijection between the isomorphism classes of finite-dimensional indecomposable representations of $\mathscr{M}_{\Lambda}$ and the positive roots of $(\Gamma_{\Lambda}, \mathbf{d})$ and also since by using \cite[Pages 269-270]{Bourbaki}, the positive roots of $(\Gamma_{\Lambda}, \mathbf{d})$ can be expressed as combinations of simple roots as follows:
 \begin{itemize}
 \item[$(1)$] $\sum_{i \leq k < j}e_k$ \hspace{1cm} $(1 \leq i < j \leq n)$;
 \item[$(2)$] $\sum_{i \leq k < j}e_k + 2 \sum_{j \leq k < n}e_k + e_n$ \hspace{1cm} $(1 \leq i < j \leq n)$;
 \item[$(3)$]  $\sum_{i \leq k < n}e_k +e_n$ \hspace{1cm} $(1 \leq i \leq n)$,
 \end{itemize}
thus by Proposition \ref{4}, the Coxeter valued quiver $(\mathscr{C}_{\Lambda}, \mathbf{m})$ of $\Lambda$ is
\begin{center}
	\definecolor{qqqqff}{rgb}{0.,0.,1.}
	\begin{tikzpicture}[line cap=round,line join=round,>=triangle 45,x=1.0cm,y=1.0cm]
	\draw (2.,9.)-- (3.,9.);
	\draw [dash pattern=on 2pt off 2pt](3.,9.)--(4.,9.);
	\draw (4.,9.)-- (5.,9.);
	\draw (5.,9.)-- (6.,9.);
	\draw [dash pattern=on 2pt off 2pt](6.,9.)--(7.,9.);
	\draw (7.,9.)-- (8.,9.);
		\begin{scriptsize}
	\draw [fill=white] (2.,9.) circle (1.5pt);
	\draw[color=black] (2.,8.7) node {$n$};
	\draw[color=black] (2.5,9.2) node {$4$};
	\draw[color=black] (2.86,9.) node {$>$};
	\draw [fill=white] (3.,9.) circle (1.5pt);
	\draw[color=black] (3.,8.7) node {$n-1$};
	\draw [fill=white] (4.,9.) circle (1.5pt);
	\draw[color=black] (4.,8.7) node {$t-1$};
	\draw[color=black] (4.86,9.) node {$>$};
	\draw [fill=white] (5.,9.) circle (1.5pt);
	\draw[color=black] (5.,8.7) node {$t$};
	\draw [fill=white] (6.,9.) circle (1.5pt);
	\draw[color=black] (6.,8.7) node {$t+1$};
	\draw[color=black] (5.14,9.) node {$<$};
	\draw [fill=white] (7.,9.) circle (1.5pt);
	\draw[color=black] (7.,8.7) node {$2$};
	\draw[color=black] (7.14,9.) node {$<$};
	\draw [fill=white] (8.,9.) circle (1.5pt);
	\draw[color=black] (8.,8.7) node {$1$};
	
	\end{scriptsize}
	\end{tikzpicture}
\end{center}
where $1 \leq t \leq n-1$.\\
$(\Leftarrow)$. Assume that $d_4({_{n-1}M_n}) = (2,1,2,1)$ is a dimension sequence and the Coxeter valued quiver $(\mathscr{C}_{\Lambda}, \mathbf{m})$ of $\Lambda$ is the quiver $$ 
	\definecolor{qqqqff}{rgb}{0.,0.,1.}
	\begin{tikzpicture}[line cap=round,line join=round,>=triangle 45,x=1.0cm,y=1.0cm]
	\draw (2.,9.)-- (3.,9.);
	\draw [dash pattern=on 2pt off 2pt](3.,9.)--(4.,9.);
	\draw (4.,9.)-- (5.,9.);
		\begin{scriptsize}
	\draw [fill=white] (2.,9.) circle (1.5pt);
	\draw[color=black] (2.,8.7) node {$n$};
	\draw[color=black] (2.5,9.2) node {$4$};
	\draw[color=black] (2.14,9.) node {$<$};
	\draw [fill=white] (3.,9.) circle (1.5pt);
	\draw[color=black] (3.,8.7) node {$n-1$};
	\draw [fill=white] (4.,9.) circle (1.5pt);
	\draw[color=black] (4.,8.7) node {$2$};
	\draw[color=black] (4.14,9.) node {$<$};
	\draw [fill=white] (5.,9.) circle (1.5pt);
	\draw[color=black] (5.,8.7) node {$1$};
	
	\end{scriptsize}
	\end{tikzpicture}
$$
Then by Proposition \ref{r2}, the valued quiver $(\Gamma_{\Lambda}, \mathbf{d})$ of $\Lambda$ is the valued Dynkin quiver
$$ 
	\definecolor{qqqqff}{rgb}{0.,0.,1.}
	\begin{tikzpicture}[line cap=round,line join=round,>=triangle 45,x=1.0cm,y=1.0cm]
	\draw (2.,9.)-- (3.,9.);
	\draw [dash pattern=on 2pt off 2pt](3.,9.)--(4.,9.);
	\draw (4.,9.)-- (5.,9.);
		\begin{scriptsize}
	\draw[color=black] (.8,9.) node {$\mathbb{B}_n:$};
	\draw [fill=white] (2.,9.) circle (1.5pt);
	\draw[color=black] (2.,8.7) node {$n$};
	\draw[color=black] (2.5,9.2) node {$(1,2)$};
	\draw[color=black] (2.14,9.) node {$<$};
	\draw [fill=white] (3.,9.) circle (1.5pt);
	\draw[color=black] (3.,8.7) node {$n-1$};
	\draw [fill=white] (4.,9.) circle (1.5pt);
	\draw[color=black] (4.,8.7) node {$2$};
	\draw[color=black] (4.14,9.) node {$<$};
	\draw [fill=white] (5.,9.) circle (1.5pt);
	\draw[color=black] (5.,8.7) node {$1$};
	
	\end{scriptsize}
	\end{tikzpicture}
$$
Thus by using \cite[Pages 267-268]{Bourbaki} and Proposition \ref{4}, $\mathscr{M}_{\Lambda}$ has the multilplicity-free top property. Therefore by Propositions \ref{r2} and \ref{r3}, $\Lambda$ is a right K\"othe ring. \\
Now we assume that $d_4({_nM_{n-1}}) = (2,1,2,1)$ is a dimension sequence and the Coxeter valued quiver $(\mathscr{C}_{\Lambda}, \mathbf{m})$ of $\Lambda$ is the quiver $$ 
	\definecolor{qqqqff}{rgb}{0.,0.,1.}
	\begin{tikzpicture}[line cap=round,line join=round,>=triangle 45,x=1.0cm,y=1.0cm]
	\draw (2.,9.)-- (3.,9.);
	\draw [dash pattern=on 2pt off 2pt](3.,9.)--(4.,9.);
	\draw (4.,9.)-- (5.,9.);
	\draw (5.,9.)-- (6.,9.);
	\draw [dash pattern=on 2pt off 2pt](6.,9.)--(7.,9.);
	\draw (7.,9.)-- (8.,9.);
		\begin{scriptsize}
	\draw [fill=white] (2.,9.) circle (1.5pt);
	\draw[color=black] (2.,8.7) node {$n$};
	\draw[color=black] (2.5,9.2) node {$4$};
	\draw[color=black] (2.86,9.) node {$>$};
	\draw [fill=white] (3.,9.) circle (1.5pt);
	\draw[color=black] (3.,8.7) node {$n-1$};
	\draw [fill=white] (4.,9.) circle (1.5pt);
	\draw[color=black] (4.,8.7) node {$t-1$};
	\draw[color=black] (4.86,9.) node {$>$};
	\draw [fill=white] (5.,9.) circle (1.5pt);
	\draw[color=black] (5.,8.7) node {$t$};
	\draw [fill=white] (6.,9.) circle (1.5pt);
	\draw[color=black] (6.,8.7) node {$t+1$};
	\draw[color=black] (5.14,9.) node {$<$};
	\draw [fill=white] (7.,9.) circle (1.5pt);
	\draw[color=black] (7.,8.7) node {$2$};
	\draw[color=black] (7.14,9.) node {$<$};
	\draw [fill=white] (8.,9.) circle (1.5pt);
	\draw[color=black] (8.,8.7) node {$1$};
	
	\end{scriptsize}
	\end{tikzpicture}
$$
where $1 \leq t \leq n-1$. Thus by Proposition \ref{r2}, the valued quiver $(\Gamma_{\Lambda}, \mathbf{d})$ of $\Lambda$ is the valued Dynkin quiver
$$ 
	\definecolor{qqqqff}{rgb}{0.,0.,1.}
	\begin{tikzpicture}[line cap=round,line join=round,>=triangle 45,x=1.0cm,y=1.0cm]
	\draw (2.,9.)-- (3.,9.);
	\draw [dash pattern=on 2pt off 2pt](3.,9.)--(4.,9.);
	\draw (4.,9.)-- (5.,9.);
	\draw (5.,9.)-- (6.,9.);
	\draw [dash pattern=on 2pt off 2pt](6.,9.)--(7.,9.);
	\draw (7.,9.)-- (8.,9.);
		\begin{scriptsize}
	\draw[color=black] (.8,9.) node {$\mathbb{C}_n:$};
	\draw [fill=white] (2.,9.) circle (1.5pt);
	\draw[color=black] (2.,8.7) node {$n$};
	\draw[color=black] (2.5,9.2) node {$(2,1)$};
	\draw[color=black] (2.86,9.) node {$>$};
	\draw [fill=white] (3.,9.) circle (1.5pt);
	\draw[color=black] (3.,8.7) node {$n-1$};
	\draw [fill=white] (4.,9.) circle (1.5pt);
	\draw[color=black] (4.,8.7) node {$t-1$};
	\draw[color=black] (4.86,9.) node {$>$};
	\draw [fill=white] (5.,9.) circle (1.5pt);
	\draw[color=black] (5.,8.7) node {$t$};
	\draw [fill=white] (6.,9.) circle (1.5pt);
	\draw[color=black] (6.,8.7) node {$t+1$};
	\draw[color=black] (5.14,9.) node {$<$};
	\draw [fill=white] (7.,9.) circle (1.5pt);
	\draw[color=black] (7.,8.7) node {$2$};
	\draw[color=black] (7.14,9.) node {$<$};
	\draw [fill=white] (8.,9.) circle (1.5pt);
	\draw[color=black] (8.,8.7) node {$1$};
	
	\end{scriptsize}
	\end{tikzpicture}
$$
where $1 \leq t \leq n-1$. It follows that by using \cite[Pages 269-270]{Bourbaki} and Proposition \ref{4}, $\mathscr{M}_{\Lambda}$ has the multilplicity-free top property. Consequently by Propositions \ref{r2} and \ref{r3}, $\Lambda$ is a right K\"othe ring.
 \end{proof}

\begin{Lem}\label{m7}
If $\Lambda$ is a basic hereditary artinian ring of the Dynkin type $\mathbb{F}_4$, then $\Lambda$ is not a right K\"othe ring.
\end{Lem}
\begin{proof}
Let $F$ and $G$ be division rings and $M$ be an $F$-$G$-bimodule such that ${\mathscr{M}}^{'}= (F, G, M)$ is a subspecies of ${\mathscr{M}}_{\Lambda}$ and $\tiny{\begin{pmatrix}
F & M\\
0& G
\end{pmatrix}}$ has only $4$ pairwise non-isomorphic indecomposable right modules. Assume to the contrary that $\Lambda$ is a right K\"othe ring. Then by Theorem \ref{3.7}, $d_4(M)=(2,1,2,1)$. It follow from Proposition \ref{r2} that the valued quiver $(\Gamma, \mathbf{d})$ of $\Lambda$ is one of the valued Dynkin quivers
\begin{center}
	\definecolor{qqqqff}{rgb}{0.,0.,1.}
	\begin{tikzpicture}[line cap=round,line join=round,>=triangle 45,x=1.0cm,y=1.0cm]
	\draw (2.,9.)-- (3.,9.);
	\draw (3.,9.)-- (4.,9.);
	\draw (4.,9.)-- (5.,9.);
		\begin{scriptsize}
	\draw[color=black] (.8,9.) node {$\mathbb{F}_4:$};
	\draw [fill=white] (2.,9.) circle (1.5pt);
	\draw[color=black] (2.87,9.) node {$>$};
	\draw[color=black] (3.14,9.) node {$<$};
	\draw [fill=white] (3.,9.) circle (1.5pt);
	\draw[color=black] (3.5,9.2) node {$(1,2)$};
	\draw [fill=white] (4.,9.) circle (1.5pt);
	\draw[color=black] (4.14,9.) node {$<$};
	\draw [fill=white] (5.,9.) circle (1.5pt);
	\draw[color=black] (5.1,9.) node {;};
	
	\end{scriptsize}
	\end{tikzpicture}
\end{center}
\begin{center}
	\definecolor{qqqqff}{rgb}{0.,0.,1.}
	\begin{tikzpicture}[line cap=round,line join=round,>=triangle 45,x=1.0cm,y=1.0cm]
	\draw (2.,9.)-- (3.,9.);
	\draw (3.,9.)-- (4.,9.);
	\draw (4.,9.)-- (5.,9.);
		\begin{scriptsize}
	\draw[color=black] (.8,9.) node {$\mathbb{F}_4:$};
	\draw [fill=white] (2.,9.) circle (1.5pt);
	\draw[color=black] (2.87,9.) node {$>$};
	\draw[color=black] (3.14,9.) node {$<$};
	\draw [fill=white] (3.,9.) circle (1.5pt);
	\draw[color=black] (3.5,9.2) node {$(1,2)$};
	\draw [fill=white] (4.,9.) circle (1.5pt);
	\draw [fill=white] (5.,9.) circle (1.5pt);
	\draw[color=black] (5.1,9.) node {;};
	\draw[color=black] (4.87,9.) node {$>$};
	
	\end{scriptsize}
	\end{tikzpicture}
\end{center}
\begin{center}
	\definecolor{qqqqff}{rgb}{0.,0.,1.}
	\begin{tikzpicture}[line cap=round,line join=round,>=triangle 45,x=1.0cm,y=1.0cm]
	\draw (2.,9.)-- (3.,9.);
	\draw (3.,9.)-- (4.,9.);
	\draw (4.,9.)-- (5.,9.);
		\begin{scriptsize}
	\draw[color=black] (.8,9.) node {$\mathbb{F}_4:$};
	\draw [fill=white] (2.,9.) circle (1.5pt);
	\draw[color=black] (2.14,9.) node {$<$};
	\draw[color=black] (3.14,9.) node {$<$};
	\draw [fill=white] (3.,9.) circle (1.5pt);
	\draw[color=black] (3.5,9.2) node {$(1,2)$};
	\draw [fill=white] (4.,9.) circle (1.5pt);
	\draw [fill=white] (5.,9.) circle (1.5pt);
	\draw[color=black] (4.87,9.) node {$>$};
	\draw[color=black] (5.1,9.) node {;};
	
	\end{scriptsize}
	\end{tikzpicture}
\end{center}
\begin{center}
	\definecolor{qqqqff}{rgb}{0.,0.,1.}
	\begin{tikzpicture}[line cap=round,line join=round,>=triangle 45,x=1.0cm,y=1.0cm]
	\draw (2.,9.)-- (3.,9.);
	\draw (3.,9.)-- (4.,9.);
	\draw (4.,9.)-- (5.,9.);
		\begin{scriptsize}
	\draw[color=black] (.8,9.) node {$\mathbb{F}_4:$};
	\draw [fill=white] (2.,9.) circle (1.5pt);
	\draw[color=black] (2.14,9.) node {$<$};
	\draw[color=black] (3.14,9.) node {$<$};
	\draw [fill=white] (3.,9.) circle (1.5pt);
	\draw[color=black] (3.5,9.2) node {$(1,2)$};
	\draw [fill=white] (4.,9.) circle (1.5pt);
	\draw [fill=white] (5.,9.) circle (1.5pt);
	\draw[color=black] (4.14,9.) node {$<$};
	
	\end{scriptsize}
	\end{tikzpicture}
\end{center}
Since by using \cite[Page 287]{Bourbaki}, there exists a finite-dimensional indecomposable representation $X$ of ${\mathscr{M}}_{\Lambda}$ with the dimension vector $\textbf{dim}~X= {\bsm 2342\esm}$, by Proposition \ref{4}, ${\mathscr{M}}_{\Lambda}$ has no the  multiplicity-free top property, which is a contradiction, by Proposition \ref{r3}. Therefore $\Lambda$ is not right K\"othe.
\end{proof}

\begin{Lem}\label{m33}
If $\Lambda$ is a basic hereditary artinian ring of the Coxeter type $\mathbb{H}_3$, then $\Lambda$ is not a right K\"othe ring.
\end{Lem}
 \begin{proof}
Let $F$ and $G$ be division rings and $M$ be an $F$-$G$-bimodule such that ${\mathscr{M}}^{'}= (F, G, M)$ is a subspecies of ${\mathscr{M}}_{\Lambda}$ and $\tiny{\begin{pmatrix}
F & M\\
0& G
\end{pmatrix}}$ has only $5$ pairwise non-isomorphic indecomposable right modules. Assume to the contrary that $\Lambda$ is a right K\"othe ring. Then by Proposition \ref{r3}, every finite-dimensional indecomposable representation of the species ${\mathscr{M}}_{\Lambda}$ has multiplicity-free top. Thus by Theorem \ref{3.7}, $d_5(M) = (3,1,2,2,1)$. If $(\mathscr{C}_{\Lambda}, \mathbf{m})$ is the quiver
\begin{center}
	\definecolor{qqqqff}{rgb}{0.,0.,1.}
	\begin{tikzpicture}[line cap=round,line join=round,>=triangle 45,x=1.0cm,y=1.0cm]
	\draw (2.,9.)-- (3.,9.);
	\draw (3.,9.)-- (4.,9.);
		\begin{scriptsize}
	\draw[color=black] (.8,9.) node {};
	\draw [fill=white] (2.,9.) circle (1.5pt);
	\draw[color=black] (2.,8.7) node {$1$};
	\draw[color=black] (2.85,9.) node {$>$};
	\draw [fill=white] (3.,9.) circle (1.5pt);
	\draw[color=black] (3.,8.7) node {$2$};
	\draw[color=black] (3.5,9.2) node {$5$};
	\draw[color=black] (3.85,9.) node {$>$};
	\draw [fill=white] (4.,9.) circle (1.5pt);
	\draw[color=black] (4.,8.7) node {$3$};
	
	\end{scriptsize}
	\end{tikzpicture}
\end{center}
then $\textbf{dim}~\textbf{P}_5 = s_1^-s_2^-s_3^-s_4^-s_5^-(0,0,1) = s_1^-s_2^-s_3^-s_4^-(0,2,1)= s_1^-s_2^-s_3^-(0,2,3)= s_1^-s_2^-(2,2,3)= s_1^-(2,3,3)= (2,3,6)$. Hence by Proposition \ref{4}, ${\mathbf{F}}_1 \oplus {\mathbf{F}}_1 \subseteq {\rm top}(\textbf{P}_5)$, which is a contradiction. If $(\mathscr{C}_{\Lambda}, \mathbf{m})$ is the quiver
 \begin{center}
	\definecolor{qqqqff}{rgb}{0.,0.,1.}
	\begin{tikzpicture}[line cap=round,line join=round,>=triangle 45,x=1.0cm,y=1.0cm]
	\draw (2.,9.)-- (3.,9.);
	\draw (3.,9.)-- (4.,9.);
		\begin{scriptsize}
	\draw[color=black] (.8,9.) node {};
	\draw [fill=white] (2.,9.) circle (1.5pt);
	\draw[color=black] (2.,8.7) node {$1$};
	\draw[color=black] (2.85,9.) node {$>$};
	\draw [fill=white] (3.,9.) circle (1.5pt);
	\draw[color=black] (3.,8.7) node {$2$};
	\draw[color=black] (3.5,9.2) node {$5$};
	\draw[color=black] (3.14,9.) node {$<$};
	\draw [fill=white] (4.,9.) circle (1.5pt);
	\draw[color=black] (4.,8.7) node {$3$};
	
	\end{scriptsize}
	\end{tikzpicture}
\end{center}
then $\textbf{dim}~\textbf{P}_4 = s_1^-s_2^-s_3^-s_4^-(0,0,1) = s_1^-s_2^-s_3^-(0,2,1)= s_1^-s_2^-(0,2,1)= s_1^-(2,2,1)= (2,3,1)$. Thus ${\mathbf{F}}_1 \oplus {\mathbf{F}}_1 \subseteq {\rm top}(\textbf{P}_4)$, which is a contradiction. Suppose that $(\mathscr{C}_{\Lambda}, \mathbf{m})$ is the quiver
\begin{center}
	\definecolor{qqqqff}{rgb}{0.,0.,1.}
	\begin{tikzpicture}[line cap=round,line join=round,>=triangle 45,x=1.0cm,y=1.0cm]
	\draw (2.,9.)-- (3.,9.);
	\draw (3.,9.)-- (4.,9.);
		\begin{scriptsize}
	\draw[color=black] (.8,9.) node {};
	\draw [fill=white] (2.,9.) circle (1.5pt);
	\draw[color=black] (2.,8.7) node {$1$};
	\draw[color=black] (2.14,9.) node {$<$};
	\draw [fill=white] (3.,9.) circle (1.5pt);
	\draw[color=black] (3.,8.7) node {$2$};
	\draw[color=black] (3.5,9.2) node {$5$};
	\draw[color=black] (3.85,9.) node {$>$};
	\draw [fill=white] (4.,9.) circle (1.5pt);
	\draw[color=black] (4.,8.7) node {$3$};
	
	\end{scriptsize}
	\end{tikzpicture}
\end{center}
Then $\textbf{dim}~\textbf{P}_5 = s_1^-s_2^-s_3^-s_4^-s_5^-(0,1,0) = s_1^-s_2^-s_3^-s_4^-(0,1,2)= s_1^-s_2^-s_3^-(1,1,2)= s_1^-s_2^-(1,2,2)= s_1^-(1,2,4)= (1,2,4)$. Therefore ${\mathbf{F}}_2 \oplus {\mathbf{F}}_2 \subseteq {\rm top}(\textbf{P}_5)$, which is a contradiction. If $(\mathscr{C}_{\Lambda}, \mathbf{m})$ is the quiver
\begin{center}
	\definecolor{qqqqff}{rgb}{0.,0.,1.}
	\begin{tikzpicture}[line cap=round,line join=round,>=triangle 45,x=1.0cm,y=1.0cm]
	\draw (2.,9.)-- (3.,9.);
	\draw (3.,9.)-- (4.,9.);
		\begin{scriptsize}
	\draw[color=black] (.8,9.) node {};
	\draw [fill=white] (2.,9.) circle (1.5pt);
	\draw[color=black] (2.,8.7) node {$1$};
	\draw[color=black] (2.14,9.) node {$<$};
	\draw [fill=white] (3.,9.) circle (1.5pt);
	\draw[color=black] (3.,8.7) node {$2$};
	\draw[color=black] (3.5,9.2) node {$5$};
	\draw[color=black] (3.14,9.) node {$<$};
	\draw [fill=white] (4.,9.) circle (1.5pt);
	\draw[color=black] (4.,8.7) node {$3$};
	
	\end{scriptsize}
	\end{tikzpicture}
\end{center}
then $\textbf{dim}~\textbf{P}_7 = s_1^-s_2^-s_3^-s_4^-s_5^-s_6^-s_7^-(0,1,0) = s_1^-s_2^-s_3^-s_4^-s_5^-s_6^-(1,1,0)= s_1^-s_2^-s_3^-s_4^-s_5^-(1,1,2) =$\\$ s_1^-s_2^-s_3^-s_4^-(1,4,2)= s_1^-s_2^-s_3^-(3,4,2)= s_1^-s_2^-(3,4,2) = s_1^-(3,5,2)= (2,5,2)$. Thus by Proposition \ref{4}, ${\mathbf{F}}_3 \oplus {\mathbf{F}}_3 \subseteq {\rm top}(\textbf{P}_7)$, which is a contradiction. Therefore $\Lambda$ is not a right K\"othe ring.
 \end{proof}

As an immediate consequence of Lemma \ref{m33}, we have the following corollary.

\begin{Cor}\label{m34}
If $\Lambda$ is a basic hereditary artinian ring of the Coxeter type $\mathbb{H}_4$, then $\Lambda$ is not a right K\"othe ring.
\end{Cor}

The pair $(\Gamma, \mathbf{m})$ is called a {\it {\rm (}general{\rm )} Coxeter valued quiver} if $\Gamma=(\Gamma_0, \Gamma_1)$ is a finite quiver and $\mathbf{m}:\Gamma_1 \rightarrow \Bbb{N} \cup \lbrace \infty\rbrace$ is a function such that $\mathbf{m}(\alpha) \geq 3$, for any arrow $\alpha \in \Gamma_1$. Notice that each Coxeter valued quiver is uniquely determined by the underlying Coxeter valued graph and some selection of its orientation. 

Let $\Lambda$ be a basic hereditary ring and ${\mathscr{M}}_{\Lambda}=(F_i, {_iM_j})_{i,j \in I}$ be the species of $\Lambda$. Note that the species ${\mathscr{M}}_{\Lambda}$ doesn't necessarily have the property ``${{_iM_j}\neq 0}$ implies that ${_jM_i} = 0$".  We say that the bimodule $_iM_j$ belongs to a fixed connected component $(\Gamma, \mathbf{m})$ of the Coxeter valued quiver of $\Lambda$, provided we have $i, j \in \Gamma_0$. Recall that a ring $\Lambda$ is called {\it indecomposable} if $\Lambda$ is not a direct product of two non-zero rings. Now, we are ready to give a characterization of basic hereditary right K\"othe rings in terms of their Coxeter valued quivers.

\begin{The}\label{m1}
Let $\Lambda$ be a basic hereditary ring. Then $\Lambda$ is right K\"othe if and only if $\Lambda$ is an artinian ring such that the Coxeter valued quiver $(\mathscr{C}_{\Lambda}, \mathbf{m})$ of $\Lambda$ is a finite disjoint union of the following Coxeter valued quivers:
\begin{itemize}
\item[$(1)$] $\mathbb{A}_n$ with any orientation;
\item[$(2)$]$\mathbb{B}_n$ with the orientation $$ 
	\definecolor{qqqqff}{rgb}{0.,0.,1.}
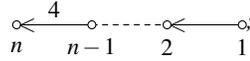
	\begin{tikzpicture}[line cap=round,line join=round,>=triangle 45,x=1.0cm,y=1.0cm]
	\draw (2.,9.)-- (3.,9.);
	\draw [dash pattern=on 2pt off 2pt](3.,9.)--(4.,9.);
	\draw (4.,9.)-- (5.,9.);
		\begin{scriptsize}
	\draw[color=black] (.8,9.) node {};
	\draw [fill=white] (2.,9.) circle (1.5pt);
	\draw[color=black] (2.,8.7) node {$n$};
	\draw[color=black] (2.5,9.2) node {$4$};
	\draw[color=black] (2.14,9.) node {$<$};
	\draw [fill=white] (3.,9.) circle (1.5pt);
	\draw[color=black] (3.,8.7) node {$n-1$};
	\draw [fill=white] (4.,9.) circle (1.5pt);
	\draw[color=black] (4.,8.7) node {$2$};
	\draw[color=black] (4.14,9.) node {$<$};
	\draw [fill=white] (5.,9.) circle (1.5pt);
	\draw[color=black] (5.,8.7) node {$1$};
	\draw [fill=black] (5.1,9.) node {;};
	\end{scriptsize}
	\end{tikzpicture}
$$

\item[$(3)$] $\mathbb{B}_n$ with the orientation $$ 
	\definecolor{qqqqff}{rgb}{0.,0.,1.}
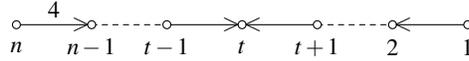
	\begin{tikzpicture}[line cap=round,line join=round,>=triangle 45,x=1.0cm,y=1.0cm]
	\draw (2.,9.)-- (3.,9.);
	\draw [dash pattern=on 2pt off 2pt](3.,9.)--(4.,9.);
	\draw (4.,9.)-- (5.,9.);
	\draw (5.,9.)-- (6.,9.);
	\draw [dash pattern=on 2pt off 2pt](6.,9.)--(7.,9.);
	\draw (7.,9.)-- (8.,9.);
		\begin{scriptsize}
	\draw [fill=white] (2.,9.) circle (1.5pt);
	\draw[color=black] (2.,8.7) node {$n$};
	\draw[color=black] (2.5,9.2) node {$4$};
	\draw[color=black] (2.86,9.) node {$>$};
	\draw [fill=white] (3.,9.) circle (1.5pt);
	\draw[color=black] (3.,8.7) node {$n-1$};
	\draw [fill=white] (4.,9.) circle (1.5pt);
	\draw[color=black] (4.,8.7) node {$t-1$};
	\draw[color=black] (4.86,9.) node {$>$};
	\draw [fill=white] (5.,9.) circle (1.5pt);
	\draw[color=black] (5.,8.7) node {$t$};
	\draw [fill=white] (6.,9.) circle (1.5pt);
	\draw[color=black] (6.,8.7) node {$t+1$};
	\draw[color=black] (5.14,9.) node {$<$};
	\draw [fill=white] (7.,9.) circle (1.5pt);
	\draw[color=black] (7.,8.7) node {$2$};
	\draw[color=black] (7.14,9.) node {$<$};
	\draw [fill=white] (8.,9.) circle (1.5pt);
	\draw[color=black] (8.,8.7) node {$1$};
	\end{scriptsize}
	\end{tikzpicture}
$$
with $1 \leq t \leq n-1$;
\item[$(4)$] $\mathbb{D}_n$ with the following conditions:
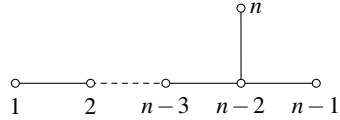
\begin{center}
 \definecolor{qqqqff}{rgb}{0.,0.,1.}
	\begin{tikzpicture}[line cap=round,line join=round,>=triangle 45,x=1.0cm,y=1.0cm]
	\draw (2.,9.)-- (3.,9.);
	\draw [dash pattern=on 2pt off 2pt](3.,9.)--(4.,9.);
	\draw (4.,9.)-- (5.,9.);
	\draw (5.,9.)-- (6.,9.);
	\draw (5.,9.)-- (6.,10.9999999999987,2.6);
		\begin{scriptsize}
	\draw [fill=white] (2.,9.) circle (1.5pt);
	\draw[color=black] (2.,8.7) node {$1$};
	\draw [fill=white] (3.,9.) circle (1.5pt);
	\draw[color=black] (3.,8.7) node {$2$};
	\draw [fill=white] (4.,9.) circle (1.5pt);
	\draw[color=black] (4.,8.7) node {$n-3$};
	\draw [fill=white] (5.,9.) circle (1.5pt);
	\draw[color=black] (5.,8.7) node {$n-2$};
	\draw [fill=white] (6.,9.) circle (1.5pt);
	\draw[color=black] (6.,8.7) node {$n-1$};
	\draw [fill=white] (6.,10.9999999999987,2.6) circle (1.5pt);
	\draw[color=black] (6.2,10.9999999999987,2.6) node {$n$};
	
	\end{scriptsize}
	\end{tikzpicture}
 \end{center}
 \begin{itemize}
\item[$(a)$] $|{(n-2)}^{+}| \leq 2$;
\item[$(b)$] For each $i\leq n-3$, there exists at most one arrow with the source $i$;
\end{itemize}
\item[$(5)$] $\mathbb{E}_6$ with the following conditions:
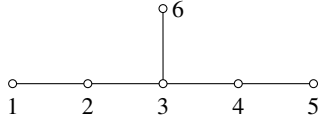
\begin{center}
 \definecolor{qqqqff}{rgb}{0.,0.,1.}
	\begin{tikzpicture}[line cap=round,line join=round,>=triangle 45,x=1.0cm,y=1.0cm]
	\draw (2.,9.)-- (3.,9.);
	\draw (3.,9.)-- (4.,9.);
	\draw (4.,9.)-- (5.,9.);
	\draw (4.,9.)-- (5.,10.9999999999987,2.6);
	\draw (5.,9.)-- (6.,9.);
		\begin{scriptsize}
	\draw [fill=white] (2.,9.) circle (1.5pt);
	\draw[color=black] (2.,8.7) node {$1$};
	\draw [fill=white] (3.,9.) circle (1.5pt);
	\draw[color=black] (3.,8.7) node {$2$};
	\draw [fill=white] (4.,9.) circle (1.5pt);
	\draw[color=black] (4.,8.7) node {$3$};
	\draw [fill=white] (5.,9.) circle (1.5pt);
	\draw[color=black] (5.,8.7) node {$4$};
	\draw [fill=white] (5.,10.9999999999987,2.6) circle (1.5pt);
	\draw[color=black] (5.2,10.9999999999987,2.6) node {$6$};
	\draw [fill=white] (6.,9.) circle (1.5pt);
	\draw[color=black] (6.,8.7) node {$5$};
	\end{scriptsize}
	\end{tikzpicture}
 \end{center}
\begin{itemize}
\item[$(a)$] $1 \leq |3^{+}| \leq 2$, $|4^{+}| \leq 1$ and $|2^{+}| \leq 1$;
\item[$(b)$] For each $y \in 3^{-}$, there exists at least one arrow with the target $y$;
\end{itemize}
\item[$(6)$] $\mathbb{E}_7$ with the orientation
$$
	\definecolor{qqqqff}{rgb}{0.,0.,1.}
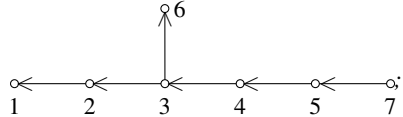
	\begin{tikzpicture}[line cap=round,line join=round,>=triangle 45,x=1.0cm,y=1.0cm]
	\draw (2.,9.)-- (3.,9.);
	\draw (3.,9.)-- (4.,9.);
	\draw (4.,9.)-- (5.,9.);
	\draw (4.,9.)-- (5.,10.9999999999987,2.6);
	\draw (5.,9.)-- (6.,9.);
	\draw (6.,9.)-- (7.,9.);
		\begin{scriptsize}
	\draw[color=black] (.8,9.) node {};
	\draw [fill=white] (2.,9.) circle (1.5pt);
	\draw[color=black] (2.,8.7) node {$1$};
	\draw[color=black] (2.14,9.) node {$<$};
	\draw [fill=white] (3.,9.) circle (1.5pt);
	\draw[color=black] (3.,8.7) node {$2$};
	\draw[color=black] (3.14,9.) node {$<$};
	\draw [fill=white] (4.,9.) circle (1.5pt);
	\draw[color=black] (4.,8.7) node {$3$};
	\draw[color=black] (4.14,9.) node {$<$};
	\draw [fill=white] (5.,9.) circle (1.5pt);
	\draw[color=black] (5.,8.7) node {$4$};
    \draw[color=black] (5.14,9.) node {$<$};
	\draw [fill=white] (5.,10.9999999999987,2.6) circle (1.5pt);
	\draw[color=black] (5.2,10.9999999999987,2.6) node {$6$};
	\draw[color=black] (4.,9.85) node {$\wedge$};
	\draw [fill=white] (6.,9.) circle (1.5pt);
	\draw[color=black] (6.,8.7) node {$5$};
	\draw [fill=white] (7.,9.) circle (1.5pt);
	\draw[color=black] (7.,8.7) node {$7$};
    \draw[color=black] (6.17,9.) node {$<$};
	\draw [fill=black] (7.1,9.) node {;};
	\end{scriptsize}
	\end{tikzpicture}$$
	\item[$(7)$] $\mathbb{G}_2$ with the orientation $$
	\definecolor{qqqqff}{rgb}{0.,0.,1.}
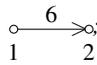
	\begin{tikzpicture}[line cap=round,line join=round,>=triangle 45,x=1.0cm,y=1.0cm]
	\draw (2.,9.)-- (3.,9.);
		\begin{scriptsize}
	\draw [fill=white] (2.,9.) circle (1.5pt);
	\draw[color=black] (2.,8.7) node {$1$};
	\draw[color=black] (2.5,9.2) node {$6$};
	\draw[color=black] (2.87,9.) node {$>$};
	\draw [fill=white] (3.,9.) circle (1.5pt);
	\draw[color=black] (3.,8.7) node {$2$};
	\draw [fill=black] (3.1,9.) node {;};
	\end{scriptsize}
	\end{tikzpicture}
$$

\item[$(8)$] $\mathbb{I}_2(p)$ with the orientation $$
	\definecolor{qqqqff}{rgb}{0.,0.,1.}
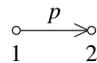
	\begin{tikzpicture}[line cap=round,line join=round,>=triangle 45,x=1.0cm,y=1.0cm]
	\draw (2.,9.)-- (3.,9.);
		\begin{scriptsize}
	\draw [fill=white] (2.,9.) circle (1.5pt);
	\draw[color=black] (2.,8.7) node {$1$};
	\draw[color=black] (2.5,9.2) node {$p$};
	\draw[color=black] (2.87,9.) node {$>$};
	\draw [fill=white] (3.,9.) circle (1.5pt);
	\draw[color=black] (3.,8.7) node {$2$};
	
	\end{scriptsize}
	\end{tikzpicture}
$$ with $p=5$ or $7 \leq p < \infty$
\end{itemize}
where additionally in the cases (2), (3), (7) and (8) the dimension sequences of the unique nontrivial bimodules ${_iM_j}$ belong to these components have very restrictive shapes given respectively as follows:
\begin{itemize}
\item[$(2)$] $d_4({_{n-1}M_n}) = (2,1,2,1)$;
\item[$(3)$] $d_4({_nM_{n-1}}) = (2,1,2,1)$;
\item[$(7)$] $d_6({_1M_2})=(4,1,2,2,2,1)$;
\item[$(8)$] $d_p({_1M_2})=(p-2, 1, 2, \ldots, 2, 1)$.
\end{itemize}
\end{The}
\begin{proof}
 Since $\Lambda$ is a basic hereditary ring,  it has a ring product decomposition $\Lambda= \Lambda_1 \oplus \cdots \oplus \Lambda_s$, where each $\Lambda_i$ is an indecomposable basic hereditary ring. Then the Coxeter valued quiver $(\mathscr{C}_{\Lambda}, \mathbf{m})$ of $\Lambda$ is a disjoint union of the Coxeter valued quivers $(\mathscr{C}_{\Lambda_i}, \mathbf{m})$ of $\Lambda_i$. \\
$(\Rightarrow)$. Assume that $\Lambda$ is a right K\"othe ring. Then $\Lambda_i$ is an indecomposable basic hereditary right K\"othe ring, for each $1 \leq i \leq s$.  Thus without loss of generality, we can assume that $\Lambda$ is an indecomposable basic hereditary right K\"othe ring.  Then by Proposition \ref{r3}, $\Lambda$ is a representation-finite ring and every finite-dimensional indecomposable representation of ${\mathscr{M}}_{\Lambda}$ is multiplicity-free top. Thus by Proposition \ref{r2} and Lemmas \ref{m6}, \ref{m7}, \ref{m33} and Corollary \ref{m34}, the underlying Coxeter valued graph of the Coxeter valued quiver $(\mathscr{C}_{\Lambda}, \mathbf{m})$ of $\Lambda$ is one of the Coxeter diagrams $\mathbb{A}_n$, $\mathbb{B}_n$, $\mathbb{D}_n$, $\mathbb{E}_6$, $\mathbb{E}_7$, $\mathbb{G}_2$ and $\mathbb{I}_2(p)$ presented in Table B. Therefore by Propositions \ref{2.12}, \ref{2.1}, \ref{3.4}, \ref{3.5}, \ref{3.1} and Theorem \ref{3.7}, the proof complete.\\
$(\Leftarrow)$. By Propositions \ref{2.12}, \ref{2.1}, \ref{3.4}, \ref{3.5}, \ref{3.1} and Theorem \ref{3.7}, for each $1 \leq i \leq s$, every right $\Lambda_i$-module is a direct sum of cyclic modules. Therefore $\Lambda$ is a right K\"othe ring.
\end{proof}

\section{A characterization of right K\"othe rings with radical square zero}
Let $U$ and $T$ be two rings, $M$ be a $T$-$U$-bimodule and $R = \tiny{\begin{pmatrix}
T & M\\
0& U
\end{pmatrix}}$. Let $\mathscr{D}_{R}$ be the category whose objects are triples $(X, Y, f)$, where $X$ is a right $T$-module, $Y$ is a right $U$-module and $f \in {\rm Hom}_U(X\otimes_TM, Y)$. If $\alpha \in {\rm Hom}_{\mathscr{D}_{R}}((X, Y, f), (X^{'}, Y^{'}, f^{'}))$, then $\alpha = (\alpha_1, \alpha_2)$, where $\alpha_1 \in {\rm Hom}_T(X, X^{'})$ and $\alpha_2 \in {\rm Hom}_U(Y, Y^{'})$ such that $\alpha_2f = f^{'}(\alpha_1 \otimes id_M)$. The  functor $\mathrm{G}:\mathscr{D}_{R} \rightarrow $Mod-$R$ is defined in \cite{Fossum} (see also \cite{3}) as follows. Let $(X, Y, f)$ be an object in the category $\mathscr{D}_{R}$. For $(x,y) \in X \oplus Y$ and $\tiny{\begin{pmatrix}
t & m\\
0& u
\end{pmatrix}} \in R$,  define $$(x,y)\tiny{\begin{pmatrix}
t & m\\
0& u
\end{pmatrix}} = (xt, f(x \otimes m ) + yu).$$
It is easy to see that $X \oplus Y$ is a right $R$-module. We define $\mathrm{G}((A, B, f))$ to be $A \oplus B$. Let $\alpha = (\alpha_1, \alpha_2) \in {\rm Hom}_{\mathscr{D}_{R}}((X, Y, f), (X^{'}, Y^{'}, f^{'}))$. We set $\mathrm{G}(\alpha)= \alpha_1 \oplus \alpha_2$.  The reader may easily verify that $\mathrm{G}(\alpha)$ is an $R$-homomorphism. It is well known that the functor $\mathrm{G}$ is an equivalence. From now on we will identify the categories $\mathscr{D}_{R}$ and Mod-$R$.\\

 Let $\Lambda$ be a basic artinian ring with radical square zero and $R = \tiny{\begin{pmatrix}
\Lambda/J & J\\
0& \Lambda/J
\end{pmatrix}}$. Let $\mathscr{A}$ denote the full subcategory of Mod-$R$ whose objects are $(X, Y, f)$, where $X$ and $Y$ are two right $\Lambda/J$-modules and $f \in {\rm Hom}_{\Lambda/J}(X \otimes_{\Lambda/J}J, Y)$ is an epimorphism. Then we have the natural functor $\mathrm{H}:$ Mod-$\Lambda \rightarrow \mathscr{A}$ which is defined in \cite{Fossum} (see also \cite{3}) as follows. Let $M$ be a right $\Lambda$-module. Then $\mathrm{H}(M) = (M/MJ, MJ, f_M)$, where $f_M: M/MJ \otimes_{\Lambda/J} J \rightarrow MJ$ is induced from the multiplication map $M \otimes_{\Lambda} J \rightarrow MJ$. It is well known that the functor $\mathrm{H}$ is full and dense and $M \in$ Mod-$\Lambda$ is indecomposable  if and only if $\mathrm{H}(M)$ in $\mathscr{A}$ is indecomposable.

\begin{Pro}\label{m2}
Let $\Lambda$ be a basic artinian ring with radical square zero. Then $\Lambda$ is a right K\"othe ring if and only if $R=\tiny{\begin{pmatrix}
{\Lambda} / J & J\\
0 & {\Lambda} / J
\end{pmatrix}}$ is a right K\"othe ring.
\end{Pro}
\begin{proof}
$(\Rightarrow)$. Assume that ${\Lambda}$ is a right K\"othe ring. Since $\Lambda$ is a basic artinian ring, by \cite[Proposition 1.8]{Goodearl}, $R$ is a basic artinian ring. Let $(A, B, f)$ be a finitely generated indecomposable right $R$-module. If $(A, B, f) \not\in \mathscr{A}$, then $f$ is not an epimorphism. Since $B$ is semisimple, $B = B^{'} \oplus {\rm Im}f$ for some $\Lambda/J$-submodule $B^{'}$ of $B$. Thus $(A, B, f) \cong (A, {\rm Im}f, f) \oplus (0, B^{'},f)$. Since $(A, B, f)$ is indecomposable, $(A, B, f) \cong (0, B^{'},f)$ and hence $B^{'}$ is a simple right $\Lambda$-module. Therefore by \cite[Corollary 2.2]{Ha}, $(A, B, f)$ is multiplicity-free top. Now assume that $(A, B, f) \in \mathscr{A}$. Then there exists an indecomposable right ${\Lambda}$-module $M$ such that $\mathrm{H}(M) \cong (A, B, f)$. It follows that $(A, B, f) \cong (M/MJ, MJ, f_M)$,  where $f_M: M/MJ \otimes_{\Lambda/J} J \rightarrow MJ$ is induced from the multiplication map $M \otimes_{\Lambda} J \rightarrow MJ$. Since $\Lambda$ is right K\"othe, by Proposition \ref{r3}, $M/MJ= S_1 \oplus \cdots \oplus S_t$, where $t \in {\Bbb{N}}$, each $S_i$ is a simple right $\Lambda$-module and  $S_i\ncong S_j$ for each $i \neq j$. Therefore by \cite[Corollary 2.2]{Ha}, ${\rm top}((A, B, f)) =  (S_1, 0, 0) \oplus \cdots \oplus (S_t, 0, 0)$, where for $i \neq j$, $(S_i, 0, 0) \ncong (S_j, 0, 0)$. This prove that every finitely generated indecomposable right $R$-module is multiplicity-free top. Therefore by Proposition \ref{r3}, $R$ is a right K\"othe ring. \\
$(\Leftarrow)$. Let $N$ be a finitely generated indecomposable right ${\Lambda}$-module. Then by \cite[Exercise 1C]{Goodearl}, $\mathrm{H}(N)= (N/NJ, NJ, f_N)$ is a finitely generated indecomposable right $R$-module. Since by \cite[Corollary 2.2]{Ha}, ${\rm top}(\mathrm{H}(N))= (S_1, 0, 0) \oplus \cdots \oplus (S_t, 0, 0)$, where  $t \in {\Bbb{N}}$ and $N/NJ= S_1 \oplus \cdots \oplus S_t$,  by Proposition \ref{r3}, for each $i \neq j$, $S_i\ncong S_j$. Therefore ${\Lambda}$ is a right K\"othe ring.
\end{proof}

Let $\Lambda$ be a basic artinian ring with radical square zero and ${\mathscr{M}}_{\Lambda}= \left( F_i,{_iM_j} \right)_{i, j \in I} $ be the species of $\Lambda$. Note that the species ${\mathscr{M}}_{\Lambda}$ doesn't necessarily have the property ``${{_iM_j}\neq 0}$ implies that ${_jM_i} = 0$".  Let $(\Gamma_{\Lambda}, \mathbf{d})$ be the valued quiver of $\Lambda$. We recall from \cite{DS} that a {\it separated quiver} of $\Lambda$ is the valued graph $(\Gamma^s_{{\Lambda}}, \mathbf{d}^s)$ with the vertex set $\lbrace (i,l)~:~i \in I, l=0,1 \rbrace$ and the arrows \begin{center}
	\definecolor{qqqqff}{rgb}{0.,0.,1.}
	\begin{tikzpicture}[line cap=round,line join=round,>=triangle 45,x=1.0cm,y=1.0cm]
	\draw (2.,9.)-- (3.,9.);
		\begin{scriptsize}
	\draw[color=black] (.8,9.) node {};
	\draw [fill=white] (2.,9.) circle (1.5pt);
	\draw[color=black] (2.,8.7) node {$(i,0)$};
	\draw[color=black] (2.5,9.3) node {$(d_{ij},d_{ji})$};
	\draw[color=black] (2.85,9.) node {$>$};
	\draw [fill=white] (3.,9.) circle (1.5pt);
	\draw[color=black] (3.,8.7) node {$(j,1)$};
	
	\end{scriptsize}
	\end{tikzpicture}
\end{center}
 precisely when ${_iM_j}\neq 0$. If $d_{ij} = d_{ji}=1$, we write simply
 \begin{center}
	\definecolor{qqqqff}{rgb}{0.,0.,1.}
	\begin{tikzpicture}[line cap=round,line join=round,>=triangle 45,x=1.0cm,y=1.0cm]
	\draw (2.,9.)-- (3.,9.);
		\begin{scriptsize}
	\draw[color=black] (.8,9.) node {};
	\draw [fill=white] (2.,9.) circle (1.5pt);
	\draw[color=black] (2.,8.7) node {$(i,0)$};
	\draw[color=black] (2.85,9.) node {$>$};
	\draw [fill=white] (3.,9.) circle (1.5pt);
	\draw[color=black] (3.,8.7) node {$(j,1)$};
	
	\end{scriptsize}
	\end{tikzpicture}
\end{center}
It is easy to see that the separated quiver of $\Lambda$ coincides with the valued quiver of $R$ (see \cite{DS}). Moreover, it is well known that $R = \tiny{\begin{pmatrix}
\Lambda/J & J\\
0& \Lambda/J
\end{pmatrix}}$ is a basic hereditary artinian ring (see \cite{DS} and also \cite{Goodearl}). \\

 We conclude this section with the following result which is a characterization of basic right K\"othe rings with  radical square zero in terms of their separated quivers.

\begin{The}\label{m3}
Let $\Lambda$ be a basic ring with radical square zero and ${\mathscr{M}}_{\Lambda}= (F_i, {_iM_j})_{i,j \in I}$ be the species of $\Lambda$ (note that ${\mathscr{M}}_{\Lambda}$ doesn't necessarily have the property ``${{_iM_j}\neq 0}$ implies that ${_jM_i} = 0$").  Then the following statements are equivalent:
\begin{itemize}
\item[$(i)$] $\Lambda$ is a right K\"othe ring;
\item[$(ii)$] $R = \tiny{\begin{pmatrix}
\Lambda/J & J\\
0& \Lambda/J
\end{pmatrix}}$ is an artinian ring such that the Coxeter valued quiver $(\mathscr{C}_{R}, \mathbf{m})$ of $R$ is a finite disjoint union of the Coxeter valued quivers presented in Theorem \ref{m1};
\item[$(iii)$] $\Lambda$ is a representation-finite ring and the separated quiver $(\Gamma^s_{\Lambda}, \mathbf{d})$ of $\Lambda$ is a finite disjoint union of the following valued Dynkin quivers:
\begin{itemize}
\item[$(1)$] $\mathbb{A}_n$ with any orientation;
\item[$(2)$] $\mathbb{B}_n$ with the orientation $$ 
	\definecolor{qqqqff}{rgb}{0.,0.,1.}
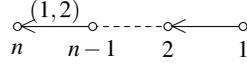
	\begin{tikzpicture}[line cap=round,line join=round,>=triangle 45,x=1.0cm,y=1.0cm]
	\draw (2.,9.)-- (3.,9.);
	\draw [dash pattern=on 2pt off 2pt](3.,9.)--(4.,9.);
	\draw (4.,9.)-- (5.,9.);
		\begin{scriptsize}
	\draw[color=black] (.8,9.) node {};
	\draw [fill=white] (2.,9.) circle (1.5pt);
	\draw[color=black] (2.,8.7) node {$n$};
	\draw[color=black] (2.5,9.2) node {$(1,2)$};
	\draw[color=black] (2.14,9.) node {$<$};
	\draw [fill=white] (3.,9.) circle (1.5pt);
	\draw[color=black] (3.,8.7) node {$n-1$};
	\draw [fill=white] (4.,9.) circle (1.5pt);
	\draw[color=black] (4.,8.7) node {$2$};
	\draw[color=black] (4.14,9.) node {$<$};
	\draw [fill=white] (5.,9.) circle (1.5pt);
	\draw[color=black] (5.,8.7) node {$1$};
	
	\end{scriptsize}
	\end{tikzpicture}
$$
where there exist precisely $4$ pairwise non-isomorphic finitely generated indecomposable right $\tiny{\begin{pmatrix}
F_{n-1} & {_{n-1}M_n}\\
0 & F_n
\end{pmatrix}}$-modules;
\item[$(3)$] $\mathbb{C}_n$ with the orientation $$ 
	\definecolor{qqqqff}{rgb}{0.,0.,1.}
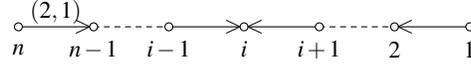
	\begin{tikzpicture}[line cap=round,line join=round,>=triangle 45,x=1.0cm,y=1.0cm]
	\draw (2.,9.)-- (3.,9.);
	\draw [dash pattern=on 2pt off 2pt](3.,9.)--(4.,9.);
	\draw (4.,9.)-- (5.,9.);
	\draw (5.,9.)-- (6.,9.);
	\draw [dash pattern=on 2pt off 2pt](6.,9.)--(7.,9.);
	\draw (7.,9.)-- (8.,9.);
		\begin{scriptsize}
	\draw [fill=white] (2.,9.) circle (1.5pt);
	\draw[color=black] (2.,8.7) node {$n$};
	\draw[color=black] (2.5,9.2) node {$(2,1)$};
	\draw[color=black] (2.86,9.) node {$>$};
	\draw [fill=white] (3.,9.) circle (1.5pt);
	\draw[color=black] (3.,8.7) node {$n-1$};
	\draw [fill=white] (4.,9.) circle (1.5pt);
	\draw[color=black] (4.,8.7) node {$i-1$};
	\draw[color=black] (4.86,9.) node {$>$};
	\draw [fill=white] (5.,9.) circle (1.5pt);
	\draw[color=black] (5.,8.7) node {$i$};
	\draw [fill=white] (6.,9.) circle (1.5pt);
	\draw[color=black] (6.,8.7) node {$i+1$};
	\draw[color=black] (5.14,9.) node {$<$};
	\draw [fill=white] (7.,9.) circle (1.5pt);
	\draw[color=black] (7.,8.7) node {$2$};
	\draw[color=black] (7.14,9.) node {$<$};
	\draw [fill=white] (8.,9.) circle (1.5pt);
	\draw[color=black] (8.,8.7) node {$1$};
	
	\end{scriptsize}
	\end{tikzpicture}
$$
where $1 \leq i \leq n-1$ and there exist precisely $4$ pairwise non-isomorphic finitely generated indecomposable right $\tiny{\begin{pmatrix}
F_n & {_nM_{n-1}}\\
0 & F_{n-1}
\end{pmatrix}}$-modules;
\item[$(4)$] $\mathbb{D}_n$ with the following conditions:
 \begin{center}
 \definecolor{qqqqff}{rgb}{0.,0.,1.}
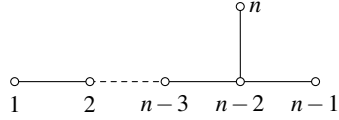
	\begin{tikzpicture}[line cap=round,line join=round,>=triangle 45,x=1.0cm,y=1.0cm]
	\draw (2.,9.)-- (3.,9.);
	\draw [dash pattern=on 2pt off 2pt](3.,9.)--(4.,9.);
	\draw (4.,9.)-- (5.,9.);
	\draw (5.,9.)-- (6.,9.);
	\draw (5.,9.)-- (6.,10.9999999999987,2.6);
		\begin{scriptsize}
	\draw [fill=white] (2.,9.) circle (1.5pt);
	\draw[color=black] (2.,8.7) node {$1$};
	\draw [fill=white] (3.,9.) circle (1.5pt);
	\draw[color=black] (3.,8.7) node {$2$};
	\draw [fill=white] (4.,9.) circle (1.5pt);
	\draw[color=black] (4.,8.7) node {$n-3$};
	\draw [fill=white] (5.,9.) circle (1.5pt);
	\draw[color=black] (5.,8.7) node {$n-2$};
	\draw [fill=white] (6.,9.) circle (1.5pt);
	\draw[color=black] (6.,8.7) node {$n-1$};
	\draw [fill=white] (6.,10.9999999999987,2.6) circle (1.5pt);
	\draw[color=black] (6.2,10.9999999999987,2.6) node {$n$};
	\end{scriptsize}
	\end{tikzpicture}
 \end{center}
 \begin{itemize}
\item[$(a)$] $|{(n-2)}^{+}| \leq 2$;
\item[$(b)$] For each $i\leq n-3$, there exists at most one arrow with the source $i$;
\end{itemize}
\item[$(5)$] $\mathbb{E}_6$ with the following conditions:
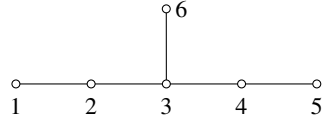
\begin{center}
 \definecolor{qqqqff}{rgb}{0.,0.,1.}
	\begin{tikzpicture}[line cap=round,line join=round,>=triangle 45,x=1.0cm,y=1.0cm]
	\draw (2.,9.)-- (3.,9.);
	\draw (3.,9.)-- (4.,9.);
	\draw (4.,9.)-- (5.,9.);
	\draw (4.,9.)-- (5.,10.9999999999987,2.6);
	\draw (5.,9.)-- (6.,9.);
		\begin{scriptsize}
	\draw [fill=white] (2.,9.) circle (1.5pt);
	\draw[color=black] (2.,8.7) node {$1$};
	\draw [fill=white] (3.,9.) circle (1.5pt);
	\draw[color=black] (3.,8.7) node {$2$};
	\draw [fill=white] (4.,9.) circle (1.5pt);
	\draw[color=black] (4.,8.7) node {$3$};
	\draw [fill=white] (5.,9.) circle (1.5pt);
	\draw[color=black] (5.,8.7) node {$4$};
	\draw [fill=white] (5.,10.9999999999987,2.6) circle (1.5pt);
	\draw[color=black] (5.2,10.9999999999987,2.6) node {$6$};
	\draw [fill=white] (6.,9.) circle (1.5pt);
	\draw[color=black] (6.,8.7) node {$5$};
	\end{scriptsize}
	\end{tikzpicture}
 \end{center}
\begin{itemize}
\item[$(a)$] $1 \leq |3^{+}| \leq 2$, $|4^{+}| \leq 1$ and $|2^{+}| \leq 1$;
\item[$(b)$] For each $y \in 3^{-}$, there exists at least one arrow with the target $y$;
\end{itemize}
\item[$(6)$] $\mathbb{E}_7$ with the orientation
$$
	\definecolor{qqqqff}{rgb}{0.,0.,1.}
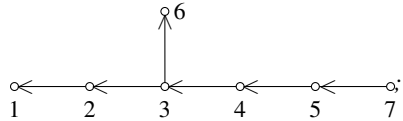
	\begin{tikzpicture}[line cap=round,line join=round,>=triangle 45,x=1.0cm,y=1.0cm]
	\draw (2.,9.)-- (3.,9.);
	\draw (3.,9.)-- (4.,9.);
	\draw (4.,9.)-- (5.,9.);
	\draw (4.,9.)-- (5.,10.9999999999987,2.6);
	\draw (5.,9.)-- (6.,9.);
	\draw (6.,9.)-- (7.,9.);
		\begin{scriptsize}
	\draw[color=black] (.8,9.) node {};
	\draw[color=black] (.8,9.) node {};
	\draw [fill=white] (2.,9.) circle (1.5pt);
	\draw[color=black] (2.,8.7) node {$1$};
	\draw[color=black] (2.14,9.) node {$<$};
	\draw [fill=white] (3.,9.) circle (1.5pt);
	\draw[color=black] (3.,8.7) node {$2$};
	\draw[color=black] (3.14,9.) node {$<$};
	\draw [fill=white] (4.,9.) circle (1.5pt);
	\draw[color=black] (4.,8.7) node {$3$};
	\draw[color=black] (4.14,9.) node {$<$};
	\draw [fill=white] (5.,9.) circle (1.5pt);
	\draw[color=black] (5.,8.7) node {$4$};
    \draw[color=black] (5.14,9.) node {$<$};
	\draw [fill=white] (5.,10.9999999999987,2.6) circle (1.5pt);
	\draw[color=black] (5.2,10.9999999999987,2.6) node {$6$};
	\draw[color=black] (4.,9.85) node {$\wedge$};
	\draw [fill=white] (6.,9.) circle (1.5pt);
	\draw[color=black] (6.,8.7) node {$5$};
	\draw [fill=white] (7.,9.) circle (1.5pt);
	\draw[color=black] (7.,8.7) node {$7$};
    \draw[color=black] (6.17,9.) node {$<$};
	\draw [fill=black] (7.1,9.) node {;};
	\end{scriptsize}
	\end{tikzpicture}$$
	
\item[$(7)$] The valued quiver
$$
	\definecolor{qqqqff}{rgb}{0.,0.,1.}
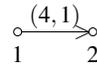
	\begin{tikzpicture}[line cap=round,line join=round,>=triangle 45,x=1.0cm,y=1.0cm]
	\draw (2.,9.)-- (3.,9.);
		\begin{scriptsize}
	\draw[color=black] (.8,9.) node {};
	\draw [fill=white] (2.,9.) circle (1.5pt);
	\draw[color=black] (2.,8.7) node {$1$};
	\draw[color=black] (2.5,9.2) node {$(4,1)$};
	\draw[color=black] (2.87,9.) node {$>$};
	\draw [fill=white] (3.,9.) circle (1.5pt);
	\draw[color=black] (3.,8.7) node {$2$};
	
	\end{scriptsize}
	\end{tikzpicture}
$$
where there exist precisely $6$ pairwise non-isomorphic finitely generated indecomposable right $\tiny{\begin{pmatrix}
F_1 & {_1M_2}\\
0 & F_2
\end{pmatrix}}$-modules;

\item[$(8)$] The valued quiver
$$
	\definecolor{qqqqff}{rgb}{0.,0.,1.}
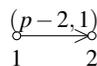
	\begin{tikzpicture}[line cap=round,line join=round,>=triangle 45,x=1.0cm,y=1.0cm]
	\draw (2.,9.)-- (3.,9.);
		\begin{scriptsize}
	\draw[color=black] (.8,9.) node {};
	\draw [fill=white] (2.,9.) circle (1.5pt);
	\draw[color=black] (2.,8.7) node {$1$};
	\draw[color=black] (2.5,9.2) node {$(p-2,1)$};
	\draw[color=black] (2.87,9.) node {$>$};
	\draw [fill=white] (3.,9.) circle (1.5pt);
	\draw[color=black] (3.,8.7) node {$2$};
	
	\end{scriptsize}
	\end{tikzpicture}
$$
where $p=5$ or $7 \leq p < \infty$ and there exist precisely $p$ pairwise non-isomorphic finitely generated indecomposable right $\tiny{\begin{pmatrix}
F_1 & {_1M_2}\\
0 & F_2
\end{pmatrix}}$-modules.
\end{itemize}
\end{itemize}
\end{The}
\begin{proof}
$(i)\Leftrightarrow (ii)$. It follows from Theorems \ref{m1} and \ref{m2}.\\
$(ii) \Rightarrow (iii)$. It follows from Proposition \ref{r2} and \cite[Lemma 3.1]{Sim5}. \\
$(iii) \Rightarrow (i)$. Since $\Lambda$ is a basic artinian ring with radical square zero, $R$ is a basic hereditary artinian ring. Thus $R=R_1 \oplus \cdots \oplus R_m$, where each $R_i$ is an indecomposable basic hereditary artinian ring. Since $R$ is a right K\"othe ring if and only if each $R_i$ is a right K\"othe ring and the valued quiver of $R$ is a disjoint union of the valued quivers of $R_i$,  without loss of generality, we can assume that $R$ is an indecomposable basic hereditary artinian ring. Let the separated quiver $(\Gamma^s_{{\Lambda}}, \mathbf{d})$ of $\Lambda$ be one of the valued Dynkin quivers $\mathbb{A}_n$, $\mathbb{B}_n$, $\mathbb{C}_n$,  $\mathbb{D}_n$, $\mathbb{E}_6$ and $\mathbb{E}_7$ presented in $(iii)$. Since $R$ is a representation-finite ring, by Proposition \ref{r2}, the underlying Coxeter valued graph of the Coxeter valued quiver $(\mathscr{C}_{\mathrm{R}}, \mathbf{m})$ of $\mathrm{R}$ is one of the Coxeter diagrams $\mathbb{A}_n$, $\mathbb{B}_n$, $\mathbb{D}_n$, $\mathbb{E}_6$ and $\mathbb{E}_7$ presented in Table B. Hence by assumption and Theorem \ref{m1},  $R$ is a right K\"othe ring. Now assume that the separated quiver $(\Gamma^s_{{\Lambda}}, \mathbf{d})$ of $\Lambda$ is the quiver
$$
	\definecolor{qqqqff}{rgb}{0.,0.,1.}
	\begin{tikzpicture}[line cap=round,line join=round,>=triangle 45,x=1.0cm,y=1.0cm]
	\draw (2.,9.)-- (3.,9.);
		\begin{scriptsize}
	\draw[color=black] (.8,9.) node {};
	\draw [fill=white] (2.,9.) circle (1.5pt);
	\draw[color=black] (2.,8.7) node {$1$};
	\draw[color=black] (2.5,9.2) node {$(4,1)$};
	\draw[color=black] (2.87,9.) node {$>$};
	\draw [fill=white] (3.,9.) circle (1.5pt);
	\draw[color=black] (3.,8.7) node {$2$};
	
	\end{scriptsize}
	\end{tikzpicture}
$$
 and there exist precisely $6$ pairwise non-isomorphic finitely generated indecomposable right $\tiny{\begin{pmatrix}
F_1 & {_1M_2}\\
0 & F_2
\end{pmatrix}}$-modules. Then by \cite[Corollary 3.5]{Sim5}, $d_6({_1M_2})$ is a dimension sequence and hence by \cite[Proposition $2^{'}$]{DRS}, $d_6({_1M_2}) = (4,1,2,2,2,1)$. Thus by Theorem \ref{m1}, $R$ is a right K\"othe ring.  If the separated quiver $(\Gamma^s_{{\Lambda}}, \mathbf{d})$ of $\Lambda$ is the quiver presented in (8), by the similar argument we can see that $R$ is a right K\"othe ring.  Therefore by Proposition \ref{m2}, $\Lambda$ is a right K\"othe ring.
\end{proof}

\section*{acknowledgements} The authors would like to thank Daniel Simson for his comments on the earlier version of this paper. Also, we are indebted to the anonymous reviewers’ insightful comments, suggestions and feedback on this paper and their time and effort spent on it. The research of the first author was in part supported by a grant from IPM. Also, the research of the second author was in part supported by a grant from IPM (No. 98170412).

\end{document}